\documentclass[12pt]{amsart}
\usepackage{amssymb}
\usepackage{amscd}
\usepackage[all]{xy}
\usepackage[normalem]{ulem}

\DeclareMathAlphabet{\matheur}{U}{eur}{m}{n}

\def\assign{\mathrel{\raise.095ex\hbox{\text{\rm :}}\mkern-4.2mu=}}
\def\assigN{\mathrel{=\mkern-4.2mu\raise.095ex\hbox{\text{\rm :}}}}
\newcommand{\Fq}{{\mathbb F}_q}
\newcommand{\ok}{\overline{k}}

\newcommand{\C}{C_\infty}
\newcommand{\Ct}{{\C}\{\tau\}}
\newcommand{\bsd}{\boldsymbol{\delta}}
\newcommand{\bsl}{\boldsymbol{\lambda}}

\newcommand{\bss}{\boldsymbol{\sigma}}
\newcommand{\bsU}{\boldsymbol{U}}
\newcommand{\bsV}{\boldsymbol{V}}

\newcommand{\Ga}{{\mathbb G}_a}

\newcommand{\ZZ}{\mathbb Z}

\newcommand{\bm}{{\mathbf m}}
\newcommand{\bn}{{\mathbf n}}
\newcommand{\bx}{{\mathbf x}}
\newcommand{\bu}{{\mathbf u}}
\newcommand{\bz}{{\mathbf z}}

\newcommand{\bU}{{\mathbf U}}

\DeclareMathOperator{\Lie}{Lie}
\DeclareMathOperator{\Ker}{Ker}

\DeclareMathOperator{\Mat}{Mat}

\DeclareMathOperator{\End}{End}
\DeclareMathOperator{\res}{res}
\DeclareMathOperator{\Exp}{Exp}
\DeclareMathOperator{\Gal}{Gal}
\DeclareMathOperator{\Ext}{Ext}
\DeclareMathOperator{\DR}{DR}

\DeclareMathOperator{\Emb}{Emb}

\newtheorem{theorem}{Theorem}[subsection]
\newtheorem{lemma}[theorem]{Lemma}
\newtheorem{proposition}[theorem]{Proposition}
\newtheorem{corollary}[theorem]{Corollary}

\theoremstyle{definition}
\newtheorem*{definition}{Definition}

\newtheorem{remark_num}[theorem]{Remark}
\newtheorem*{remark}{Remark}

\numberwithin{equation}{section}

\newcommand{\abs}[1]{\lvert#1\rvert}

\newcommand{\ibar}{\overline{\infty}}
\newcommand{\tpi}{\widetilde{\pi}}
\newcommand{\PP}{\mathbb{P}}
\newcommand{\bX}{\mathbf{X}}
\newcommand{\fA}{\mathfrak{A}}
\newcommand{\fB}{\mathfrak{B}}
\newcommand{\fU}{\mathfrak{U}}
\newcommand{\fX}{\mathfrak{X}}
\newcommand{\cI}{\mathcal{I}}

\newcommand{\cM}{\mathcal{M}}
\newcommand{\cN}{\mathcal{N}}
\newcommand{\cO}{\mathcal{O}}
\newcommand{\cR}{\mathcal{R}}
\newcommand{\cS}{\mathcal{S}}
\newcommand{\cU}{\mathcal{U}}
\newcommand{\bd}{\boldsymbol{\delta}}

\newcommand{\eA}{\matheur{A}}
\newcommand{\eB}{\matheur{B}}
\newcommand{\eC}{\matheur{C}}
\newcommand{\eF}{\matheur{F}}
\newcommand{\eG}{\matheur{G}}
\newcommand{\eS}{\matheur{S}}
\newcommand{\eK}{\matheur{K}}
\newcommand{\eL}{\matheur{L}}
\newcommand{\eM}{\matheur{M}}
\newcommand{\ea}{a}
\newcommand{\eb}{b}
\newcommand{\ec}{c}
\newcommand{\ef}{f}
\newcommand{\eg}{g}
\newcommand{\ek}{\matheur{k}}
\newcommand{\eum}{m}
\newcommand{\es}{s}
\newcommand{\eu}{u}

\newcommand{\ey}{y}

\DeclareMathOperator{\Der}{Der}
\DeclareMathOperator{\Hom}{Hom}
\DeclareMathOperator{\Spec}{Spec}
\DeclareMathOperator{\rank}{rank}
\DeclareMathOperator{\Frob}{Frob}
\DeclareMathOperator{\RLie}{{\cR}Lie}
\DeclareMathOperator{\RExp}{{\cR}Exp}
\DeclareMathOperator{\dv}{div}

\newcommand{\Fd}{{F_{\bd}}}
\newcommand{\RF}{{\cR\Fd}}
\newcommand{\Gmm}[1]{\Gamma\! \left( #1 \right)}
\newcommand{\Lt}{L\{ \tau \}}
\newcommand{\power}[2]{#1[\![ #2 ]\!]}
\newcommand{\laurent}[2]{#1(\!( #2 )\!)}

\begin{document}

%%%%%%%%%%%%%%%%%%%%%%%%%%%%%%%%%%%%%%%%%
%               Title, Etc.             %
%%%%%%%%%%%%%%%%%%%%%%%%%%%%%%%%%%%%%%%%%

\title[Linear independence of gamma values]{Linear independence of
  gamma values in positive~characteristic}
\author{W.~Dale Brownawell}
\address{Department of Mathematics\\
  Penn State University\\
  University Park, PA 16802} \email{wdb@math.psu.edu} \author{Matthew
  A.~Papanikolas}
\address{Department of Mathematics\\
  Penn State University\\
  University Park, PA 16802} \email{map@math.psu.edu}
%\date{June 5, 2001} %
\thanks{The research of the first author was partially supported by an
  NSF grant and the second author by K.~Ono's Sloan Research
  Fellowship.}

\begin{abstract}
  We investigate the arithmetic nature of special values of Thakur's
  function field Gamma function at rational points.  Our main result
  is that all linear dependence relations over the field of algebraic
  functions are consequences of the Anderson-Deligne-Thakur bracket
  relations.
\end{abstract}


%%%%%%%%%%%%%%%%%%%%%%%%%%%%%%%%%%%%%%%%%
%               Document Text           %
%%%%%%%%%%%%%%%%%%%%%%%%%%%%%%%%%%%%%%%%%

\maketitle

\tableofcontents

\section{Introduction and Statement of Results} \label{S:int}

\subsection{Transcendence of Gamma Values}

Let $\Fq$ be the field of $q$ elements, where $q$ is a power of a
prime $p$.  Let $A \assign \Fq[\theta]$, $k \assign \Fq(\theta)$ for a
variable $\theta$.  Let $\C$ be the completion of the algebraic
closure of the completion $\Fq(\!(1/\theta)\!)$ with respect to the
non-archimedean absolute value on $k$ for which $\abs{\theta} = q$.
Let $A_{+} \assign \{ a \in A : a\ \text{\rm is monic}\}$ be the
``positive integers'' of $A$.

In this setting, D.~Thakur studied a Gamma function
\[
\Gamma (z) \assign \frac{1}{z} \prod_{n \in A_{+}} \left( 1 +
  \frac{z}{n} \right)^{-1},
\]
which is meromorphic on $\C$ with poles at zero and the ``negative''
integers $-n \in -A_{+}$.  One recognizes immediately a strong analogy
with the classical Euler Gamma function
\[
\Gamma(z) = \frac{e^{-\gamma z}}{z} \prod_{n = 1}^\infty \left( 1 +
  \frac{z}{n}\right)^{-1}e^{z/n}.
\]
Thakur's function shares other striking features with the Euler Gamma
function, such as various natural functional equations, and rational
(or infinite) values at the integers.  This Gamma function is a
one-variable specialization of the two-variable Gamma function defined
by D.\ Goss \cite{Goss88}.

After isolated results by Thakur in \cite{Th91}, S.~Sinha established the
first transcendence results for general classes of values of the Gamma
function.

\begin{theorem}[Sinha {\cite{SS95}; \cite{SS97d}, {\S}6.2}]  Let $a$, $f
  \in A_{+}$, $b \in A$, with $(a,f) = 1$, $\deg a < \deg f$. Then $
  \Gamma(\frac{a}{f} + b)$ is transcendental over $k$.
\end{theorem}

When $q=2$, this result had been observed by Thakur in \cite{Th91}.
Our goal is to extend Sinha's result, treating several values at once
and evaluating $\Gamma$ at more general arguments.  There are however
some natural dependencies among the values.

\subsection{Dependence of Gamma Values} \label{SS:brk}

In \cite{Th91} Thakur established algebraic relations on Gamma values
in this setting which are analogues of well-known relations of
Deligne-Koblitz-Ogus in \cite{K-O} for the classical Gamma function.
These relations express certain ratios of Gamma values at rational
arguments as algebraic multiples of powers of the Carlitz period,
\[
\tpi = \theta \sqrt[q-1]{-\theta} \prod_{i=1}^{\infty} \left( 1 -
  \theta^{1-q^i} \right)^{-1} \quad \in \laurent{\Fq}{1/\theta} \cdot
  \!\sqrt[q-1]{-\theta},
\]
where $\sqrt[q-1]{-\theta}$ is some fixed $(q-1)$-st root of
$-\theta$.  The quantity $\tpi$ arises as the fundamental period
of the Carlitz module (see \cite{Goss98}), whence its name.

The Anderson-Deligne-Thakur relations on Gamma values (see
\cite{SS97a} or \cite{Th91}), derived from the formalism of
\emph{brackets}, may be expressed as follows.  Fix $f \in A_+$ and
define
\begin{gather*}
\cN_f \assign \{ a \in A : a \not\equiv 0 \bmod{f} \}, \qquad
\cU_f \assign \{ u \in A : (u,f) = 1 \}, \\
\cM_f \assign \{ a \in \cN_f : a \equiv m \bmod{f}, m \text{\rm\ monic},
  \deg m < \deg f \}.
\end{gather*}
Then for any $u \in \cU_f$, multiplication by $u$ gives an injection
of the set $\cN_f$ into itself.  This action of $\cU_f$ on $\cN_f$
induces an action of $\cU_f$ on $\oplus_{\cN_f} \ZZ$, which we denote
$(u,\bm) \mapsto u \ast \bm$.  

Our restatement of the bracket relations involves the sum of those
coordinates of $\bm \assign (m_a)_{a \in {\cN_f}} \in \oplus_{\cN_f}
\ZZ$ having monic indices modulo $f$:
\[
\Sigma_{+}(\bm) \assign \sum_{a \in \cM_f} m_a.
\]
The following theorem gives the known (and, it is believed, all)
$\ok$-algebraic relations on Gamma values at rational points.

\begin{definition}
For non-zero $x$, $y \in \C$, we write $x \thicksim y$ if $x/y \in \ok$.
\end{definition}

\begin{theorem}[{Thakur \cite{Th91}, {\S}7.8}] \label{T:brk}
  Let $f \in A_{+}$.  Suppose that $\bm = (m_a) \in \oplus_{\cN_f}
  \ZZ$ and that $\Sigma_+(\bm) = \Sigma_+ (u\ast \bm)$, for any
  {\text{\rm (}}\!and thus every{\text{\rm )}} choice of representatives
  $u$ of elements from $(A/f)^{\times}$.  Then
\[
\prod_{a \in \cN_f} \Gamma \left(\frac{a}{f} \right)^{m_a} 
\thicksim \quad \tpi^{\Sigma_+(\bm)}.
\]
\end{theorem}

An equivalent version of this theorem was first formulated by
G.\ Anderson \cite{Th91}, {\S}7, and then proven by Thakur.  Moreover
Sinha obtained information on the quotient of the two sides in the
above theorem in \cite{SS97a}, Thm.~2.2.4 and Rmk.~3.3.6.

The Gamma function satisfies several functional equations, directly
analogous to the functional equations of the classical Gamma function,
which imply the following algebraic relations on special values
\cite{Th91}:  For all $r \in k\smallsetminus A$, $a \in A$, $g$ in
$A_+$, $\deg g = d$,
\begin{align}
\Gmm{r + a} &\thicksim \Gmm{r}, \label{E:fneq1} \\
\prod_{\theta \in \Fq^{\times}} \Gmm{\theta r} &\thicksim \tpi, \\
\prod_{\alpha \in A/(g)} \Gamma \Bigl( \frac{r +
\alpha}{g} \Bigr) &\thicksim \tpi^{\frac{q^d - 1}{q-1}} \Gmm{r}.
\label{E:fneq3}
\end{align}
The polynomial relations \eqref{E:fneq1}--\eqref{E:fneq3} are special
cases of the bracket relations.  See the theorem of Thakur in
\cite{SS95}, Thm.~VII.1, or \cite{SS97a}, Thm~2.1.3, for
the extent to which the bracket relations are more general than
\eqref{E:fneq1}--\eqref{E:fneq3}.

\begin{definition}
  We adopt the notation $\Gamma(a/f) \approx \Gamma(b/f)$ to indicate
  that the relation $\Gamma(a/f) \thicksim \Gamma(b/f)$ follows {from}
  the bracket relation of Theorem~\ref{T:brk} with $m_a = 1, m_b = -
  1$, and $m_c = 0$ for all the other entries in $\bm$.
\end{definition}

\subsection{Main Theorem}

The main result of this paper is that all the $\ok$-linear relations
on $\tpi$ and the values $\Gamma(r),$ for $r \in k \smallsetminus A$
are consequences of the bracket relations.

\begin{theorem} \label{T:MT}
  Let $q > 2$, and let $r_1, \dots, r_n \in k \smallsetminus A$.  Then
  the values $1$, $\tpi$, $\Gamma(r_1), \dots, \Gamma(r_n)$ are
  $\ok$-linearly independent unless for some $1 \le i < j \le n$,
  $\Gamma(r_i) \approx \Gamma(r_j)$.
\end{theorem}

There are a number of noteworthy corollaries.  In particular, we
extend Sinha's result to all possible rational arguments~$r$.

\begin{corollary} \label{C:nonintegral} For each $r \in k
  \smallsetminus A$, $\Gamma(r)$ is transcendental.
\end{corollary}

Another easily stated consequence is the following:

\begin{corollary} \label{C:powerdenom}
  Let $q > 2$.  Let $r_1, \dots, r_n \in k$ be distinct with prime
  power denominators and with the numerator of each $r_i$ having
  degree less than that of the denominator of $r_i$.  Then the values
  $1$, $\tpi$, $\Gamma(r_1), \dots, \Gamma(r_n)$ are $\ok$-linearly
  independent.
\end{corollary}

This corollary can be extended in the following manner:

\begin{corollary} \label{C:moredenoms}
  Let $q > 2$.  Let $f = \prod f_i^{e_i}$ be the decomposition of $f$
  in $A$ in terms of distinct irreducible factors $f_i$.  If no $f_i$
  divides any $f_j - 1$, then the numbers
\[
  1,\, \tpi,\, \Gamma(a/f), \quad a \in A,\ 0 \leq \deg a < \deg f,
\]
are $\ok$-linearly independent.
\end{corollary}

Note that we specifically allow $(a,f) \ne 1$.

In the case that $q=2$, Thakur has observed in \cite{Th91}, {\S}6, that
the bracket relations impose strong conditions on Gamma values and
that in fact $\Gamma(r) \thicksim \tpi$ for all $r \in k
\smallsetminus A$.

The formulation of our main result clearly resembles that of Satz 4 of
\cite{ww85}, where the analogous result is proven for values of the
classical Beta function at rational points.  There the analogue of the
bracket relations of Anderson-Deligne-Thakur are the relations on the
Beta values which arise from the Deligne-Koblitz-Ogus relations for
values of the classical Gamma function.  The analogue of our result
for values of the classical Gamma function itself is still unknown
except for very special cases due to Th.~Schneider and
G.V.~Chudnovsky.

Seen on a large enough scale, the proofs here and in \cite{ww85} also
run somewhat parallel, based as they are on J.\ Yu's Theorem of the
Sub-$t$-module \cite{Yu97}, reproduced below as Theorem
\ref{T:Yusubt}, and G.\ W\"ustholz's Theorem of the Subgroup
\cite{wue89}, respectively.  Luckily, as stated above, the theory of
bracket relations even provides the analogue of the
Deligne-Koblitz-Ogus characterization for the algebraicity of the
product of values of the (normalized) classical Gamma function at
rational points predicted by the classical relations.
 
However some of the crucial tools of \cite{ww85} were not available
for application to $t$-modules.  In particular, we lacked analogues of
the following:
\begin{enumerate}
\item[(a)] Poincar\'e's complete decomposability (up to isogeny) of abelian
  varieties into products of simple ones.
\item[(b)] The Shimura-Tanayama criterion for the explicit
  decomposition of the Jacobian of the Fermat curve into simple
  varieties of CM-type; indeed the very notion of a Jacobian is
  missing from our context.
\item[(c)] An interpretation of all Beta values at rational points as
  abelian integrals.  (Thakur and later Sinha provide a full analogue
  only when $q = 2$.)
\end{enumerate}
In the next section, we describe in general terms how we proceed in
this paper.

\subsection{Outline of the Paper}
In order to apply the transcendence machine embodied in the Theorem of
the Sub-$t$-module, we obviously require appropriate $t$-modules.

In Section \ref{S:tmod} we review some basic definitions.  We then
introduce $t$-modules of CM-type and show that, up to isogeny, they
are always powers of simple $t$-modules of CM-type, thus developing a
serviceable version of Poincar\'e's theorem.

We also give natural criteria, in terms of the underlying CM
structure, for determining the simplicity of $t$-modules of CM-type
and for determining whether the simple $t$-modules underlying two
given $t$-modules of CM-type are isogenous.  These criteria, although
of a vastly different nature, play a role in our transcendence
considerations somewhat analogous to the above mentioned
Shimura-Taniyama criterion.

In Section \ref{S:qp} we define biderivations for arbitrary
$t$-modules $E$ and construct their associated quasi-periodic
extensions $Q$.  We show that, when the biderivations represent
linearly independent classes modulo the inner biderivations, the
resulting quasi-periodic extension is minimal in a precise sense.  The
exponential function of $Q$ comprises the components of the
exponential function of $E$ as well as the quasi-periodic functions
associated to the representative biderivations (plus new variables
corresponding to $\Ga^j$).

We note that an isogeny $T: E_1 \to E_2$ between $t$-modules lifts to
an isogeny $T_\ast :Q_1 \to Q_2$ between the minimal quasi-periodic
extensions of maximal dimension.  As a result, using our special
version of Poincar\'e's theorem, we can determine explicitly when
quasi-periodic extensions of $t$-modules of CM-type are themselves
isogenous.

Section \ref{S:sol} is the heart of the paper. The main work of
\cite{SS95}, \cite{SS97d} involves the construction and investigation
of a $t$-module $E_f$ whose periods have coordinates which are
algebraic multiples of the quantities
\begin{equation} \label{E:sinper}
\Gamma(a/f), \quad \text{$a \in A_{+}$, $\deg(a) < \deg(f)$, $(a,f) =
        1$,}
\end{equation}
where $f \in A_{+}$ is fixed.  Sinha does this by generalizing
Drinfeld's shtuka version of Drinfeld modules to higher dimensional
$t$-modules.  Starting with Anderson's two-variable ``soliton''
function, which generalizes some one-variable results of R.~Coleman,
Sinha creates a module $E_f$ which is of CM-type.

We observe that the bracket relations of Section~\ref{SS:brk} reflect
precisely the CM-structure of $E_f$.  At this point we could deduce
that the set of $\Gamma(a/f)$, $a$, $f \in A_{+}$, $a \in \cU_f$, is
$\overline{k}$-linearly independent if and only if no pair of these
Gamma values is $\overline{k}$-linearly dependent.
 
To provide a setting in which the remaining Gamma values occur, we
consider certain quasi-periodic $t$-module extensions $Q_f$ of $E_f$
by $\Ga^j$:
\begin{equation} \label{E:Qj}
0 \to \Ga^j \to Q_f \to E_f \to 0.
\end{equation}
We compute the periods of these $Q_f$ and see that their components
are algebraic multiples of all the values
\[
\Gamma(a/f), \quad \text{$\deg(a) < \deg(f)$, $(a,f) = 1$,}
\]
where $a$, $f \in A$.  This removes the restriction that
$a$ and $f$ be monic.

In Section \ref{S:indep} we prove independence results about periods
of quasi-periodic extensions of $t$-modules of CM-type.  Basic
properties of minimal extensions enable us to apply Yu's theorem to a
product $Q_{f_1} \times \dots \times Q_{f_m}$ of quasi-periodic
$t$-modules $Q_{f_j}$ of the sort in \eqref{E:Qj} .  Each such
quasi-periodic extension $Q_f$ is minimal and is in fact, up to
isogeny, the power of a minimal extension of a $t$-module of CM-type.
Therefore, the question of linear independence of the coordinates of a
period of $Q_{f_1} \times \dots \times Q_{f_m}$ is reduced to the
question of the corresponding periods of $E_{f_1} \times \dots \times
E_{f_m}$.  Theorem \ref{T:MT} then follows directly from the isogeny
criterion and Yu's theorem.

In Section \ref{S:exmp} we present a few examples.

\subsection*{Acknowledgements}  The authors have greatly benefited from
the encouragement of various colleagues.  In particular, we would like
to thank the referee, as well as Greg Anderson, David Goss, Marius van
der Put, and Dinesh Thakur, who have provided advice and invaluable
assistance at various turns.

\section{$t$-Modules and $t$-Motives of CM-Type} \label{S:tmod}

\subsection{General Definitions} \label{SS:tmod}

We review some basic facts about $t$-modules and $t$-motives in order
to establish notation.  Complete accounts of the material sketched in
this section are contained in the standard references \cite{And86} and
Chapter 5 of \cite{Goss98}.  Continuing with the notation of
Section~\ref{S:int}, let $k \subset L \subset \C$ with $L$
algebraically closed.

\begin{remark}
  We will also need another copy of the pair $A$, $k$ which we keep
  separate, as they will be associated with ``operators'' rather than
  the scalars of $\C$.  We denote the new variable by $t$ and the
  polynomial ring and the fields by the Euler fonts: $\eA \assign
  \Fq[t]$, $\ek \assign \Fq(t)$.  We let $\iota : \ek \to k$ be the
  isomorphism fixing $\Fq$ and sending $t \mapsto \theta$, and we fix
  an extension $\iota \colon \overline{\ek} \to \overline{k}$.  In
  general, the fonts, $\eA$, $\eB$, $\eK$ and so on will be reserved
  for rings of operators corresponding under $\iota$ to $A, B, K$ and
  so on.  We will not maintain this font distinction for elements, and
  so whether we consider $f \in A$ or $f \in \eA$ will depend on the
  context.  Nevertheless, we will persist in the distinction between
  $\theta \in A$ and $t \in \eA$.
\end{remark}

Let $\tau$ denote the $q$-th power Frobenius map: $x \mapsto x^q$, $q
= p^r$.

Let $E$ be an algebraic group defined over $L$ isomorphic to $\Ga^d$.
We take $\Lie(E)$ for its tangent space at the origin and note that,
after choosing a basis for this isomorphism, $\Lie(E) \simeq L^d$.
Similarly, if $\End_L^q(E)$ is the ring of $\Fq$-linear endomorphisms
of $E$ as an algebraic group over $L$, then selecting a basis induces
an isomorphism
\[
\End_L^q(E) \simeq \Mat_{d \times d}(\Lt) \assigN
\Mat_{d \times d}(L)\{ \tau \},
\]
where we denote by $L\{ \tau \}$ the non-commutative ring of
\emph{twisted polynomials in $\tau$}, for which $\alpha \tau^i \beta
\tau^j = \alpha \beta^{q^i}\tau^{i+j}$ when $\alpha, \beta \in L$.

\subsubsection{$t$-Modules and $t$-Motives}
A \emph{$t$-module over $L$} is an algebraic group $E$, isomorphic to
$\Ga^d$ over $L$, for which there is an $\Fq$-linear homomorphism
\[
\Phi : \eA \to \End^q_L(E)
\]
and an $\ell >0$ such that the endomorphism $\Phi(t) - \theta\tau^0$
satisfies 
\[
(\Phi(t) - \theta\tau^0)^{\ell} \Lie(E) = \{0\}.
\]

We refer to $d$ as the \emph{dimension} of $E$.  A sub-$t$-module $H$
of $E$ is a connected sub-algebraic group of $E$ which is also
invariant under the action of $\Phi(t)$.  A $t$-module is said to be
\emph{simple} if it does not contain any proper sub-$t$-modules.  A
\emph{morphism} $\Theta : E_1 \to E_2$ of $t$-modules over $L$ is a
morphism of algebraic groups over $L$ which commutes with the actions
of $t$ on $E_1$ and $E_2$.  We let $\Hom(E_1,E_2)$ denote the group of
$t$-module morphisms $E_1 \to E_2$ and also take $\End(E) \assign
\Hom(E,E)$.  A morphism $\Theta$ will be called an \emph{isogeny} if
the dimensions of $E_1$ and $E_2$ are the same and if $\Theta$ has a
finite kernel as a map of algebraic groups.  We will write $E_1
\thicksim E_2$ to denote that $E_1$ is isogenous to $E_2$; isogeny is
an equivalence relation (see \cite{Yu97}, Lem.~1.1).

By choosing an isomorphism $E \simeq \Ga^d$, we obtain an $\Fq$-linear
homomorphism
\[
\Phi : \eA \to \Mat_{d \times d}(L)\{ \tau \}
\]
which provides the action of $\eA$ on $E$.  If $\Phi$ is defined by
\[
\Phi(t) = M_0\tau^0 + M_1\tau + \dots + M_n \tau^n
\]
with the $M_i \in \Mat_{d \times d}(L)$, then the induced action on
$\Lie(E)$ is given by $d\Phi(t) \assign M_0$.  The definition of a
$t$-module dictates that $d\Phi(t) = (\theta I_d + N)$ for the $d
\times d$ identity matrix $I_d$ and a nilpotent matrix $N$.  To
signify that we have chosen a system of coordinates for $E$ we will
write $E = (\Phi,\Ga^d)$.  Thus if $E_1 = (\Phi_1,\Ga^{d_1})$, $E_2 =
(\Phi_2,\Ga^{d_2})$ are two $t$-modules, a morphism $\Theta : E_1
\rightarrow E_2$ is a matrix of twisted polynomials, $\Theta \in
\Mat_{d_2\times d_1}(L\{\tau\})$, satisfying
\[
\Theta \Phi_1(t) = \Phi_2(t)\Theta.
\]

The dual notion of a $t$-module is that of a $t$-motive.  Set
$L[t,\tau] \assign L\{\tau\}[t]$, the ring of commuting polynomials in
the variable $t$ over the non-commutative ring $L\{\tau\}$.  Then a
\emph{$t$-motive} $M$ is a left $L[t,\tau]$-module which is free and
finitely generated as an $L\{\tau\}$-module and for which
\[
(t - \theta)^{\ell} (M / \tau M) = \{ 0 \},
\]
for some $\ell > 0$.  Morphisms of $t$-motives are morphisms of
left $L[t,\tau]$-modules.

To every $t$-module $E = (\Phi,\Ga^d)$ over $L$, there corresponds a
unique $t$-motive over~$L$:
\[
M \assign M(E) \assign \Hom_L^q(E,\Ga),
\]
where $\Hom_L^q(A,B)$ denotes the $L$-module of $\Fq$-linear morphisms
of the algebraic groups $A$, $B$ over $L$.  In this setting, the
action of $f t^i$, $f \in L\{\tau\}$ on $m \in M$ is
\[
(f t^i, m) \mapsto f \circ m \circ \Phi(t^i).
\]
Projections on the $d$ coordinates of $E \simeq \Ga^d$ form an
$L\{\tau\}$-basis for $M$; $d = \rank_{L\{\tau\}} M$; and we can take
$\ell \leq d$.

It is a fundamental theorem of Anderson that the functor $E \mapsto
M(E)$ gives an anti-equivalence from the category of $t$-modules over
$L$ to the category of $t$-motives over $L$.  Given a $t$-motive $M$
together with an $L\{\tau\}$-basis $m_1,\dots,m_d$ for $M$, 
we can express the $t$-action with respect to this basis:
\begin{equation} \label{E:Phi}
t\cdot 
\begin{pmatrix}
  m_1 \\ \vdots \\ m_d
\end{pmatrix} = 
\Phi(t) \begin{pmatrix}
m_1\\ \vdots \\ m_d
\end{pmatrix},
\end{equation}
where $\Phi(t) \in \Mat_{d\times d}(L\{\tau\})$.  

We verify this equivalence of categories in the following way:
Elements $m$ of $M$ correspond uniquely to $\bsU = (U_1,\dots,U_d) \in
\Mat_{1\times d}(L\{\tau\})$ via
\[
m = (U_1,\dots,U_d)\cdot 
\begin{pmatrix}
  m_1 \\  \vdots \\
  m_d
\end{pmatrix}.
\] 
According to the
commutativity of $t$ with elements of $L\{\tau\}$, 
\[
t \cdot \bsU \begin{pmatrix} m_1 \\ \vdots \\ m_d
  \end{pmatrix} = 
\bsU \cdot t \begin{pmatrix} m_1 \\ \vdots \\ m_d
  \end{pmatrix} = 
\bsU \Phi(t) \begin{pmatrix} m_1 \\ \vdots \\ m_d
  \end{pmatrix}.
\]
In other words, we can take the action of $t$ on $\Mat_{1\times
  d}(L\{\tau\})$ to be
\[
t\cdot \bsU = \bsU \Phi(t).
\]
If we now define $E \assign (\Phi, \Ga^d)$, then the action of $t$ on
$M$ with respect to the basis $m_1, \dots, m_d$ is that of $M(E)$.  In
particular, when the $t$-motive $M$ is defined over $L$, so is its
corresponding $t$-module $E$.  This remark will be essential for our
transcendence considerations below.

\subsubsection{Exponential Function and Uniformization} \label{SS:Exp}
%wdb
Given a $t$-module $E = (\Phi,\Ga^d)$ 
defined over $L$, there is an
associated {\em exponential function}
\[
\Exp = \Exp_E : \Lie(E) \to E(\C),
\]
given by the unique power series which satisfies the relations:
\begin{enumerate}
\item[(a)] $\Exp(d\Phi(t)\bz) = \Phi(t)\Exp(\bz)$.
\item[(b)] $\partial \Exp(\bz) = I_d$.
\end{enumerate}
Here $\partial \Exp(\bz)$ denotes the matrix of coefficients of linear
terms in $\Exp(\bz)$.  The function $\Exp$ is an entire $\Fq$-linear
function $\Exp: \C^d \to \C^d$ with coefficients from $L$.  (A power
series is $\Fq$-\emph{linear} if it is a sum of power series in single
variables in which each exponent is a power of $q$.)

Elements of $\Lambda \assign \Ker \Exp$ are called \emph{periods} of
$E$.  A major effort of the paper is devoted to the definition of
certain $t$-modules whose periods involve the special Gamma values.

\subsubsection{Abelian $t$-Modules and Uniformization} \label{SS:atm}

The $t$-module $E$ is called \emph{uniformizable} if the image of
$\Exp$ is all of $\Ga^d(\C)$.  Conditions for surjectivity of the
exponential function is a quite fascinating topic; we hope to give new
criteria for uniformizability in a future note.

When a $t$-motive $M$ is finitely generated over $L[t]$, it and the
corresponding $t$-module $E$ are said to be {\em abelian}.  Anderson's
definition of $t$-motive in \cite{And86}, {\S}1.1.2, includes the
abelian condition just given.  For the sake of clarity, we point out
that in this matter we follow the terminology in Goss \cite{Goss98},
{\S}5.4, to allow consideration of the more general situation when
appropriate.

Now $M$ is abelian precisely when it is free of finite rank over
$L[t]$; we call this rank the {\em rank} of $M$ and $E$ and denote it
by $r(M) = r(E)$.  By a theorem of Anderson \cite{And86}, Thm.~4,
and \cite{Goss98}, Thm.~5.9.14, the uniformizability of an abelian
$t$-module $E$ is equivalent to the condition
\begin{equation} \label{E:gwacrit}
        \rank_{\eA} (\Ker \Exp_E) = \rank_{\C[t]}M(E) \assigN r(E).
\end{equation}
Furthermore, if $E$ is uniformizable, then a theorem of
Anderson~\cite{And86}, Cor.~3.3.6, shows that $\Lambda$ spans
$\Lie(E)$ over $\C$.

\begin{remark_num} \label{R:AbUn}
  It should be noted that if $E$ is abelian and uniformizable, then
  every sub-$t$-module $H$ is also abelian and uniformizable.  Indeed,
  it is easy to see that if $H$ is a sub-$t$-module of $E$, then $E$
  is abelian if and only if both $H$ and the quotient $t$-module $E/H$
  are abelian as well.  By Yu~\cite{Yu89}, Prop.~5.3, we know that $H$
  must be uniformizable if $E$ is.
\end{remark_num}

Let $E = (\Ga^d,\Phi)$ be an abelian $t$-module.  If $T$ is a finite
subgroup of $E(\C)$, which is invariant under the action of $\Phi(t)$,
then the quotient $E/T$ (as algebraic groups) is naturally given the
structure of a $t$-module (see \cite{Goss98}, {\S}5.6).  The map of
$t$-modules $E \to E/T$ is a surjective map of algebraic groups.  If
$E$ is uniformizable, it then follows that $E/T$ is also
uniformizable.

\begin{lemma} \label{L:Elatt}
  Let $E_1 = (\Phi_1,\Ga^d)$ be a uniformizable abelian $t$-module
  with period lattice $\Lambda_1$.  Let $\Lambda_2 \subset \Lie(E_1)$
  be an $\eA$-lattice in which $\Lambda_1$ has finite index.  Then
  $\Lambda_2$ is the period lattice of a uniformizable abelian
  $t$-module $E_2 = (\Phi_2,\Ga^d)$, and there is a natural isogeny $E_1
  \to E_2$.
\end{lemma}

\begin{proof}
  Let $T \subset E_1(\C)$ be the image of $\Lambda_2$ under
  $\Exp_{E_1}$.  Then $T$ is finite and invariant under $\Phi_1(t)$.
  Let $E_2 \assign E_1/T$.  We can identify the tangent spaces of $E_1$ and
  $E_2$, and by the functoriality of exponential functions, the
  exponential function of $E_2$ is the map
\[
\Exp_{E_2} : \Lie(E_1) \overset{\Exp_{E_1}}{\longrightarrow} E_1(\C)
\to E_2(\C),
\]
which is surjective and has kernel $\Lambda_2$.
\end{proof}

\subsection{Hilbert-Blumenthal-Drinfeld Modules} \label{SS:HBD}

For fields $\eK \subset \eL$, let $\Emb(\eL/\eK)$ denote the set of
embeddings $\sigma \colon \eL \to \C$ such that $\sigma|_{\eK} =
\iota|_{\eK}$.  Let $\eK_+/\ek$ be a finite separable extension with
$[\eK_+:\ek] = d$ such that the place $\infty$ of $\ek$ is totally
split in $\eK_+$, and take $\eB_+$ for the integral closure of $\eA$
in $\eK_+$.  Let $\{\sigma_1, \dots, \sigma_d\} = \Emb(\eK_+/\ek)$.
  Define the \emph{conjugate action} $\bss \assign \sigma_1 \oplus
  \cdots \oplus \sigma_d$ of $\eK_+$ on $\C^d$ in the following
  manner: For $\eb \in \eK_+$ and $\bz \assign (z_1,\dots,z_d) \in
  \C^d$,
\begin{equation} \label{E:conj}
\bss(\eb): (z_i) \mapsto \bss(\eb)(\bz) \assign 
(\sigma_i(b) z_i).  
\end{equation}

Suppose now that $E = (\Phi,\Ga^d)$ is a uniformizable abelian
$t$-module and that $\Phi$ extends to a map of $\eA$-algebras
\[
\Phi : \eB_+ \to \End(E)
\]
in such a way that the action of $d\Phi$ on $\Lie(E)$ is given by
\[
d\Phi(\eb) = \bss(\eb), \ \forall b \in \eB_+.
\]
Then $E$ is called a \emph{Hilbert-Blumenthal-Drinfeld module} (H-B-D
module) with multiplications by $\eB_+$.  By definition, the period
lattice $\Lambda$ of $E$ is invariant under $\bss(\eB_+)$.  Working
with Hilbert-Blumenthal-Drinfeld modules is greatly facilitated by the
following equivalence.

\begin{theorem}[Anderson {\cite{And86}, Thm.~7}] \label{T:HBD}
The following two categories are equivalent.
\begin{enumerate}
\item[(a)] Objects: H-B-D modules with multiplications by $\eB_+$.

      \noindent Morphisms: $t$-module homomorphisms which are
      $\eB_+$-equivariant.

\item[(b)] Objects: Lattices contained in $\C^d$ invariant under
      $\bss(\eB_+)$.

      \noindent Morphisms: $\C$-linear maps on $\C^d$ which carry one
      lattice into the other and commute with $\bss(\eB_+)$.
\end{enumerate}
\end{theorem}

\subsection{$t$-Modules of CM-type} \label{SS:CM}

One particular sort of H-B-D module, which we call \emph{$t$-modules
  of CM-type}, are of prime importance in this paper.  We continue
with the notation of the previous section.  Let $\eK$ be a finite
separable extension of $\eK_+$ which is totally ramified at each
infinite place of $\eK_+$.  In this situation, we say that $\eK_+$ is
the {\em maximal real subfield of $\eK$}.

Choose extensions of $\sigma_1, \dots, \sigma_d$ to embeddings $\eK
\hookrightarrow \C$ which we also denote by $\sigma_1, \dots,
\sigma_d$.  Letting 
\begin{equation} \label{E:S}
\cS \assign \{ \sigma_1, \dots, \sigma_d \} \subset \Emb(\eK/\ek)
\end{equation} 
be the set of these extensions, we denote the extension of the
conjugate action $\bss$ to $\eK$ by $\bss_{\cS}$.

Let $\Lambda$ be a discrete $\eA$-submodule of $\C^d$ of rank
$[\eK:\ek]$ which is invariant under the action, via $\bss_{\cS}$, of
some order of $\eK$.  We denote by $\eB$ the maximal such order.  As a
$\eB$-module, $\Lambda$ is isomorphic to an ideal of $\eB$.  We think
of $\Lambda$ as having {\em real multiplications} by the conjugate
action $\bss(\eB_+ \cap \eB)$ and {\em complex multiplications}
by the conjugate action $\bss_{\cS}(\eB)$, {\em CM by $\eB$} for
short.

By Anderson's Theorem~\ref{T:HBD}, $\Lambda$ contains a sublattice of
finite index which is the period lattice of a H-B-D module with
multiplications by $\eB_+$.  By Lemma~\ref{L:Elatt}, $\Lambda$ is then
itself the period lattice of a uniformizable abelian $t$-module $E =
(\Phi, \Ga^d)$.  Because $\bss_{\cS}(\eB) \Lambda \subset \Lambda$ and
the action of $\eB$ via $\bss_{\cS}$ commutes with the action of
$\eB_+ \cap \eB$ via $\bss$, Anderson's theorem also gives an
extension of $\Phi$ to all of $\eB$ such that
\[
d\Phi(b) = \bss_{\cS}(b)
\]
for all $b \in \eB$.  In other words, we have an injection of
$\eA$-algebras
\begin{equation} \label{E:PhiB}
\Phi : \eB \hookrightarrow \End(E),
\end{equation}
which extends the $t$-module homomorphism $\Phi$ on $\eA$.  We call
$\eK$ the \emph{CM-field} of $E$ and say that $E$ has \emph{CM-type
$(\eK,\cS)$.}

\begin{remark_num} 
  In the case that $\Lambda$ is isomorphic to a principal ideal of
  $\eB$, fixing a generator $\bsl = (\lambda_1,\dots,\lambda_d)^{tr}
  \in \C^d$ allows us to identify $\eB$ with $\Lambda \subset \C^d$
  via $\bss$:
\begin{equation} \label{E:BLmb}
\begin{aligned}
        \eB &\stackrel{\sim}{\longrightarrow} \Lambda \subset \C^d \\
        \eb &\mapsto \begin{pmatrix} \sigma_1(\eb)\lambda_1 \\ \vdots
        \\ \sigma_d(\eb)\lambda_d \end{pmatrix} \mathrel{=:}
        \bss_{\cS}(\eb)(\bsl).
\end{aligned}
\end{equation}
Yu \cite{Yu89}, Lem.~6.2, shows that all the coordinates $\lambda_i$ of
$\bsl$ are non-zero.
\end{remark_num}

\begin{remark_num} \label{R:CMisog}
  Suppose that the $t$-module $E_1$ has CM by $\eB_1$, with period
  lattice $\Lambda_1$ and CM-type $(\eK,\cS)$.  Let $\beta \in \eB_1$
  and $\bsl \in \Lambda_1$ be non-zero.  Then, as $\bss(\eB_1) \bsl$
  has finite index in $\Lambda_1$, there is a non-zero $a \in \eA$
  such that $a \Lambda_1 \subset \bss(\eB_1) \bsl$.  Thus, if $\eB$
  denotes the ring of integers of $\eK$, the $\eB_1$-lattice indices
  of the tower
\[
a\Lambda_1 \subset \bss(\eB_1) \bsl \subset \bss(\eB) \bsl
\]
are finite.  As the left- and right-most lattices are isomorphic to
$\Lambda_1$ and $\bss(\eB)$, respectively, we see by Anderson's
Theorem \ref{T:HBD} above and by Lemma~\ref{L:Elatt}, that $E_1$ is
isogenous to a $t$-module $E$ with lattice isomorphic to the full ring
of integers $\eB$ of $\eK$ (and thus of CM-type $(\eK,\cS)$ and with
CM by $\eB$).
\end{remark_num}

\begin{remark_num}  In particular, if two $t$-modules have the
  same CM-type $(\eK,\cS)$, then they are necessarily isogenous.
\end{remark_num}

\subsection{Endomorphism Rings} \label{SS:End}

Throughout this section we will assume that $E = (\Phi, \Ga^d)$ is an
abelian $t$-module with period lattice $\Lambda$.  Let $\End^0(E)
\assign \End(E) \otimes_{\eA} \ek$.  Certainly $\Phi$ extends to a map
$\Phi : \ek \to \End^0(E)$.  Our goal of this section is to determine
the structures of $\End(E)$ and $\End^0(E)$.

\begin{lemma}
  Let $E$ be simple.  Then the ring $\End^0(E) \assign \End(E)
  \otimes_{\eA} \ek$ is a division algebra and $\ek$ lies in its
  center.
\end{lemma}

\begin{proof}
$\End(E)$ has no non-zero zero divisors because $E$ is assumed to be
simple.  Now if $f \in \End(E)$, then the kernel of $f$ is a
sub-algebraic group invariant under the action of $t$.  Its connected
component will be a sub-$t$-module of $E$, and so the kernel of $f$
must be finite.  Thus Lemma~1.1 of~\cite{Yu97} gives $\ea \in \eA$ and
$g \in \End(E)$ so that, (as elements of $\End(E)$),
\[
        fg = \Phi(\ea).
\]
Therefore $\End^0(E)$ is a division algebra, and its center clearly
contains $\ek$.
\end{proof}

\begin{proposition}
Suppose $E$ is isogenous to a product $E_1^{n_1} \times \cdots \times
E_k^{n_k}$ of powers of pair-wise non-isogenous, simple abelian
$t$-modules.  Then
\[
        \End^0(E) \simeq \bigoplus_{i=1}^k \Mat_{n_i \times n_i}
                                ( \End^0(E_i)).
\]
\end{proposition}

\begin{proof} Replace the word ``isomorphism'' in the statement and
  proof of Proposition XVII.1.2 of \cite{lang} by ``isogeny''.
\end{proof}

For the rest of this section we assume that $E$ is simple, let $s
\assign \rank_{\eA} \Lambda$ and assume $s > 0$.  Unless $E$ is
uniformizable, this rank $s$ differs from the rank $r$ of $E$ as an
abelian $t$-module.

Because $E$ is simple, the image of each non-zero element $f \in
\End(E)$ is all of $E$.  Therefore the kernel of such an $f$ is
finite, and $f$ must take $\Lambda$ to a sublattice which has finite
index in $\Lambda$.  Hence there is an injection of $\eA$-algebras
$\End(E) \hookrightarrow \Mat_{s \times s}(\eA)$ such that the image
of every non-zero endomorphism is invertible in $\Mat_{s\times
s}(\ek)$.  This map extends uniquely to a map of $\ek$-algebras
\begin{equation} \label{E:ratl}
  \End^0(E) \hookrightarrow \Mat_{s \times s}(\ek),
\end{equation}
which is called the \emph{rational representation} of $\End^0(E)$.  We
see right away that $\End(E)$ is a free $\eA$-module of rank over
$\eA$ at most $s^2$ and that $\dim_{\ek} \End^0(E) \leq s^2$.

\begin{proposition} \label{P:End}
  Let $E$ be a simple abelian $t$-module with period lattice of rank
  $s$ over $\eA$.  Let $\eK_0$ be the center of $\End^0(E)$; $g =
  [\eK_0:\ek]$; and $h^2 = [\End^0(E):\eK_0]$.  Then $gh^2 \mid s$.
\end{proposition}

\begin{proof}
As the rational representation of \eqref{E:ratl} is a faithful
representation of the division algebra $\End^0(E)$, it follows that
$[\End^0(E):\ek]$ must divide $s$ (see \cite{lang}, Prop.~XVII.4.7).
\end{proof}

\begin{corollary} \label{C:EHBD}
  Let $E = (\Phi, \Ga^d)$ be a simple Hilbert-Blumenthal-Drinfeld
  module of CM-type $(\eK,\cS)$ with complex multiplications by $\eB$.
  Then complex multiplication gives an isomorphism $\Phi : \eB
  \stackrel{\sim}{\longrightarrow} \End(E)$ which lifts to an
  isomorphism $\eK \stackrel{\sim}{\longrightarrow} \End^0(E)$. 
\end{corollary}

\begin{proof}
  Let $\eK_0$, $g$, $h$, $s$ be as in Proposition \ref{P:End}.  Now
  $\Phi(\eK)$ is a subfield of $\End^0(E)$, and $[\Phi(\eK):\Phi(\ek)]
  = [\eK:\ek] = r$, and $r = s$, as $E$ is uniformizable.  Since the
  centralizer of $\Phi(\eK)$ contains $\eK_0$, we conclude that
  $\Phi(\eK)\eK_0$ is a subfield of $\End^0(E)$.  Because each maximal
  commutative subfield of $\End^0(E)$ has dimension $h$ over $\eK_0$
  (see \cite{Goss98}, Cor.~4.11.15), we conclude that
  $[\Phi(\eK)\eK_0:\Phi(\ek)] \leq gh$.  By Proposition~\ref{P:End} we
  know $gh^2 \mid r = [\eK : \ek]$, and so
\[
  r = [\Phi(\eK):\Phi(\ek)] \leq [\Phi(\eK)\eK_0:\Phi(\ek)] \leq gh
  \leq gh^2 \leq r.
\]
Hence $r = gh = gh^2$ and so $h = 1$; thus $\End^0(E) \simeq \eK$.
Now $\End(E)$ is isomorphic to an order in $\eK$, and since $\eB$ is
the largest order ${\mathcal O}$ for which $\bss_{\cS}({\mathcal O})$
leaves $\Lambda$ invariant, we must have $\End(E) = \Phi(\eB)$.
\end{proof}

\begin{remark}
It is possible to prove Corollary~\ref{C:EHBD} more directly.  If $f
\in \End(E)$, then $df$ leaves $\Lambda$ invariant.  Because
$\bss_{\cS}(\eK)\Lambda = \bss_{\cS}(\ek)\Lambda$, we see that
$df$ leaves the $1$-dimensional $\eK$-vector space
$\bss_{\cS}(\eK)\Lambda$ invariant.  Any non-zero element $\bsl \in
\Lambda$ is a generator for this vector space, and so there exists a
unique $b \in \eK^{\times}$ such that
\[
df(\bsl) = \bss_{\cS}(b)(\bsl).
\]
We know that $\bss(b) \in \End^0(E)$ and we want to show, among other
things, that $b \in \eB$.  However we know that, for some denominator
$a \in \eA$, we have $ab \in \eB$.  For $g \assign \Phi(a)f - \Phi(ab)
\in \End(E)$, the displayed line shows that $dg$ has a non-trivial
kernel ($dg \bsl = 0$).  Because $E$ is simple, every non-zero
endomorphism is an isogeny and so is an isomorphism on $\Lie(E)$.
Thus $g = 0$, i.e.\ $\Phi(a)f = \Phi(ab),$ and thus 
$\Phi(b) = f \in \End(E)$. As $\bss_{\cS}(b)\Lambda = 
d\Phi(b)\Lambda = df\Lambda \subset \Lambda,$ we see that $b \in \eB$.
\end{remark}

\subsection{Sub-$t$-modules and Isogenies} \label{SS:CMsub}

In this section, we investigate criteria for determining when two
$t$-modules of CM-type $E_1$ and $E_2$ have a non-trivial $t$-module
morphism between them or, what will be equivalent, when the two have
isogenous non-trivial sub-$t$-modules.  When $\eK_+/\ek$ is Galois,
we describe completely the sub-$t$-module structure of $E$ and thus
ascertain what sorts of $t$-module morphisms exist between $t$-modules
of CM-type.

Let $E$ be a $t$-module of CM-type as in Section \ref{SS:HBD}, with
multiplication rings $\eB$ and $\eB_+$ and fraction fields $\eK$ and
$\eK_+$ as described there, with the added assumption of this section
that $\eK_+$ is a Galois extension of $\ek$.  Let $E$ have CM-type
$(\eK,\cS)$ with $ \cS \assign \{ \sigma_1, \dots, \sigma_d \} \subset
\Emb(\eK/\ek)$.  The following lemma is the basis for the remainder of
our discussion.

\begin{lemma} \label{L:crit}
  Let $\eL \subset \eK$ be a subfield, and set $\eL_+ \assign \eL \cap
  \eK_+$.  Assume that $\eK_+/\ek$ is Galois.  Then the following
  are equivalent:
\begin{enumerate}
\item[(a)] $[\eL:\eL_+] = [\eK:\eK_+]$; and for all $i$, $j$, if
  $\sigma_i|_{\eL_+} = \sigma_j|_{\eL_+}$, then $\sigma_i|_{\eL} =
  \sigma_j|_{\eL}$.
\item[(b)] $\cS = \cup_i \Emb(\alpha_i(\eK)/\alpha_i(\eL)) \circ \alpha_i$
  for certain $\alpha_i \in \Emb(\eK/\ek)$ with distinct restrictions
  $\alpha_i|_L$.
\end{enumerate}
\end{lemma}

\begin{proof}  
  The proof is based on two straightforward remarks:
  (1)~$\eK_+\eL/\eL$ is a Galois extension with $\Gal(\eK_{+}\eL/\eL)
  \simeq \Gal(\eK_+/\eL_+)$ (cf.\ \cite{lang}, Thm.~VI.1.12).  (2)~As
  $\cS$ consists of the extensions to $\eK$ of the distinct elements
  of $\Gal(\eK_+/\ek)$, precisely $[\eK_+:\eL_+]$ distinct elements of
  $\cS$ restrict to each of the $[\eL_+:\ek]$ embeddings of $\eL_+$
  into $\eK_+$.

(a) $\Rightarrow$ (b): By the first part of (a) and~(1), $[\eK:\eL] =
[\eK_+:\eL_+] = [\eK_+\eL :\eL]$, so $\eK_+\eL = \eK$.
By the second part of (a), elements of $\cS$ in distinct $\cS_\alpha
\assign \Emb(\alpha(\eK)/\alpha(\eL))\circ \alpha$ restrict to give
distinct embeddings of $\eL_+$.  According to (2), $\cS$ contains
precisely $[\eK_+:\eL_+]$ elements restricting to the same embedding
of of $\eL_+/\ek$.  According to (a), these elements restrict to the
same embedding of $\eL/\ek$.  As $[\eK:\eL] = [\eK_+:\eL_+]$, these
are all the embeddings of $\eK$ restricting to the given embedding of
$\eL/\ek$, i.e.\ these elements constitute a subset of the form
$\cS_\alpha$.  This proves (b).

(b) $\Rightarrow$ (a): Since elements of a set $\cS_{\alpha}$ restrict
to the same embedding of $\eL_+$ into $\eK_+$, (2) shows that $\cS$
contains at least $[\eL_+:\ek]$ such sets $\cS_{\alpha}$, and
\[
[\eK:\eL][\eL_+:\ek] \le |\cS| = [\eK_+:\ek] =
[\eK_+:\eL_+][\eL_+:\ek]. 
\]
Since, by (1), $[\eK_+:\eL_+] = [\eK_+\eL:\eL]$, we find $[\eK:\eL]
\le [\eK_+\eL:\eL]$.  Thus we must have that $\eK_+\eL =
\eK$.  Therefore $[\eK:\eL] = [\eK_+:\eL_+]$, and it follows that
$[\eL:\eL_+] = [\eK:\eK_+]$, which is the first part of (a).

Now we know that each $\cS_{\alpha}$ contains $[\eK:\eL] =
[\eK_+:\eL_+]$ elements, each of which restricts to the same embedding
of $\eL_+$ into $\eK_+$.  By (2) then, elements from distinct
$\cS_{\alpha}$ must restrict to distinct embeddings of $\eL_+$ into
$\eK_+$.  This is the second part of~(a).
\end{proof}

The following proposition provides the sub-$t$-module structure of $E$
and determines when $E$ is simple.  It is similar to a well-known
theorem about abelian varieties of CM-type (see
Mumford \cite{mum70}, {\S}22, pp.~213-14).

\begin{theorem} \label{T:Esub}
Let $E$ be a $t$-module of CM-type $(\eK,\cS) = (\eK,\{ \sigma_1,
\dots, \sigma_d\})$.

\begin{enumerate}
\item[(a)] $E$ is isogenous to a power of a simple sub-$t$-module.
\item[(b)] Assume that $\eK_+/\ek$ is Galois.  Then $E$ is itself
  simple if and only if no proper subfield $\eL \subset \eK$ satisfies
  either of the equivalent criteria of Lemma~\ref{L:crit}.
\end{enumerate}
\end{theorem}

\begin{remark_num} \label{R:subCM}
The simple sub-$t$-module provided by this proposition is itself a
$t$-module of CM-type.  Indeed the conditions of Lemma~\ref{L:crit}
ensure that it has CM-type $(\eL,\cS|_{\eL})$.
\end{remark_num}

\begin{remark_num}
  There are abelian $t$-modules which are not semi-simple even up to
  isogeny.  For example, one can construct non-trivial non-torsion
  extensions of abelian $t$-modules over an algebraically closed
  field, and so for this reason there is no general Poincar\'e
  decomposability theorem for $t$-modules.  However, it is not known
  whether pure $t$-modules are indeed semi-simple up to isogeny.
\end{remark_num}

\begin{proof}[Proof of Theorem~\ref{T:Esub}]
We first prove (a).  Let $H \subset E$ be a simple sub-$t$-module.
Let $\eC \subset \eB$ be the largest $\eA$-algebra for which
$\Phi(\eC)$ leaves $H$ invariant.  Let $\eL$ be the fraction field of
$\eC$; thus $\eC$ is an order of $\eL$.  Let $d_0 = \dim H$; $m =
[\eK:\eL]$; and $n = [\eL:\ek]$.  For some $\beta_1, \dots, \beta_m
\in \eB$, the direct sum
\[
\beta_1\eC + \dots + \beta_m\eC \subset \eB,
\]
is an $\eA$-submodule of finite index.  Now for all $i$,
$\Phi(\beta_i) H$ is a simple sub-$t$-module isogenous to $H$.\

As a sub-$t$-module of a uniformizable abelian $t$-module, $H$ is also
abelian and uniformizable (Remark \ref{R:AbUn}), and for $W \assign
\Lie(H)$ we have $\Exp_H = \Exp_E|_W$.  The action of $\eC$ on $\C^d$
induced by $\bss_{\cS}(\eB)$ leaves both $W$ and the period lattice
$\Lambda_H \subset W$ of $H$ invariant.  For each $i$,
$\Lie(\Phi(\beta_i)H) = \bss_{\cS}(\beta_i)W =: W_i$, and
$\bss_{\cS}(\beta_i)\Lambda_H$ has finite index in the period lattice
of $\Phi(\beta_i)H$.

By a theorem of Anderson \cite{And86}, Cor. 3.3.6, we know both that
$\Lambda$ spans $\Lie(E)$ over $\C$ and that $\bss_{\cS}
(\beta_i)\Lambda_H$ spans $W_i$ over $\C$.  As
$\bss_{\cS}(\beta_1)\Lambda_H + \dots + \bss_{\cS}(\beta_m) \Lambda_H$
has finite index in $\Lambda$, we find that $\Lie(E) = W_1 + \dots +
W_m$.  Moreover,
\begin{equation} \label{E:Esum}
  E = \Phi(\beta_1)H + \dots + \Phi(\beta_m)H.
\end{equation}
Now as each $\Phi(\beta_i)H$ is simple, we can re-order the $\beta_i$
so that for some $\ell \leq m$, we have an isogeny
\begin{equation} \label{E:Epow}
  E \thicksim \Phi(\beta_1)H \times \dots \times \Phi(\beta_{\ell})H.
\end{equation}
Thus $E$ is isogenous to $H^{\ell}$, which proves part (a).  We will
now show that $\ell = m$.

First, we show that the centralizer of $\Phi(\eK)$ in $\End^0(E)$ is
$\Phi(\eK)$.  As a 1-dimensional $\eK$-vector space,
$\bss_{\cS}(\eK)\Lambda = \bss_{\cS}(\eK)\bsl$, for some $\bsl \in
\Lambda$.  For any $f \in \End^0(E)$, there is a $c_f \in \eK$ so that
$df(\bsl) = \bss_{\cS}(c_f)\bsl.$ If $f$ commutes with all of
$\Phi(\eK)$, then in fact, $(df - \bss_{\eS}(c_f))|_{\Lambda} = 0$.
Because $\Lambda$ spans $\Lie(E)$ over $\C$, it follows that $df =
\bss_{\eS}(c_f)$.  For some non-zero $a \in \eA$, $g \assign \Phi(a)
f - \Phi(a) \Phi(c_f) \in \End(E)$ has $dg = 0$, and therefore $g =
0$.  In other words, $f = \Phi(c_f)$.

Now let $Z_0$ be the center of $\End^0(H)$.  From the diagonal mapping
$Z_0 \hookrightarrow \End^0(E) \simeq \Mat_{\ell \times
\ell}(\End^0(H))$, we see that the image of $Z_0$ is in the center of
$\End^0(E)$ and in particular in the centralizer of $\Phi(\eK)$.
Therefore, $Z_0 \subset \Phi(\eK)$ by the above paragraph.  Thus there
is a field $\eL_0 \subset \eK$ such that $Z_0 = \Phi(\eL_0)$.

Let $g = [Z_0 : \Phi(\ek)] = [\eL_0:\ek]$ and $h^2 = [\End^0(H):Z_0]$.
By the general theory of semi-simple algebras (see \cite{langCM},
pp.~11--12), since (1) $\End^0(H)$ is a division algebra and (2)
$\Phi(\eK) \supset \Phi(\eL_0)$ is a subfield of $\End^0(E) \simeq
\Mat_{\ell \times \ell}(\End^0(H))$ which is its own centralizer, it
follows that
\[
  [\eK : \eL_0] = \ell h.
\]
Thus we have
\[
\ell \rank_{\eA}(\Lambda_H) = \rank_{\eA} \Lambda = [\eK:\ek] =
[\eK:\eL_0][\eL_0:\ek] = \ell h [\eL_0:\ek],
\]
and so
\[
  \rank(H) = \rank_{\eA}(\Lambda_H) = gh.
\]
Since $gh^2 \mid \rank_{\eA}(\Lambda_H)$ by Proposition~\ref{P:End},
it follows that $h = 1$.  Therefore $\End^0(H) = \Phi(\eL_0)$ and
$\rank_{\eA}(\Lambda_H) = [\eL_0 : \ek]$.

On the other hand,
\[
  \Phi(\eL) \subset \End^0(H) = \Phi(\eL_0) \subset \Phi(\eK).
\]
Now $\End(H) \cap \Phi(\eB)$ stabilizes $H$ and is an order in
$\Phi(\eL_0)$.  Our choice that $\eC$ be the largest $\eA$-subalgebra
of $\eB$ for which $\Phi(\eC)$ leaves $H$ invariant then implies that
$\eC \supset \End(H) \cap \Phi(\eB)$ and so $\eL \supset \eL_0$.
Therefore, $\eL_0 = \eL$ and
\[
  \rank_{\eA}(\Lambda_H) = \rank_{\eA}(\eC).
\]
{From} \eqref{E:Epow} we have $\rank_{\eA}(\Lambda) = \ell
\rank_{\eA}(\eC)$, and of course $\rank_{\eA}(\Lambda) =
\rank_{\eA}(\eB) = m \rank_{\eA}(\eC)$.  Hence $\ell = m$ and $E$ is
isogenous to $H^m$.

For part (b) we prove the contrapositive of both directions.  If $E$
is not simple, it contains a proper simple sub-$t$-module $H$.
Picking up the considerations and notation of part (a), we see that
$d_0 = d/m$.  As the action of $\eC$ on $\Lambda_H$ is obtained
through $\bss_{\cS}$, we find from part (a) that we can partition
\[
  \{ \sigma_1, \dots, \sigma_d \} = \bigcup_{i = 1}^m \{ \sigma_{i_1},
  \dots, \sigma_{i_{d_0}} \}
\]
such that for all $1 \leq i, j \leq m$ and for all $c \in \eC$,
\[
  ( \sigma_{i_1}(c), \dots, \sigma_{i_{d_0}}(c)) =
  ( \sigma_{j_1}(c), \dots, \sigma_{j_{d_0}}(c)).
\]
Thus for each $k$, $1 \leq k \leq d_0$, the embeddings $\{
\sigma_{i_k} \}_{i=1}^m$ agree on $\eL$.  Since $m = [\eK:\eL]$, the
set $\{ \sigma_{i_k} \}_{i=1}^m$ can be written as $\cS_\alpha \assign
\Emb(\alpha(\eK)/\alpha(\eL))\circ \alpha$, for $\alpha$ equal to any
of the elements $\sigma_{i_k}$, $ i = 1, \dots, m$.  Thus $\cS$ is the
union of sets $\cS_\alpha$, and the criteria of Lemma~\ref{L:crit} are
satisfied.

Suppose now that there is a field $\eL \subset \eK$ satisfying the
equivalent criteria in Lemma~\ref{L:crit}.  By Remark~\ref{R:CMisog}
it suffices to consider the case that the period lattice of $E$ is
isomorphic to the maximal order $\eB$ of $\eK$.  Let $d_0 \assign [
\eL_+ : \ek ]$, $m \assign [\eK:\eL]$, and let $\eC$ and $\eC_+$ be
the integral closures of $\eA$ in $\eL$ and $\eL_+$ respectively.  By
hypothesis, we know that $\cS$ is the union of $d_0$ subsets of
$\Emb(\eK/\ek)$, each consisting of all $[\eK:\eL]$ embeddings of
$\eK$ with given restriction to $\eL$, given by, say,
$\sigma_1|_{\eL},\dots, \sigma_{d_0}|_{\eL}$, i.e.\
\begin{equation} \label{E:betaC}
\cS = \cS_{\sigma_1} \cup \cdots \cup \cS_{\sigma_{d_0}}.
\end{equation}
We may index the elements of $\cS = \{ \sigma_i \}$ so that
\begin{equation} \label{E:sigm}
\sigma_{(j-1)d_0 + k} \in \cS_{\sigma_k}, \quad 1 \leq j \leq m,
\quad 1 \leq k \leq d_0.
\end{equation}
We define $\Lambda_{\eC}$ to be the image in $\C^{d_0}$ of $\eC$ under
the conjugate embedding
\[
\bss_{\eL} \assign \sigma_1|_{\eL} \oplus \dots \oplus
\sigma_{d_0}|_{\eL},
\]
as in \eqref{E:conj} and \eqref{E:BLmb}, taking $\bsl_{\eC} \assign
(1, \dots, 1)$.  By Anderson's Theorem~\ref{T:HBD}, $\Lambda_{\eC}$ is
the period lattice of a $t$-module $H_{\eC} = (\Ga^{d_0},\Phi_{\eC})$
of CM-type $(\eL,\{ \sigma_1|_{\eL}, \dots, \sigma_{d_0}|_{\eL}\})$
with complex multiplications by $\eC$.  Let $\Exp_{\eC}: \C^{d_0} \to
\C^{d_0}$ be the exponential function of $H_{\eC}$.

Let $\bsl = (\lambda_1, \dots, \lambda_d)$ be given as in
\eqref{E:BLmb} so that $\Lambda = \bss_{\cS}(\eB)(\bsl)$, and let $\eM
\assign \beta_1\eC + \dots + \beta_m \eC \subset \eB$ be a direct sum
so that $\eM$ is a free $\eC$-module of maximal rank $m$ lying inside
$\eB$.  Let $\Lambda_{\eM} \assign \bss_{\cS}(\eM)(\bsl)$ which has
finite index in $\Lambda$.  Then for $\eb = \sum_{i=1}^m \beta_i\ec_i
\in \eM$, by definition we have
\begin{equation} \label{E:bss}
\begin{aligned}
\bss_{\cS}(\eb)(\bsl) &= \sum_{i=1}^m \bss_{\cS}(\beta_i \ec_i)(\bsl) \\
&= \sum_{i=1}^m \begin{pmatrix} \bss_{\eL}(\ec_i) & \cdots & 0 \\
        \vdots & \ddots & \vdots \\ 0 & \cdots & \bss_{\eL}(\ec_i)
        \end{pmatrix} \begin{pmatrix} \sigma_1(\beta_i) \lambda_1 \\
        \vdots \\ \sigma_d(\beta_i) \lambda_d \end{pmatrix}.
\end{aligned}
\end{equation}
Keeping the convention \eqref{E:sigm} in mind, for $1 \leq i, j \leq m$
we let $U_{ji} \in \Mat_{d_0 \times d_0}(\C)$ be the matrix
\[
U_{ji} \assign
\begin{pmatrix}
\sigma_{(j-1)d_0 + 1}(\beta_i) \lambda_{(j-1)d_0 + 1} & \cdots & 0 \\
\vdots & \ddots & \vdots \\
0 & \cdots & \sigma_{jd_0}(\beta_i)\lambda_{jd_0}
\end{pmatrix},
\]
and then set
\[
U \assign (U_{ji}) \in \Mat_{d \times d}(\C).
\]
{From} \eqref{E:bss} we see that
\[
\Lambda_{\eM}= \bss_{\cS}(\eB)(\bsl) = U (\sigma_1(c_1), \dots,
\sigma_{d_0}(c_1), \dots, \sigma_1(c_m), \dots, \sigma_{d_0}(c_m))^{tr},
\]
where $c_1, \dots, c_m$ run through $\eC$.  Since the image of $U$
spans $\Lie(E)$, we know that $U$ is invertible.  It follows that for
$z_1, \dots, z_d \in \C$,
\[
\bz \in \Lambda_{\eM} \Longleftrightarrow U^{-1}\bz \in(\sigma_1(c_1),
\dots, \sigma_{d_0}(c_1), \dots, \sigma_1(c_m), \dots,
\sigma_{d_0}(c_m))^{tr}, 
\]
where $\bz = (z_1, \dots, z_d)^{tr}$.  Let $\pi_i$ denote projection
onto the $i$-th component of $\bss_{\eL}(\eC) \oplus \cdots \oplus
\bss_{\eL}(\eC)$ so that
\begin{equation} \label{E:elin}
I_d = U \begin{pmatrix}
        \pi_1 & & 0 \\ & \ddots & \\ 0 & & \pi_{d_0}
        \end{pmatrix} U^{-1}.
\end{equation}
Now define an analytic map $e_{\eM} : \C^d \to \C^d$ by
\begin{equation} \label{E:edef}
e_{\eM}(\bz) \assign U \, \begin{pmatrix}
\Exp_{\eC}(\pi_1 U^{-1}(\bz)) \\ \vdots \\
\Exp_{\eC}(\pi_{d_0} U^{-1}(\bz))
\end{pmatrix}.
\end{equation}
The function $e_{\eM}(\bz)$ then (1) is entire; (2) vanishes exactly
on $\Lambda$ with simple zeros; (3) as a power series in $z_1,
\dots, z_d$ satisfies $\partial e_{\eM}(\bz) = I_d$ according to
\eqref{E:elin}; and (4) inherits a functional equation from
$\Exp_{\eC}$.  These properties make $e_{\eM}(\bz)$ the exponential
function of the $t$-module $E_{\eM} = (\Ga^d,\Phi_{\eM})$ where
$\Phi_{\eM} : \eA \to \Mat_{d \times d} (\Ct)$ is given by
\[
  \Phi(a) = U \begin{pmatrix}
    \Phi_{\eC}(a) & & 0 \\ & \ddots & \\ 0 & & \Phi_{\eC}(a)
    \end{pmatrix} U^{-1}.
\]
As $(\pi_i U^{-1}(U \bz))_i = (\bz_1, \dots, \bz_m)$ for arbitrary
$\bz = (\bz_1, \dots, \bz_m) \in (\C^{d_0})^m$, by \eqref{E:edef} this
$t$-module is uniformizable, since $H_{\eC}$ is.  Moreover, it is
abelian since $H_{\eC}$ is.  By construction $E_{\eM}$ is isomorphic
to $H_{\eC}^m$.

Now $\Lambda_{\eM} \subset \Lambda$ has finite index.  By
Lemma~\ref{L:Elatt} there is a uniformizable abelian $t$-module
$E_{\Lambda}$ isogenous to $E_{\eM}$, which has $\Lambda$ as its
period lattice.  However, uniformizable abelian $t$-modules are
determined by their period lattices by a theorem of Anderson
\cite{And86}, Cor.~2.12.2, and so $E_{\Lambda} = E$.  Therefore,
$E$ is isogenous to $H_{\eC}^m$ and so contains a sub-$t$-module
isogenous to $H_{\eC}$.
\end{proof}

\begin{remark}
The identity \eqref{E:edef} explicitly determines the way $E_{\eM}$
decomposes as a product $H_{\eC}^m$.  In practice, if $E = E_{\eM}$,
this is a fruitful method for determining the exponential functions of
such non-simple $t$-modules of CM-type, as shown by the examples in
Section~\ref{S:exmp}.
\end{remark}

\begin{theorem} \label{T:subt}
  Let $E_1$ and $E_2$ be $t$-modules of CM-types $(\eK_1,\cS_1)$ and
  $(\eK_2,\cS_2)$ with maximal real subfields which are Galois over
  $\ek$.  There exists a non-zero $t$-module morphism $E_1 \to E_2$ if
  and only if there are subfields $\eL_1 \subset \eK_1$ and $\eL_2
  \subset \eK_2$ and a $\ek$-isomorphism $\rho \colon \eL_1 \to \eL_2$
  such that each $\eL_i$ satisfies the criteria of Lemma~\ref{L:crit}
  for $\eK_i$ and $\cS_1|_{\eL_1} = (\cS_2|_{\eL_2})\circ \rho$.
\end{theorem}

\begin{proof}
  Suppose $E_1$ is isogenous to $H_1^{m_1}$ and $E_2$ is isogenous to
  $H_2^{m_2}$ with $H_1$ and $H_2$ both simple.  From
  Theorem~\ref{T:Esub}a it follows that the group $\Hom(E_1,E_2)$ of
  $t$-module morphisms is non-trivial if and only if $H_1$ and $H_2$
  are isogenous.  By Remark \ref{R:subCM} we know that $H_1$ and $H_2$
  are $t$-modules of CM-type.  Because for the concerns of this
  theorem we fix our sub-$t$-modules only up to isogeny, by
  Remark~\ref{R:CMisog} we can assume that $H_1$ and $H_2$ have period
  lattices isomorphic to the full rings of integers in their
  respective CM-fields.

Now suppose that subfields $\eL_1 \subset \eK_1, \eL_2 \subset \eK_2$
exist as in the statement of the theorem.  Let $\cS_1|_{\eL_1} = \{
\sigma_1, \dots, \sigma_d \}$ and $\cS_2|_{\eL_2} = \{ \tau_1, \dots,
\tau_d \}$. Because $\cS_1|_{\eL_1} = (\cS_2|_{\eL_2})\circ\rho$, we can
reorder the coordinates of either $H_1$ or $H_2$ so that
\begin{equation} \label{E:srt}
\sigma_i|_{\eL_1} = (\tau_i|_{\eL_2})\circ\rho, \quad 1 \le i \le d.
\end{equation}
Let $\eC_i$ be the ring of integers in $\eL_i$, $i = 1,2$.  We can
assume without loss of generality that for, say, $\bsl =
(1,\dots,1)^{tr} \in \C^d$ the period lattices $\Lambda_1$ and
$\Lambda_2$ of $H_1$ and $H_2$ are
\[
\Lambda_1 = \bss_{\cS_1}|_{\eL_1}(\eC_1)(\bsl) \quad \text{and} \quad
\Lambda_2 = \bss_{\cS_2}|_{\eL_2}(\eC_2)(\bsl).
\]
By \eqref{E:srt}, for $c \in \eC_1$, $\bss_{\cS_1|_{\eL_1}}(c) =
\bss_{\cS_2}|_{\eL_2}\circ \rho(c)$.  So it follows that $\Lambda_1 =
\Lambda_2$, and, as uniformizable abelian $t$-modules are determined
by their period lattices, the original $H_1$ and $H_2$ are isogenous
(by \cite{And86}, Cor.~2.12.2).

For the other direction, we start with a given isogeny $\psi \colon
H_1 \to H_2$, where each $H_i$ has CM by the full ring of integers
$\eC_i$ in $\eL_i$.  By Lemma 1.1 of \cite{Yu97}, there is an $a \in
\eA$ and an isogeny $\phi \colon H_2 \to H_1$ so that $\phi \circ \psi
= \Phi_1(a)$.  These isogenies induce a map of $\eA$-modules
$\End(H_1) \to \End(H_2)$ via
\[
f \longmapsto \psi \circ f \circ \phi.
\]
This map lifts to a $\ek$-linear map $\psi_\ast \colon \End^0(H_1) \to
\End^0(H_2 )$,
\[
\psi_\ast \colon f \otimes b \longmapsto \psi \circ f \circ \phi
\otimes \frac{b}{a},
\]
which is also multiplicative, as the following calculation shows:
\begin{align*}
  ((\psi \circ f \circ \phi) \otimes \frac{b}{a})\cdot ((\psi \circ g
  \circ \phi) \otimes \frac{c}{a})  %
  & = \psi \circ f \circ  \Phi_1(a) \circ g \circ \phi \otimes \frac{bc}{a^2} \\
&   = \psi \circ f \circ g \circ \phi \circ \Phi_2(a) \otimes
  \frac{bc}{a^2} %
= \psi \circ f \circ g \circ \phi \otimes \frac{bc}{a}.
\end{align*}

By Corollary \ref{C:EHBD}, since the $H_i$ are simple, the maps
$\Phi_i$ extend to isomorphisms $\Phi_i \colon \eL_i
\stackrel{\sim}{\longrightarrow} \End^0(H_i)$.  Thus
\[
\rho \assign \Phi_2^{-1} \circ \psi_{\ast} \circ \Phi_1 \colon \eL_1
\stackrel{\sim}{\longrightarrow} \eL_2.
\]
Note also that the natural map $d\colon \End(H_i) \longrightarrow
\End(\Lie(H_i))$ extends naturally to $\End^0(H_i)$ via
\[
d( f \otimes b) = (df)(b I_{d_i}).
\]

For each $c \in \eL_1$, we verify that the following diagram commutes:
\begin{equation} \label{E:diagrm}
\xymatrix{ \Lie(H_1) \ar[rr]^{d{\psi}}
        \ar[d]_{d\Phi_1(\ec)\,=\,\bss_{\cS_1}|_{\eL_1} (\ec)}
        & & \Lie(H_2)
        \ar[d]^{d{\psi}_{\ast}(\Phi_1(\ec)))} \\
        \Lie(H_1) \ar[rr]_{d{\psi}}
        & & \Lie(H_2)}
\end{equation}
Indeed, writing $c = b/e$ with $b \in \eC_1$ and $e \in \eA$, we see
that
\begin{align*}
d(\psi_{\ast}\circ\Phi_1(c))d\psi & = d ((\psi \circ \Phi_1(b)\circ
\phi)\otimes \frac{1}{ea}) d\psi \\
& = d\psi\ \bss_{\cS_1}|_{\eL_1}(b)\ d\phi\ d\psi\ \frac{1}{ea}{I_{d_1}} 
 = d\psi\ \bss_{\cS_1}|_{\eL_1}(c), 
\end{align*}
as $d\phi\, d\psi = a I_{d_1}.$ Because $\psi_{\ast} \circ \Phi_1 =
\Phi_2 \circ \rho$,
\[
d(\psi_{\ast}\circ\Phi_1(c)) = d(\Phi_2(\rho(c)) =
\bss_{\cS_2}|_{\eL_2}(\rho(c)). 
\]

Therefore it follows from the diagram \ref{E:diagrm} that, for every
$c \in \eL_1$, $\bss_{\cS_1}|_{\eL_1}(c)$ and
$\bss_{\cS_2}|_{\eL_2}(\rho(c))$ are conjugate and hence have the same
eigenvalues.  By choosing a primitive element $c_0 \in \eL_1$ over
$\ek$, we see that, up to a fixed permutation of coordinates,
$\bss_{\cS_1}|_{\eL_1}(c)$ equals $\bss_{\cS_2}|_{\eL_2}(\rho(c))$ for
all $c \in \eL_1$, and therefore ${\cS_1}|_{\eL_1} =
{\cS_2}|_{\eL_2} \circ \rho$.
\end{proof}

% qpbody.tex

\section{Biderivations and Quasi-Periodic Extensions of $t$-Modules}
\label{S:qp}

The theory of biderivations and quasi-periodic extensions for Drinfeld
modules was developed by P.~Deligne, Anderson, and Yu.  E.-U.~Gekeler
gave an alternate approach to the de Rham isomorphism in the very
accessible source \cite{Gek89}.  Anderson suggested the relevance of
quasi-periods for obtaining Gamma values at those rational points
which are not the ratio of monic polynomials to Thakur, who passed the
hint along to us some years later.  Here we extend much of
\cite{wdb93} {from} the setting of $\eA$-Drinfeld modules to that of
arbitrary $t$-modules.

For the remainder of Section~\ref{S:qp} we fix a $t$-module $E =
(\Phi, \Ga^d)$.

\subsection{Biderivations} \label{SS:bider}
A \emph{$\Phi$-biderivation} is an $\Fq$-linear map $\bsd : \eA
\to \tau M(E)$ satisfying the product formula that, for all
$\ea$, $\eb \in \eA$,
\[
\bsd(\ea\eb) = \iota(\ea)\cdot\bsd(\eb) + \bsd(\ea) \Phi(\eb).
\]
 
\begin{lemma} \label{L:sr} 
  Let $m \in \tau M(E)$, i.e.\ $m$ is any element of the $t$-motive
  $M(E)$ associated to $E$ such that $m$ has no $\tau^0$ terms:
\[
m \in \Hom^q_{L}(E,\Ga), \quad dm = 0.
\]
Then the assignment $t \mapsto \bsd(t) \assign m$ induces a
$\Phi$-biderivation $\bsd_m$.
\end{lemma}

\begin{proof} 
  By $\Fq$-linearity, we need only check that the
  ``competing'' expressions given by applying the product formula
  inductively to different factorizations
\[
t^{m_1}t^{n_1} = t^{m_2}t^{n_2},\ m_i,n_i \in \ZZ_{>0}, 
\]
are equal.  This verification is straightforward.
\end{proof} 

This lemma gives a natural isomorphism between the $L$-vector spaces
$\Der(\Phi)$ and $\tau M(E)$.  Certain $\Phi$-biderivations can be
given algebraically in terms of $\Phi$ using the identification $M(E)
\simeq (L\{\tau\})^d$.  Set
\[
N^{\perp} \assign N^\perp(L) \assign \{ V \in \Mat_{1\times
d}(L) : V N = 0\},
\]
where $N$ is the nilpotent part of $d\Phi(t) = \theta I_d + N$.

Let $\bsU = (U_1,\dots,U_d) \in (\Lt)^d$ with $d\bsU \in N^\perp(L),$
where $d\bsU = (dU_1,\dots,dU_d)$ denotes the vector of coefficients
of $\tau^0$ in $\bsU$. We define $\bsd^{(\bsU)} : \eA \to M(E)$ via
\begin{equation} \label{E:deltaU}
\bsd^{(\bsU)}(\ea) \assign \bsU\Phi(\ea) - \iota(\ea) \bsU,
\end{equation}
for every $a \in \eA$.  The condition $d\bsU \in N^{\perp}$ is
equivalent to saying that $d\bsd^{(\bsU)}(t) = 0$, i.e.\ that
$\bsd^{(\bsU)}(t) \in \tau M(E)$.

Since $\bsd^{(\bsU)}(\ea\eb) = \bsU\Phi(\ea\eb) -\iota(\ea\eb) \bsU =
\iota(\ea)(\bsU\Phi(\eb) - \iota(\eb) \bsU) + (\bsU\Phi(\ea) -
\iota(\ea)\bsU)\Phi(\eb)$, $\bsd^{(\bsU)}$ is indeed a
$\Phi$-biderivation.  Such $\Phi$-biderivations will be called
\emph{inner}, and they constitute an $L$-vector space which we denote
$\Der_{in}(\Phi)$.  Note further that, in terms of the $t$-motive
$M(E)$, we are setting $\bsd^{(\bsU)}(t) = (t-\theta)\bsU$.

\begin{lemma} \label{L:M_in} If $M = M(E)$ is a torsion-free
  $t$-motive over $L$, then, as $L$-vector spaces,
\[
M_{in} \assign (t - \theta)(\tau M + N^\perp \tau^0) \simeq
\Der_{in}(\Phi)
\]
via the natural isomorphism
\[
(t-\theta) \bsU \mapsto \bsd^{(\bsU)}.
\]
\end{lemma}

\begin{proof}
It is clear that for $\bsU \in M$, 
\[
(t - \theta)\bsU \in \Der_{in}(\Phi)
\Longleftrightarrow \bsU \in \tau M + N^\perp \tau^0. 
\]
Since $M$ is torsion-free, multiplication by $t - \theta$ is
injective.
\end{proof}

\begin{definition}
  Several other distinguished subspaces of $\Der(\Phi)$ and quotient
  spaces will play a role in our discussion:
\[
\begin{gathered}
\Der_0(\Phi) \assign \{ \bsd^{(\bsU)} : \bsU \in N^\perp(L) \tau^0 \} =
  \{ \bsd^{(\bsU)} \in \Der_{in}(\Phi) : \bsU \in L^d \tau^0 \} \\
\Der_{si}(\Phi) \assign \{ \bsd^{(\bsU)} : \bsU \in \tau M \} \qquad
  \text{(strictly inner)}\\
H_{DR}(\Phi) \assign \Der(\Phi)/\Der_{si}(\Phi) \qquad \text{(de Rham)}\\
H_{sr}(\Phi) \assign \Der(\Phi)/\Der_{in}(\Phi) \qquad \text{(strictly
  reduced)}
\end{gathered}
\]
The phrases after the definitions indicate the names for the type of
biderivation involved and thus account for the subscripts employed.
\end{definition}

\begin{proposition} \label{P:qpdim}
 For any $t$-motive $M = M(E)$ over $L$, where $E = (\Phi,\Ga^d)$,
\begin{gather*}
  \Der_{in}(\Phi) = \Der_0(\Phi) \oplus \Der_{si}(\Phi)  \\
 H_{DR}(\Phi) \simeq \Der_0(\Phi) \oplus H_{sr}(\Phi)
\end{gather*}
and, if $M$ is abelian,
\begin{gather}
\dim_{L} \Der_0(\Phi) = d - \rank N \\
\dim_L H_{DR}(\Phi) = r, \\
\dim_{L} H_{sr}(\Phi) = r - d + \rank N 
\end{gather}
where $d = \dim E = \rank_{L\{\tau\} }M$, $r = \rank E = \rank_{L[t]}
M,$ and $d\Phi(t) = \theta I_d + N$. 
\end{proposition}

\begin{proof}
  The direct sum decompositions are immediate.  By definition, $\dim_L
  \Der_0 (\Phi) = \dim_L N^{\perp} = d - \rank N$.  To calculate the
  dimension of $H_{sr}(\Phi)$, we appeal to the following diagram,
  where here we assume that $E$ is abelian of rank~$r$:
\[
\vcenter{\xymatrix@H+2pc  @+2pc % @-1.2pc  This sets the scaling. 
@ur
{{M_{in} } \ar@{^{(}->}[r]^{i_d} \ar@{^{(}->}[d]_{i_r} & 
{ \tau M } \ar@<-1ex>@{^{(}->}[d]^d \\
(t -\theta)M  \ar@{^{(}->}[r]_{r} & {  M\quad} \
}}
\]

The labels indicate the $L$-dimensions of the quotients.  
By Lemma \ref{L:M_in}, $i_r$ is the codimension of $\tau M +
N^{\perp}\tau^0$ in $M$.  Thus $i_r = \rank N$ and $\dim_L
H_{sr}(\Phi) = i_d = r - d + \rank N$.
\end{proof}

\subsection{Quasi-Periodic Functions} \label{SS:qpf}

In this section we investigate quasi-periodic functions and
quasi-periods associated to $\Phi$-biderivations.  The reader is
directed to Gekeler~\cite{Gek89} for a historical motivation for this
and related terminology.

\begin{proposition} \label{P:qpfn} Given a $\Phi$-biderivation $\bsd$
  defined over $L$, there is an entire $\Fq$-linear function $F_{\bsd}
  : \C^d \to \C$ given by the unique power series satisfying:
\begin{gather}
F_{\bsd}(d\Phi(\ea)\bz) = \iota(\ea) F_{\bsd}(\bz) +
\bsd(\ea)\Exp(\bz), \label{E:Fdeq} \\
F_{\bsd}(\bz) \equiv 0 \pmod {\bz^q}, \label{E:Fdlin}
\end{gather}
where the latter condition means that every non-zero monomial in the
power series $F_{\bsd}(\bz)$ is of the form $c_{hi}z_i^{q^h}$ with $h
> 0$, $c_{hi} \in \C$.  Furthermore, $F_{\bsd}(\bz)$ has coefficients
{from} $L$.
\end{proposition}

\begin{proof}  Write $\Exp(\bz) = (e_1(\bz),\dots,e_d(\bz)),$
  $\bsd(t) = (\delta_1(t),\dots,\delta_d(t))$, with the variables chosen
  so that $N = (n_{ij})$ is upper triangular.  If 
\[
F_{\bsd}(\bz) =
  \sum_h {\mathbf c}_h\cdot \bz^{q^h} =
\sum_h c_{h1}z_1^{q^h} + \dots + c_{hd}z_d^{q^h},
\] 
then we can equate coefficients in \eqref{E:Fdeq} with $a = t$ to find
the equality
\[
(\theta^{q^h} - \theta) c_{h1} z_1^{q^h} = \text{\rm Term involving}\ 
z_1^{q^h}\ \text{\rm in}\ \sum \delta_i(t)e_i(\bz).
\]
It shows that the $c_{h1}$ is uniquely determined and that the terms
$c_{h1} z_1^{q^h} \to 0$ as $h \to \infty$, for any fixed value of
$z_1$.  Therefore, for any fixed values of $z_2$, the terms
$c_{h1}(z_2n_{12})^{q^h}$ tend to zero as $h \to \infty$ in the
following equality:
\[
(\theta^{h} - \theta) c_{h2} z^{q^h} + c_{h1}n_{12}^{q^h} z_2^{q^h} =
\text{\rm Term involving}\ z_2^{q^h}\ \text{\rm in}\ \sum
\bsd_i(t)e_i(\bz).
\]

Consequently, we find that the terms $c_{h2}z_2^{q^h}$ tend toward
zero as $h \to \infty$ for any fixed value of $z_2$.  Iterating this
argument, we find that each series
\[
\sum_h c_{hi}z_i^{q^h}
\]
is uniquely determined and entire.  Thus $F_{\bsd}(\bz)$ is uniquely
determined and everywhere convergent.

The $\Fq$-linearity of $F_{\bsd}(\bz)$ and the following recursive
calculation show that $F_{\bsd}(\bz)$ satisfies the required
functional equation with respect to $\bsd$:
\begin{align*}
F_{\bsd}(d\Phi(t^{i+1})\bz) &= \theta F_{\bsd}(d\Phi(t^{i})\bz) +
\bsd(t)\Exp(d\Phi(t^i) \bz) \\ 
&= \theta(\theta^i F_{\bsd}(\bz) + \bsd(t^i)\Exp(\bz)) +
\bsd(t)\Phi(t^i)\Exp(\bz) \\
&= \theta^{i+1}F_{\bsd}(\bz) + \bsd(t^{i+1})\Exp(\bz).
\end{align*}
\end{proof}

The unique $\Fq$-linear function given by Proposition \ref{P:qpfn} is
said to be the \emph{quasi-periodic function associated to $\bsd$.}

It is easy to check from the unicity of $F_{\bsd}$ that, if $\bsd =
\ell_1 \bsd_1 + \ell_2 \bsd_2$, with $\ell_1$, $\ell_2 \in L$, then
\[
F_{\bsd}(\bz) = \ell_1 F_{\bsd_1}(\bz) + \ell_2 F_{\bsd_2}(\bz).
\]

\begin{remark}
Note that, since the image of $\Exp(\bz)$ is dense in $\Ga^d(\C)$, if
the biderivation $\bsd$ is non-zero, then all functions satisfying the
functional equation~\eqref{E:Fdeq} of Proposition~\ref{P:qpfn} are
non-zero, even if we allow solutions which violate \eqref{E:Fdlin}.
The solutions to that relaxed functional equation comprise the set
$F_{\bsd}(\bz) + N^\perp(L) \cdot \bz$, where $F_{\bsd}(\bz)$ is the
quasi-periodic function of $\bsd$.
\end{remark}

\begin{proposition} \label{P:FU}
The quasi-periodic function related to $\bsd^{(\bsU)} \in
\Der_{in}(\Phi)$ is
\[
F^{(\bsU)}(\bz) \assign \bsU  \Exp(\bz) - d\bsU\bz.
\]
\end{proposition}

\begin{proof} Note that $F^{(\bsU)}(\bz)$ has no linear terms and
that the functional equation holds:
\begin{align*}
  F^{(\bsU)}(d\Phi(t)\bz) &= \bsU\Exp(d\Phi(t)\bz) - d\bsU d\Phi(t)\bz \\
  & = \theta(\bsU\Exp(\bz) - d\bsU \bz) + (\bsU \Phi(t) -
  \theta\bsU)\Exp(\bz).
\end{align*}
\end{proof}

If $\bsl \in \Lambda \assign \Ker \Exp$ is a period of of $E$ and
$\bsd \in \Der(\Phi)$, then $\eta \assign \eta_{\delta}(\bsl) \assign
F_{\bsd}(\bsl)$ is called the \emph{quasi-period of $\bsd$
  corresponding to $\bsl$.}  The following result is an immediate
corollary of Proposition~\ref{P:FU}.

\begin{corollary} \label{C:inqp}
The quasi-period $\eta_{\bsd^{(\bsU)}}(\bsl)$ associated to the
  inner biderivation $\bsd \assign \bsd^{(\bsU)}$ is
\[
\eta_{\bsd^{(\bsU)}}(\bsl) = d\bsU \cdot \bsl.
\]
\end{corollary}

Therefore the quasi-periods of inner biderivations defined over $L$
are $L$-linear combinations of the coordinates of the corresponding
periods.

\subsection{Quasi-Periodic Extensions} \label{SS:qpext} 

Let $\bsd_1,\dots,\bsd_j$ be $\Phi$-biderivations.  Then define the
matrix
\[
\Psi(t) \assign
\left(\begin{array}{c|c}
& \\ \Phi(t) & 0 \\
& \\ \hline
\\[-6pt]
\begin{gathered}
\bsd_1(t) \\ \vdots \\ \bsd_j(t)
\end{gathered}\ 
&
\quad \theta I_j \quad
\end{array}
\right),
\]
where $\bsd_i(t)$ is situated directly beneath $\Phi(t)$ in the
$(d+i)$th row of $\Psi(t)$, and the unique following non-zero entry in
that row is in the $(d+i)$th column, i.e.\ in the $i$th column
following the entries for $\bsd_i(t)$.  Here and below we omit the
$\tau^0$ from the linear terms for easier reading of matrices.  In
addition, define the entire mapping
\[
\Exp_\Psi : \C^{d+j} \to \C^{d+j}
\]
via
\[
\Exp_\Psi(\bz,\bu) \assign \left(\Exp(\bz), u_1 + F_1(\bz),\dots,
  u_j + F_j(\bz)\right)^{tr},
\]
where for simplicity we have written $F_i$ for $F_{\bsd_i}$.

\begin{proposition} \label{P:qpext} In the situation of the preceding
paragraph, $Q \assign (\Psi,\Ga^{d+j})$ is a $t$-module with
exponential function $\Exp_Q = \Exp_\Psi$ and with periods
\[
(\bsl, -\eta_1(\bsl), \dots,-\eta_j(\bsl)),
\]
where $\eta_i(\bsl)$ is the quasi-period of $\bsd_i$ corresponding
to the period $\bsl \in \C^d$.
\end{proposition}
 
When $\bsd_1,\dots,\bsd_j$ represent $L$-linearly independent classes in
$H_{sr}(\Phi)$, we call the corresponding extension $Q$ of $E$ a {\em
    strictly quasi-periodic extension}.

\begin{proof}
That $\Psi$ defines a $t$-module follows from the hypothesis that
$\Phi$ does and that
\[
d\Psi(t) = 
\begin{pmatrix}
d\Phi(t) & 0 \\
0 & \theta I_j
\end{pmatrix},
\]  
where $I_j$ is the $j \times j$ identity matrix.

That $\Exp_\Psi$ satisfies the appropriate functional equation 
\[
\Exp_\Psi (d\Psi(t)(\bz,\bu)) = \Psi(t) \Exp_\Psi(\bz,\bu)
\]
follows from the functional equations for $\Exp_\Phi$ and for the
$F_i$.  Moreover, since the $F_i(\bz)$ have zero linear terms, the
linear terms of $\Exp_\Psi(\bz,\bu)$ are precisely $(\bz,\bu)^{tr}$,
as required for the exponential function of a $t$-module.
\end{proof}

Thus we find that $\Psi$ gives an extension of the $t$-module $E$
by the basic $t$-module $\Ga^j$:
\[
\begin{CD}
0 @>>> {\Ga^j(\C)} @>>> Q @>>> E @>>> 0 \\
@. @A\text{id}AA @A{\Exp_\Psi}AA @A{\Exp_\Phi}AA \\
0 @>>> \Ga^j(\C) @>>> \Lie(Q) @>>> \Lie(E) @>>> 0.
\end{CD}
\]

This is a generalization of the one-dimensional (Drinfeld) case and an
exact analogue of the extensions of an elliptic curve ${\mathcal E}$
by the additive group $\Ga$.  In the latter situation we have the
exponential maps
\begin{alignat*}{2}
  \Exp_\Phi \quad &\longleftrightarrow & z &\mapsto (\wp(z),\wp^\prime(z),1) \\
  \Exp_\Psi \quad &\longleftrightarrow \quad & (u,z) &\mapsto 
  (u + b\zeta(z), \wp(z),\wp^\prime(z),1),
\end{alignat*}
involving the Weierstrass quasi-elliptic function $\zeta$ and a
constant $b \in {\mathbb C}$ (see \cite{wschmidt}, {\S}3.2.c).  The
various choices of the constant $b$ classify the possible extensions,
with $b = 0$ giving the trivial, i.e.\ split, extension.  Thus
$\Ext({\mathcal E},\Ga) \simeq {\mathbb C}$.

\begin{remark}
The quasi-periodic extension $Q$ is uniformizable if and only if $E$ is
uniformizable.  Moreover notice that $Q$ is never abelian for $j > 0$,
since the $\Psi(t)$ action on the latter $j$ coordinates is
multiplication by the scalar $\theta$.
\end{remark}

The {\it leitmotiv} of the remainder of this section is that
quasi-periodic extensions of $E$ of the same dimension depend
essentially only on the strictly reduced cohomology classes involved.
The first step is given by the following:

\begin{proposition} \label{P:srext} Let $\bsd_1, \dots, \bsd_j$ and
  $\bsd'_1, \dots, \bsd'_j$ be $\Phi$-biderivations.  If both sets
  generate the same subspace of $H_{sr}(\Phi)$, then the corresponding
  quasi-periodic extensions are isomorphic.
\end{proposition}

\begin{proof}
  This proposition follows from the following remark, whose
  verification is immediate.  The first identity there shows that the
  isomorphism class of the extension is independent of the
  representatives chosen for the classes in $H_{sr}(\Phi)$.  The
  second shows that the isomorphism class is independent of the
  generators chosen for the span in $H_{sr}(\Phi)$.
\end{proof}

\begin{lemma}\label{L:extlem} If $\bsd_1,\bsd_2$ are
  $\Phi$-biderivations and $\bsU \in M(E)$ with $d\bsU \in N^{\perp}$,
  then
\begin{gather*}
\begin{pmatrix}
I_d & 0_d \\
\bsU & 1
\end{pmatrix}
\begin{pmatrix}
\Phi(t)   & 0_d \\
\bsd_1(t) & \theta
\end{pmatrix} 
\begin{pmatrix}
I_d  & 0_d \\
-\bsU & 1
\end{pmatrix} =
\begin{pmatrix}
\Phi(t)                  & 0_d \\
\bsd_1(t) + \bsd^{(\bsU)}(t) & \theta
\end{pmatrix}
\\
\begin{pmatrix}
I_d      & 0_d & 0_d \\
0_d^{tr} &  1  & 0   \\
0_d^{tr} & c & 1   
\end{pmatrix}
\begin{pmatrix}
\Phi(t)    & 0_d & 0_d \\
\bsd_1(t)  & \theta \   & 0   \\
\bsd_2(t)  & 0   & \theta    
\end{pmatrix}
\begin{pmatrix}
I_d      & 0_d & 0_d \\
0_d^{tr} &  1  & 0   \\
0_d^{tr} & -c    & 1
\end{pmatrix} =
\begin{pmatrix}
\Phi(t)                 & 0_d & 0_d \\
\bsd_1(t)               & \theta    & 0 \\
\bsd_2(t) + c \bsd_1(t) & 0   & \theta 
\end{pmatrix},
\end{gather*} 
where $0_d$ denotes the zero (column) vector of length $d$.
\end{lemma}

One important consequence of this lemma is the following remark: 
\begin{corollary} \label{C:split} Let $Q$ be the extension of $E$
  associated to the $\Phi$-biderivation $\bsd$.  Then the extension
\[
0 \rightarrow \Ga \rightarrow Q \rightarrow E \rightarrow 0
\]
splits if and only if $\bsd$ is inner.
\end{corollary}

\begin{proof}
If $\bsd = \bsd^{(\bsU)}$ is inner, then choose $\bsd_1 = 0$ in the
first part of the preceding lemma.

Similarly, if there is a self-isogeny $\Theta$ of $Q$ such that
\[
\Theta \Psi(t) = 
\begin{pmatrix}
\Phi(t) & 0_d \\
0_d^{tr} & \theta
\end{pmatrix}
\Theta,
\] 
then comparing entries first in the lower right-hand corner gives that
the lower right-hand entry $c$ in $\Theta$ satisfies
$c \theta = \theta c$, i.e.\ $c \in L$.  Then considering entries 
along the bottom row shows that, up to a non-zero scalar
multiple $c$,  
\[
c\bsd(t) + \bsU \Phi(t) = \theta \bsU,
\]
where $(\bsU ; c\tau^0)$ is the bottom row of $\Theta$.
Now if $c = 0$, then $(t - \theta)\bsU = 0$.  However, since the $t$-motive
$M(Q)$ is a free $L[t]$-module, we would have $\bsU = 0$.  Thus the
bottom line of $\Theta$ would consist of zeros.  In that case,
$\Theta$ would have an infinite kernel and could not be an isogeny.
Therefore $c \ne 0$, and $\bsd$ is inner, as claimed.
\end{proof}

\subsection{Minimality of Quasi-Periodic Extensions} \label{SS:minqp}

Let $\Delta$ be a surjective morphism from the uniformizable
$t$-module $(\Upsilon,\Ga^e)$ to the uniformizable $t$-module
$(\Phi,\Ga^d)$.  We say that $(\Upsilon,\Ga^e)$ is a {\em minimal
extension} of $(\Phi,\Ga^d)$ if no proper sub-$t$-module of
$(\Upsilon,\Ga^e)$ surjects onto $(\Phi,\Ga^d)$.  Intuitively this
places $(\Upsilon,\Ga^e)$ at the opposite extreme from a split
extension.  That intuition will be made precise in this section. 

\begin{proposition} \label{P:minqp}
  Let $Q = (\Psi,\Ga^{d+j}) $ be the extension associated to
  $\Phi$-biderivations $\bsd_1,\dots,\bsd_j$.  Then the following are
  equivalent:
\begin{enumerate}
\item[(a)] The biderivations $\bsd_1,\dots,\bsd_j$ represent linearly
  independent classes in $H_{sr}(\Phi)$.
\item[(b)] $Q$ is a minimal extension of $E$.
\end{enumerate}
\end{proposition}

\begin{proof}  (a) $\Rightarrow$ (b).
  We adopt the notation of Section \ref{SS:qpext} for our strictly
  reduced quasi-periodic extension.  The plan is to start with a
  non-trivial algebraic relation on a proper sub-$t$-module $H$ of $Q$
  and conclude that the image of $H$ lies in a proper sub-$t$-module
  of $E$.  Choose coordinates for $E$ so that $d\Phi(t)$ is
  upper-triangular.
  
  Assume that we have a non-trivial relation holding on the
  coordinates of $H$, which is smallest with respect to the reverse
  lexicographical ordering on monomials in $x_1, \dots, x_d, u_1,
  \dots, u_j$.  Since the underlying group is $\Fq$-linear, this
  minimal relation has the following form for all $(\bx,\bu) \in H$
  (see \cite{lang}, {\S}VI.12):
\begin{equation} \label{E:reln}
R(\bx,\bu) = R_1(x_1) + \dots + R_d(x_d) + S_1(u_1) + \dots + S_j(u_j)
= 0,
\end{equation}
in which all the $R_i$, $S_i \in L\{\tau\}$.  Now we make a series of
observations based on the fact that $H$ is a sub-$t$-module, i.e.\
that, for all $(\bx,\bu) \in H$,
\begin{equation} \label{E:nreln}
R \circ \Psi(t)(\bx,\bu) = R_1^\prime(x_1) + \dots + R_d^\prime(x_d) +
S_1^\prime(u_1) + \dots + S_j^\prime(u_j) = 0
\end{equation}
as well.  We show by contradiction that the variables of $\bu$ are not
involved in this minimal relation.

Since the relation \eqref{E:reln} is minimal and the effect of
$\Psi(t)$ on the variables $\bu$ is multiplication by $\theta$ (plus a
sum affecting the variables of $\bx$), we apply $\Psi(t)$ to obtain
another algebraic relation and, on comparing maximal monomials with
respect to our ordering, we conclude that, if the variables $\bu$ were
involved in \eqref{E:reln}, then for each $i$
\[
S_i(\theta u_i) = \theta^{q^h} S_i(u_i), 
\]
for a fixed $h$.  Thus, since $\theta \notin \overline{\mathbb{F}}_q$,
each $S_i$ is a monomial of degree $q^h$.  However $L$ is a perfect
field, and we can write
\begin{equation} \label{E:Ssum}
S_1(u_1) + \dots + S_j(u_j) = (s_1 u_1 + \dots + s_j u_j)^{q^h},
\end{equation}
with the $s_i \in L$.

Again comparing terms in the relations \eqref{E:reln}, \eqref{E:nreln}
and \eqref{E:Ssum}, but this time for the variables in $\bx$, we find
that
\begin{equation}\label{E:rreln}
\theta^{q^h}\mathbf{R}\bx = \mathbf{R}\Phi(t)\bx + \tau^h\bsd(t)\bx,
\end{equation}
where $\mathbf{R} = (R_1,\dots,R_d)$ and $\bsd \assign s_1\bsd_1 +
\dots + s_j\bsd_j$.  In particular, for each variable $x_{\ell}$,
\[
\theta^{q^h} R_{\ell} x_{\ell} = \sum R_i\phi_{i,\ell}(t)x_{\ell}
+ s_i^{q^h} \tau^h \delta_{i,\ell}(t) x_{\ell},
\]
where $\Phi(t) = (\phi_{i,\ell}(t))$ and $\bsd_i(t) =
(\delta_{i,1}(t),\dots,\delta_{i,j}(t)) \in \tau M$.

Assume for the moment that $h \neq 0$.  Let $R_{\ell} = r_{\ell}
\tau^0 + \text{higher degree terms}$, $r_{\ell} \in L$.  Notice now
that the $\delta_{i,\ell}(t)$ lack linear terms. Since $d\Phi(t) =
\theta I_d + N$ is upper-triangular, $\phi_{i,1}(t)$ has no linear
terms unless $i=1$.  From \eqref{E:rreln} we see that
\[
\theta^{q^h}r_1 = r_1\theta.
\]
Because $h \neq 0$, it follows that $r_1 = 0$, i.e.\ $R_1$ does not
have a linear term.  Similarly from \eqref{E:rreln} we see that
\[
\theta^{q^h}r_2 = r_1\,d\phi_{1,2}(t) + r_2\theta.
\]
Since $r_1=0$, we see as above that $r_2=0$.  Proceeding by induction
we find that $r_{\ell} = 0$ for $\ell = 1, \dots, d$, and so none of
the $R_{\ell}$ involve linear terms.  But in that case, since the
$R_{\ell}$ are $\Fq$-linear, the relation \eqref{E:reln} would not be
minimal; we could extract the $q$th root and have a equation of
smaller degree for elements of $H$.

Therefore we are reduced to the case $h = 0$, and we find that
\eqref{E:rreln} has the form
\[
\theta{\mathbf R} = {\mathbf R}\Phi(t) + \bsd(t).
\]
However this means that $- \bsd$ and thus $\bsd$ are inner, whereas
the given $\bsd_1,\dots,\bsd_j$ are $L$-linearly independent modulo
the inner biderivations. Consequently $\bsd = 0$.  Since $s_1 = \dots
= s_j = 0$, we see that the minimal relation \eqref{E:reln} actually
involves only variables from $\bx$ after all.  Therefore the
projection of $H$ in $E$ cannot be surjective.

(b) $\Rightarrow$ (a): Assume now that $\bsd_1,\dots, \bsd_j$ do not
represent linearly independent classes in $H_{sr}(\Phi)$.  Then by an
automorphism of $Q$ involving only a linear change of coordinates, we
may assume that $\bsd_j$ is inner.  Then, as we have seen above, we
may conjugate by a matrix leaving the terms corresponding to
$\bsd_1,\dots,\bsd_{j-1}$ invariant to see that $Q$ has a direct factor
corresponding to $\bsd_j$.  If there are further inner biderivations
lying in $L\bsd_1 + \dots + L\bsd_{j-1}$, we may proceed to find
further direct factors of $Q$.

At any rate, in this case, the direct complement of these factors form
    a proper sub-$t$-module of $Q$ which projects onto $E$.  Thus $Q$
    is not a minimal extension of $E$ in this case.
\end{proof}

Therefore we have the following result:

\begin{corollary} Let the quasi-periodic extension
\[
0 \rightarrow \Ga^j \rightarrow Q \overset{\pi}\longrightarrow E
\rightarrow 0.
\]
be associated to a basis for the subspace $D \subset \Der(\Phi)$.
Then there is a unique maximal direct factor of $Q$ lying in $\Ker
\pi \simeq \Ga^j$, and it has dimension equal to the dimension of
$D \cap \Der_{in}(\Phi)$
\end{corollary}

\begin{proof}
  According to the identity of the first part of Lemma \ref{L:extlem},
  the elements of $\Ker \pi$ corresponding to $D_{in} \assign D \cap
  \Der_{in}(\Phi)$ form a direct summand of $Q$.  Let $H \subset \Ker
  \pi$ be a direct summand of $Q$.  We need to show that the space $B$
  of biderivations associated to elements of $H$ is contained in
  $D_{in}$.  The above short exact sequence becomes
\[
  0 \rightarrow H\oplus K^\prime \rightarrow H\oplus Q^\prime
  \overset{\pi}\longrightarrow E \rightarrow 0
\]
or, since $H \subset \Ker \pi$,
\[ 
  0 \rightarrow K^\prime \rightarrow Q^\prime
  \overset{\pi|_{Q^\prime}}\longrightarrow E \rightarrow 0.
\]

Choose an $L$-basis $\bsd_1,\dots,\bsd_i$ for $B\cap D_{in}$ and extend
to a basis for $D_{in}$ with the biderivations
$\bsd_{i+1},\dots,\bsd_l$.  Augment $\bsd_1,\dots,\bsd_i$ with
$\bsd_{l+1},\dots,\bsd_m$ to an $L$-basis for $B$ and finally extend
with biderivations $\bsd_{m+1},\dots,\bsd_j$ to obtain a basis
$\bsd_1,\dots,\bsd_j$ which produces $\Ker \pi$.

 Now the part of $H$ corresponding to $B\cap D_{in}$ is a direct
summand of $H$.  Quotienting the original exact sequence by it, we may
assume that $i = 0$, i.e.\ that $B\cap D_{in} = 0$, and take as our
objective to show that, in this case, $H = 0$. According to the
identity of the first part of Lemma \ref{L:extlem}, the quasi-periodic
extension of $E$ corresponding to $\bsd_{1},\dots,\bsd_l$ is also a
direct summand of $Q^\prime$.  Quotienting by it, we may assume that
also that $D_{in} = 0$, i.e.\ that $l = 0$.  But then
$\bsd_1,\dots,\bsd_l$ are representatives of linearly independent
classes in $H_{sr}(\Phi)$, and therefore by Proposition \ref{P:minqp}
in this case, the original extension is minimal.  Consequently, since
$Q^\prime$ projects onto $E$, $Q^\prime = Q$, and $H = 0$, as desired.
\end{proof}

As a special case, we state the following corollary:

\begin{corollary}
  Let $\Delta$ be the space spanned over $L$ by the
  $\Phi$-biderivations $\bsd_1,\dots,\bsd_j$ and let $Q$ be the
  corresponding quasi-periodic extension of $E$.
\begin{enumerate}  
\item[(a)] $Q$ is a split extension of $E$ if and only if $\Delta
  \subset \Der_{in}(\Phi)$.
\item[(b)] $Q$ is a minimal extension of $E$ if and only if $0 =
  \Delta \cap \Der_{in}(\Phi)$
\end{enumerate}
\end{corollary}

\subsection{Cohomology under Isogeny of Abelian $t$-modules}

Here we extend a remark of \cite{wdb96} to arbitrary abelian
$t$-modules. If $\Theta$ is a $t$-module morphism from $E_1 =
(\Phi_1,\Ga^{d_1})$ to $E_2 = (\Phi_2,\Ga^{d_2})$, then it induces a
$\C$-linear map $\Theta^\ast : \Der(\Phi_2) \to \Der(\Phi_1)$ via
\[
(\Theta^\ast \bsd)(a) \assign \bsd(a) \Theta,\ \forall a \in \eA.
\]
One sees from its functional equation that the quasi-periodic function
associated to $\Theta^\ast \bsd$ is
\begin{equation} \label{E:isoqpf}
F_{\Theta^\ast \bsd}(\bz) = F_{\bsd}(d\Theta \bz).
\end{equation}
When we choose biderivations $\bsd_1,\dots,\bsd_j$ reducing to a basis
for $H_{sr}(\Phi_1)$ and $j = r - d + \rank N$, we can write $\Theta^\ast
\bsd = \bsd^{(\bsU)} + \sum_{i = 1}^j c_i \bsd_i.$  Then by
Proposition~\ref{P:FU} we can also write
\begin{equation} \label{E:qpisogfn}
F_{\Theta^\ast \bsd}(\bz) = \bsU\Exp_{E_1}(\bz) - d\bsU \bz+
\sum_{i=1}^j c_i F_{\bsd_i}(\bz).
\end{equation}

\begin{proposition} \label{P:isobider} If $\Theta$ is an isogeny from
  the abelian $t$-module $E_1 = (\Phi_1, \Ga^d)$ to $E_2 =
  (\Phi_2,\Ga^d)$ (necessarily also abelian), then $\Theta^\ast$
  induces isomorphisms
\begin{enumerate}
\item[(a)] $\Theta^\ast_{\DR} : H_{\DR}(\Phi_2) \to H_{\DR}(\Phi_1)$ and
\item[(b)] $\Theta^\ast_{sr} : H_{sr}(\Phi_2) \to H_{sr}(\Phi_1)$.
\end{enumerate}
\end{proposition}

\begin{proof}
  By isogeny $E_1$ and $E_2$ have the same dimension and rank;
  furthermore, $N_1 = d\Phi_1(t) - \theta I_d$ and $N_2 = d\Phi_2(t) -
  \theta I_d$ have the same rank because $d\Theta N_1 = N_2 d\Theta.$
  In light of the dimensions calculated in Proposition~\ref{P:qpdim},
  it suffices to show that $\Theta^\ast_{\DR}$ and $\Theta^\ast_{sr}$
  are injective.

Since $\Theta : E_1 \to E_2$ is an isogeny, there is an isogeny
$\Omega : E_2 \to E_1$ and a non-zero element $a \in \eA$ so that
\[
\Omega\Theta = \Phi_1(a).
\]
Clearly $(\Omega\Theta)^\ast = \Theta^\ast \Omega^\ast =
\Phi_1(a)^\ast$, and so we need only show that $\Phi_1(a)^\ast_{\DR}$
and $\Phi_1(a)^\ast_{sr}$ are injective.

Suppose that $\Phi_1(a)^\ast \bsd$ is inner (or strictly inner).  Then
\[
(\Phi_1(a)^\ast \bsd)(b) = \bsd(b)\Phi_1(a) = \bsU \Phi_1(b) -
\iota(\eb) \bsU
\]
for some $\bsU \in M(E_1)$ with $d\bsU \in N_1^{\perp}$ (or $d\bsU = 0$)
and for all $b \in \eA$.  

Let $\bsV \assign (\iota(a))^{-1}(\bsU - \bsd(a))$.  Since $\bsV
\Phi_1(a) = \bsU$, we see that
%\begin{align*}
\[
\bsd^{(\bsV)}(b)\Phi_1(a) = \bsV \Phi_1(ba) - \iota(b) \bsV \Phi_1(a) 
                         = \bsU \Phi_1(b) - \iota(b) \bsU  
                         = \bsd(b)\Phi_1(a).
\]
%\end{align*}
Since $M(E_1)$ has no $t$-torsion, we can cancel the common right
factor to see that $\bsd = \bsd^{(\bsV)}$, which is strictly
inner if and only if $\Phi_1(a)^\ast \bsd$ is so.
\end{proof}

\begin{corollary} Let $\Theta$ be an isogeny from the abelian $t$-module
  $E_1 = (\Phi_1, \Ga^d)$ to $E_2 = (\Phi_2,\Ga^d)$.  Let
  $\bsd_1,\dots,\bsd_j$ be representatives of linearly independent
  classes of $H_{sr}(\Phi_2)$.  Let $Q_2$ be the quasi-periodic
  extension of $E_2$ associated to $\bsd_1, \dots, \bsd_j$, and let
  $Q_1$ be the quasi-periodic extension of $E_1$ associated to
  $\Theta^\ast\bsd_1,\dots,\Theta^\ast\bsd_j$.  Then the matrix
\[
\Theta_{\ast} \assign
\begin{pmatrix}
\Theta & 0 \\
0 & I_j
\end{pmatrix},
\]
where $I_j$ is the $j\times j$ identity matrix, is an isogeny from
$Q_1$ to $Q_2$.
\end{corollary}

\begin{proof}
According to the preceding proposition, the $\Theta^\ast\bsd_i$ are
representatives of a basis for $\Der(\Phi_1)/\!\Der_{in}(\Phi_1)$.
Since 
\[
\begin{pmatrix}
\Theta & 0 \\
0 & I_j
\end{pmatrix}
\begin{pmatrix}
\Phi_1(t) & 0 & 0 & \dots & 0 \\
\Theta^*\bsd_1(t)& \theta & 0 & \dots & 0 \\
\Theta^*\bsd_2(t)& 0 &\theta  & \dots & 0 \\ 
\vdots & & & \ddots & \vdots \\
\Theta^*\bsd_j(t)& 0 & 0 & \dots & \theta
\end{pmatrix}
=
\begin{pmatrix}
\Phi_2(t) & 0 & 0 & \dots & 0 \\
\bsd_1(t)& \theta & 0 & \dots & 0 \\
\bsd_2(t)& 0 &\theta  & \dots & 0 \\
\vdots & & & \ddots & \vdots \\
\bsd_j(t)& 0 & 0 & \dots & \theta
\end{pmatrix}
\begin{pmatrix}
\Theta & 0 \\
0 & I_j
\end{pmatrix},
\]
we see that $\Theta_\ast$ is an isogeny from $Q_1$ to $Q_2$.
\end{proof}

\begin{corollary} \label{C:qpprod} Let $Q_i$ be strictly
  quasi-periodic extensions of the abelian $t$-modules $E_i$ of
  maximal dimension, $i = 1, \dots, m$.  Let $E$ be isogenous to
  $\prod E_i$.  Then any strictly quasi-periodic extension $Q$ of the
  $t$-module $E$ of maximal dimension is isogenous to $\prod Q_i$.
\end{corollary}

\begin{proof}
Let $Q_i$ arise from the $\Phi_i$-biderivations $\bsd_{i,j}, j = 1,
\dots, j_i$.  Then according to part (b) of Proposition
\ref{P:isobider}, the $\Theta^\ast \bsd_{i,j}$ are linearly
independent modulo $\Der_{in}(\Phi)$.

As in the proof of part (a) of Proposition \ref{P:isobider},
\[
\rank \left( d\Phi(t) - \theta I_d \right) =
\sum \rank \left( d\Phi_i(t) - \theta I_{d_i} \right).
\]
Moreover, $\rank E = \sum \rank E_i,  \dim E = \sum \dim E_i.$
Thus, according to Lemma \ref{L:sr},
\begin{align*}
\dim_L H_{sr}(\Phi) & = \rank E - \dim E + \rank (d\Phi(t) - \theta
I_d) \\
                    & = \sum_i \left( \rank E_i - \dim E_i + \rank
                      (d\Phi_i(t) - \theta I_{d_i}) \right) \\
                    & = \sum j_i.
\end{align*}
Consequently, the $\Theta^\ast \bsd_{i,j}$ span $H_{sr}(\Phi)$, and we
conclude via Proposition \ref{P:srext}.
\end{proof}

\begin{remark}  Let $\Theta: E_1 \to E_2$ be an isogeny between simple 
  $t$-modules defined over $L$.  Then $\Theta$ is defined over $L$, as $L$
  is algebraically closed.
\end{remark}

\begin{corollary} \label{C:isogspan}  If $E_1$ and $E_2$ are isogenous
  abelian $t$-modules defined over $L$, then the $L$-vector space
  spanned by the coordinates of the periods of $E_i$ is independent of
  $i$.  The same is true of the $L$-vector space spanned by the
  quasi-periods and the coordinates of the periods for maximal
  strictly quasi-periodic extensions defined over $L$.
\end{corollary}

\begin{proof}
  Let $\Theta : E_1 \to E_2$ be an isogeny and let $\Lambda_1,
  \Lambda_2$ denote the period lattices for $E_1,E_2$.  Then $d\Theta
  \Lambda_1 $ is an $\eA$-sublattice of $\Lambda_2$ of finite index.
Thus the $L$-spans of the lattices are the same.
  
  We saw in the preceding result that such an isogeny $\Theta$ induces
  an isogeny $\Theta_\ast$ of the related minimal quasi-periodic
  extensions.  The claim follows on appealing to \eqref{E:qpisogfn}.
\end{proof}

% soliton.tex

\section{$t$-Modules Arising from Solitons} \label{S:sol}

In~\cite{And92}, Anderson introduced the notion of a soliton function,
which is a higher dimensional analogue of the shtuka function for rank
$1$ Drinfeld modules.  Coleman had given such meromorphic functions
explicitly for the Fermat and Artin-Schreier curves \cite{rc88}, and
Thakur explained these examples in terms of the Gamma function.
Anderson observed the parallel between these ideas and certain topics
in differential equations and thus developed ``solitons.''

Sinha used the soliton theory in his Minnesota Ph.D.  thesis
\cite{SS95}, \cite{SS97d}, \cite{SS97a} to construct $t$-modules whose
periods have coordinates which are algebraic multiples of
$\Gamma(a/f)$ with $a$ and $f$ in $A$ both monic with $\deg (a) < \deg
(f)$.

In this section we propose to construct $t$-modules whose periods
involve arbitrary $\Gamma(a/f)$ with $\deg (a) < \deg (f)$, and to
determine their sub-$t$-module structure in terms of Thakur's bracket
relations. 

In order to render a coherent account here of our own considerations,
we have found it necessary to first recapitulate in
Sections~\ref{SS:solf}--\ref{SS:per}, those concepts, constructions,
and conclusions of Sinha to which we appeal.  We believe that in this
way, the reader may discern the flow of the arguments much better.
The expert may simply glance through this section to reassure himself
or herself of the notation we have adopted.  The newcomer may well
find it a guide for the exploration of the related aspects of Sinha's
far-ranging work.

\subsection{Soliton Functions} \label{SS:solf}

We fix some notation.  In all that follows, unless otherwise
specified, we will take fiber products and tensor products over $\Fq$:
i.e.\ ${\times} \assign {\times_{\Fq}}$ and ${\otimes} \assign
{\otimes_{\Fq}}$.

Recall that the Carlitz module $(C, \Ga)$ is the $1$-dimensional
$t$-module defined by
\[
  C(t) = \theta + \tau,
\]
with exponential function $e_C : \C \to \C$.  Fix throughout $f
\in A_{+} = \{ a \in A : \text{$a$ is monic} \}$.  Let
\[
  \zeta_f \assign e_C \left( \frac{\tpi}{f} \right),
\]
and let
\[
  B \assign B_f \assign \Fq[\theta,\zeta_f]
\]
be the integral closure of $A$ in $K \assign K_f \assign k(\zeta_f)$.
The group $(A/f)^{\times}$ is isomorphic to the Galois group
$\Gal(K/k)$ via
\begin{equation} \label{E:Gal}
\sigma_a \colon e_C \left( \frac{b\tpi}{f} \right)
  \longmapsto e_C \left( \frac{ab\tpi}{f} \right), \quad \forall b \in A.
\end{equation}
The infinite place $\infty$ of $k$ is ramified in $K$, with inertia
and decomposition group the natural subgroup $\Fq^{\times} \subset
(A/f)^{\times}$.  Thus if we take $K_+ \subset K$ to be the maximal
subfield in which $\infty$ is totally split, then
\[
  [K : K_+] = q - 1,
\]
and
\[
  [K_+ : k] = \bigl| \{ a \in A_{+} : \text{$\deg(a) < \deg(f)$,
  $(a,f) = 1$} \} \bigr|.
\]
We take $B_+$ for the ring of integers of $K_+$.  Moreover, if we let
\begin{equation} \label{E:cIp}
  \cI_{+} = \{ a \in A_{+} : \text{$\deg(a) < \deg(f)$,
  $(a,f) = 1$} \}
\end{equation}
and
\begin{equation} \label{E:cI}
  \cI = \{ a \in A : \text{$\deg(a) < \deg(f)$,
  $(a,f) = 1$} \},
\end{equation}
then $\Gal(K_f/k) = \{ \sigma_a : a \in \cI \}$ and $\Gal(K_+/k) =
\{ \sigma_a|_{K_+} : a \in \cI_{+} \}$.  We outline Sinha's
definition of an important class of $t$-modules $E_f \assign
(\Phi_f, \Ga^d)$ of dimension $d \assign \abs{\cI_{+}}$ using
Anderson's soliton function as follows.

Let $X/\Fq$ be an irreducible smooth projective curve with function
field $K$, i.e.\ we choose an isomorphism $\xi : \Fq(X) \to K$, which
amounts to providing a $K$-valued point $\xi \in X(K)$.  Let $U
\assign \Spec (B) \subset X$, and take $\ibar$ for the complement of
$U$, i.e.\ the points extending $\infty$ in the natural map $X \to
\PP^1$.  Let $V \subset X$ be the open subscheme which is the
complement of the zeros of $f$.  For $a \in (A/f)^{\times}$ we take
$[a]: X \to X$ over $\Fq$ for the automorphism induced by the morphism
$\sigma_a^{-1}: B \to B$ of $A$-algebras (going in the opposite
direction).

\begin{theorem}[Anderson \cite{And92}] \label{T:And}
Let $f \in A_+$.  There exists a unique rational function $\phi$ on
$X \times X$, regular on $V \times U$, such that for all $a \in A$
with $\deg(a) < \deg(f)$ and $(a,f) = 1$ and for all positive integers
$N$,
\[
  1 - \phi (\mathrm{Frob}^N(\xi), [a]\xi)
= \prod_{\substack{n \in A_{+} \\ \deg(n) = N - 1}}
 \left(1 + \frac{a}{fn}\right)
\]
where $\Frob : X \to X$ is the $q$-th power Frobenius morphism.
\end{theorem}

In fact Anderson and Sinha completely determine the divisors of $\phi$
and $1 - \phi$ on $X \times X$, see \cite{And92}, \cite{SS95} for
details.  These divisor identities led Anderson to term $\phi$ a
\emph{soliton function for $X$.}

We will not use the full two-variable version of soliton functions,
but rather a one-variable version as follows.  We extend scalars on
$X$ by an algebraically closed field $L/k$ (e.g., $L = \ok$ or $L =
\C$).  We then obtain the curve over $L$
\[
  \bX \assign \Spec L \times X.
\]
We note that there are copies of the rings, $A$, $B$, $K$, etc.,
contained in the function field of $\bX$: once as \emph{functions on
$\bX$} and once as \emph{scalars.}  To make the distinction between
the two interpretations, we take $\theta \in L$ for the constant
function on $\bX$ and $t \in L(\bX)$ for the function on $\bX$ such
that $t(\xi) = \theta$.  On $\bX$, the function $t - \theta$ then has
divisor
\begin{equation} \label{E:divt}
  \dv(t - \theta) = \biggl( \sum_{\ea \in \cI} [\ea] \xi \biggr)
- (q-1) I,
\end{equation}
where
\[
  I \assign \Spec L \times \ibar.
\]
When appropriate, we will use the notation $\eA = \Fq[t]$, $\eB$,
$\eK$, etc., for copies of $A$, $B$ and $K$ with underlying variable
$t \in \eA$ playing the role of $\theta \in A$.  As before, define
$\iota : t \mapsto \theta$ to be the corresponding isomorphism fixing
$\Fq$. 

Via the action of the $q$-th power Frobenius $\tau : L \to L$, it is
possible to conjugate functions, divisors, etc. on $\bX$.  For a
rational (or analytic) function $r$ on $\bX$, we let $r^{(1)}$,
$r^{(2)}$, $\dots$, denote successive conjugations.  Note that
$r^{(1)}$ is obtained by raising the coefficients of $r$ to the $q$-th
power.  Likewise, for a divisor $D$ on $\bX$, the conjugate $D^{(1)}$
is obtained by raising the coordinates of the points in the support of
$D$ (in some, and thus every, coordinate system) to the $q$-th power.

\begin{definition}
We can also pull-back functions on $X \times X$ by the natural map
\[
  \bX = \Spec L \times X \to X \times X.
\]
Given $\ef \in \eA_{+}$, the single-variable \emph{Anderson-Coleman
soliton function $g_{\ef}$,} or simply \emph{soliton function,} is
defined to be the rational function in $L(\bX)$,
\[
  g \assign g_{\ef} \assign 1 - \phi_{\ef} :
  \bX \to L.
\]
\end{definition}

Anderson's theorem now leads to the following result, which provides
the connection between solitons and Gamma values at rational
arguments:

\begin{proposition} [Anderson \cite{And92}] \label{P:gint}
For all $a \in \cI$ and $N > 0$,
\[
  g_{\ef}^{(N)}([\ea]\xi) = \prod_{\substack{n \in A_{+} \\
\deg(n) = N - 1}} \left(1 + \frac{a}{fn}\right).
\]
\end{proposition}

Note that the right-hand side in the above equation is the reciprocal
of a partial product of $\Gamma(a/f)$.

\subsection{Soliton $t$-Modules} \label{SS:solt}

\subsubsection{Divisors of Solitons}
Fix $\ef \in \eA_+$.  We define the following effective divisors
on~$\bX$:
\begin{equation} \label{E:WXiI}
\begin{aligned} 
  W \assign W_{\ef} & \assign \sum_{j=0}^{\deg(\ef) - 2}
    \sum_{\substack{\ea \in \cI_+ \\ \deg(\ea) \leq j}}
    [\ea] \circ \Frob^{\deg(\ef) - j - 2}(\xi), \\
  \Xi \assign \Xi_{\ef} & \assign \sum_{\ea \in \cI_+} [\ea]\xi.
\end{aligned}
\end{equation}
The divisor of the soliton function $g$ is then shown to be
\begin{equation} \label{E:dvg}
  \dv(g) = W^{(1)} - W + \Xi - I.
\end{equation}
Moreover, Sinha~\cite{SS97d} uses this fact to construct $t$-modules
in the following manner.

\subsubsection{Soliton $t$-motives} \label{SS:soltmot}
Recall that $U \subset X$ is the open set isomorphic to $\Spec B$.  We
take $\bU \subset \bX$ to be
\[
  \bU \assign \Spec L \times U.
\]
Let $\Omega_{\bX}$ be the sheaf of meromorphic $1$-forms on $\bX$, and
define
\[
  M_{\ef} \assign \Gamma(\bU,\Omega_{\bX}(W))
\]
to be the usual $L$-vector space of $1$-forms which are regular on
$\bU$ with at worst poles (all simple) along $W$.  As defined,
$M_{\ef}$ is a $L[t]$-module under left-multiplication, since $t$ is a
rational function on $\bX$ which is regular on $\bU$.  In analogy with
Drinfeld's shtuka approach to Drinfeld modules, $M_{\ef}$ is made into
a $\Lt$-module by defining
\begin{equation} \label{E:tau}
  \tau \cdot m \assign g_{\ef}m^{(1)},
\end{equation}
for every $m \in M_{\ef}$.  Sinha then proves the following theorem.

\begin{theorem}[{Sinha \cite{SS97d}, {\S}3, {\S}5}] \label{T:Mot}
With the given $L[t,\tau]$-module structure, $M_{\ef}$ is a
uniformizable abelian $t$-motive defined over $L$.
\end{theorem}

Under the equivalence of categories between $t$-motives and
$t$-modules, we take
\[
  E_{\ef} \assign E(M_{\ef}),
\]
to be the abelian $t$-module associated to $M_{\ef}$ (see
Section~\ref{SS:tmod}).  Sinha then computes the dimension and rank of
$E_{\ef}$ as a $t$-module to be:
\begin{equation} \label{E:d}
  d \assign \dim(E_{\ef}) \assign \rank_{\Lt} M_{\ef} =
  \abs{\cI_{+}},
\end{equation}
and
\begin{equation} \label{E:r}
  r \assign \rank(E_{\ef}) \assign \rank_{L[t]} M_{\ef} =
  \abs{\cI},
\end{equation}
where $\cI_{+}$ and $\cI$ are defined as in \eqref{E:cIp} and
\eqref{E:cI}.  We put $E_{\ef} = (\Phi_{\ef},\Ga^d)$, i.e\ we take
$\Phi \assign \Phi_{\ef}$ for the homomorphism defining the
$t$-module action on $E_{\ef}$.

\begin{remark}
The dimension of $E_{\ef}$ is equal to $[\eK_+ : \ek]$, where $\eK_+$
is an extension of $\ek$ which is totally split at $\infty$, i.e.\ the
maximal real subfield of $\eK_{\ef}$.  Furthermore, the elements of
$\eB_{\ef}$ are rational functions on $\bX$ which are regular on
$\bU$, so $M_{\ef}$ is a $\eB_{\ef}$-module.  In Section
\ref{SS:per} (see Remark~\ref{R:ECM}) we will see how the action of
$d\Phi(\eB_{\ef})$ on $\Lie(E_{\ef})$ extends the conjugate action of
$\eB_+$ as in Sections \ref{SS:HBD}--\ref{SS:CM}.  Thus $E_{\ef}$ will
be seen to be a Hilbert-Blumenthal-Drinfeld module with real
multiplications by $\eB_+$ and with complex multiplications by all of
$\eB_{\ef}$, making it a $t$-module of CM-type.
\end{remark}

\subsubsection{Exponential Functions of Soliton $t$-Modules}
We proceed as in Sinha~\cite{SS97d}, {\S}3--4.  For the rest of this
section we take $L = \C$.  We begin with a discussion on analytic
functions on $\bX$.  First, as $\C(t) = \C \cdot \ek$ is the
function field of $\PP^1/\C$, we interpret the \emph{Tate algebra},
\[
  \fA \assign \C \{t/\theta\} \assign \left\{
  \sum_{i=0}^{\infty} a_i t^i \in \power{\C}{t} \mathrel{\bigg|}
  \lim_{i\to \infty} |\theta^i a_i| = 0 \right\},
\]
as the $\C$-algebra of functions on $\PP^1$ which are analytic in the
closed disk of radius $\abs{\theta}$ centered at $0$, i.e.\ the rigid
analytic disk of radius $\abs{\theta}$.  Extend $\fA$ to
\[
  \fB \assign \eB_{\ef} \otimes_{\eA} \fA = \fA[\zeta_{\ef}].
\]
Let $\fX$ be the rigid analytic variety associated to $\bX$.  Then
$\fB$ is the $\C$-algebra of functions on $\fX$ which are analytic on
$\fU$, defined to be the inverse image under $t$ of the closed disk of
radius $\abs{\theta}$ about $0$.  In particular $\fU$ contains the
zeros of $t - \theta$,
\[
  [\cI]\xi \assign \sum_{\ea \in \cI} [\ea] \xi,
\]
as in \eqref{E:divt}.  Moreover, $\fU$ is the affinoid variety
contained in $\fX$ associated to $\fB$.  Note that, as the $\C$-valued
points of $\bX$ and $\fX$ exactly coincide, we are free to consider
points and divisors on $\bX$ as points and divisors on $\fX$.

Following Sinha, we let
\begin{equation} \label{E:cR}
  \cR \assign \frac{ {\displaystyle
\Big\{ \alpha \in \Gamma(\fU,\cO_{\fX}(-W + \Xi))
\mathrel{\Big|}
\alpha - \frac{\alpha^{(1)}}{g} \in
\Gamma(\bU,\cO_{\bX}(-W + \Xi)) \Big\} }
}
{\Gamma(\bU,\cO_{\bX}(-W))}.
\end{equation}
In light of the discussion in the previous paragraph, modulo those
rational functions on $\bX$ which vanish along $W$, elements of $\cR$
are analytic functions on $\fU$ which have zeros along $W$ and
possibly poles (all simple) along $\Xi$ and which, under the operation
\begin{equation} \label{E:alph}
  \alpha \mapsto \alpha - \frac{\alpha^{(1)}}{g},
\end{equation}
become rational functions on $\bX$ (which themselves vanish along $W$
and might have poles along $\Xi$).  The operation in \eqref{E:alph}
allows us to continue $\alpha \in \cR$ meromorphically to all of
$\bU$, with poles possibly at points of $\Xi^{(i)}$ for $i \geq 0$.

Sinha then obtains the exponential function of $E_{\ef}$ by
constructing a commutative diagram via residues
\begin{equation} \label{E:diag}
\begin{array}{c}
  \xymatrix{ & \cR \ar[dr]^{\RExp\, =\, \res_{\Xi^{(0)} +
 \Xi^{(1)} + \cdots}} \ar[dl]_{\RLie\, =\,
 \res_{\Xi}} & \\ \Lie(E_{\ef}) \ar[rr]^{\Exp} & &
 E_{\ef}(\C), }
\end{array}
\end{equation}
where $\res_{\Xi}$ represents the sum of residues over all points of
$\Xi$ and the map $\RLie : \cR \to \Lie(E_{\ef})$ is an isomorphism
of $\Fq$-vector spaces.

To define $\RLie$ and $\RExp$, we first choose a basis $\{ n_1, \dots,
n_d \}$ for $M_{\ef}$ as a $\Ct$-module.  Since $M_{\ef}$ is defined
to be $\Gamma(\bU,\Omega_{\bX}(W))$, it follows that for any $\alpha
\in \cR$, the differential $\alpha n_j$ will be regular on $\bU$
except for possible poles along $\Xi^{(i)}$, $i \geq 0$.  The map
\[
  \RLie : \cR  \to  \Lie(E_{\ef}) \simeq \C^d
\]
is defined by
\begin{equation} \label{E:RLie}
  \RLie : \alpha \mapsto
        \begin{pmatrix}
                  \res_{\Xi}(\alpha n_1) \\
                  \vdots \\
                  \res_{\Xi}(\alpha n_d)
        \end{pmatrix}.
\end{equation}
Sinha shows that \eqref{E:RLie} is an isomorphism of
$\Fq$-vector spaces via the calculation of certain $\mathrm{Ext}$
groups (see \cite{SS97d}, {\S}4).  The exponential function is then
obtained from the map
\[
 \RExp : \cR \to E_{\ef}(\C)
\]
defined by
\begin{equation} \label{E:RExp}
 \RExp : \alpha \mapsto
        \begin{pmatrix}
                  \sum_{i=0}^{\infty} \res_{\Xi^{(i)}}(\alpha n_1) \\
                  \vdots \\
                  \sum_{i=0}^{\infty} \res_{\Xi^{(i)}}(\alpha n_d)
        \end{pmatrix},
\end{equation}
so that $\Exp_{E_{\ef}}(\bz)$ is obtained by composition of
$\RExp$ with the inverse of $\RLie$.

\begin{remark_num} \label{R:exp1}
  Because we will need the technique later, we demonstrate that the
  map $\RExp$ satisfies the functional equation determined by the
  $t$-module structure on $E_{\ef}$. It can further be shown that
  map $\RExp \circ \RLie^{-1}$ is analytic and that its coordinate
  functions are normalized as in Section~\ref{SS:Exp} (cf.
  Proposition~\ref{P:RF}).  For our given $\Ct$-basis for
  $M_{\ef}$, we have a representation
\[
 t n_i = \sum_j \Phi(t)_{ij} n_j,
\]
where $\Phi_f(t) = (\Phi(t)_{ij}) \in \Mat_{d \times d}(\Ct)$ is the
multiplication-by-$t$ action on $E_{\ef}$ as in \eqref{E:Phi}.
Therefore
\[
 \res_{\Xi^{(0)}+\Xi^{(1)}+\cdots} (t\alpha n_i)
        = \sum_j \res_{\Xi^{(0)}+\Xi^{(1)}+\cdots}
         (\alpha \Phi(t)_{ij} n_j).
\]
To verify that $\RExp$ satisfies the functional equation, it is enough
to show for every $m \in M_{\ef}$ that
\[
 \res_{\Xi^{(0)}+\Xi^{(1)}+\cdots} (\alpha\tau m)
 = (\res_{\Xi^{(0)}+\Xi^{(1)}+\cdots} (\alpha m))^q.
\]
Recall from the definition of $\alpha$ in \eqref{E:cR} that $g\alpha -
\alpha^{(1)} = s \in \Gamma(\bU,\cO_{\bX}(-W^{(1)}))$; thus
$sm^{(1)}$ has no poles on $\bU$.  Therefore, from \eqref{E:tau}
\begin{align}
 \res_{\Xi^{(0)}+\Xi^{(1)}+\cdots} (\alpha \tau m)
 & = \res_{\Xi^{(0)}+\Xi^{(1)}+\cdots} (\alpha g m^{(1)})
                  \notag \\
 & = \res_{\Xi^{(0)}+\Xi^{(1)}+\cdots} (\alpha^{(1)} m^{(1)})
        + \res_{\Xi^{(0)}+\Xi^{(1)}+\cdots} (s m^{(1)})
         \label{E:rtau} \\
 & = (\res_{\Xi^{(0)}+\Xi^{(1)}+\cdots} (\alpha m))^q. \notag
\end{align}
\end{remark_num}

\begin{remark_num} \label{R:exp2} Again because we will need the technique
  below, we record here the fact that the map $\RExp : \cR \to
  E_{\ef}(\C)$ can also be defined by taking residues along $I =
  \Spec L \times \ibar$: For $\alpha \in \cR$, let $r \assign \alpha -
  \alpha^{(1)}/g \in \Gamma(\bU, \cO_{\bX}(-W + \Xi))$.
\[
 \RExp : \alpha \mapsto - \begin{pmatrix}
        \res_I \left( \left( r + \frac{r^{(1)}}{g^{(0)}} +
                  \frac{r^{(2)}}{g^{(0)}g^{(1)}} + \cdots
                  \right) n_1 \right) \\
        \vdots \\
        \res_I \left( \left( r + \frac{r^{(1)}}{g^{(0)}} +
                  \frac{r^{(2)}}{g^{(0)}g^{(1)}} + \cdots
                  \right) n_d \right)
 \end{pmatrix}.
\]
This representation relies on the following lemma (see
\cite{rc89}, {\S}2, or \cite{FvdP81}, {\S}I.3.3).

\begin{lemma} \label{L:sum0}
  Let $\beta$ be a meromorphic function on $\bU$ with discrete poles,
  and let $n$ be an algebraic differential form on $\bX$.  Let $m_j$
  be a strictly increasing sequence of integers such that $\omega
  \assign \beta n$ is regular on $\{ \gamma \in \bU : \abs{t(\gamma)}
  = \abs{\theta}^{m_j} \}$.  If for every real number $\rho > 0$,
  $\lim_{j \to \infty} \,\rho^{m_j} \sup_{\abs{t(\gamma)} =
    \abs{\theta}^{m_j}} \bigl\{\abs{\beta}\bigr\} = 0$, then
\[
\sum_{\gamma \in \bU} \res_{\gamma}(\omega) = 0.
\]
\end{lemma}

First, using the identity
\begin{equation} \label{E:itwist}
\alpha - \frac{\alpha^{(N+1)}}{g^{(0)} \cdots g^{(N)}} =
r + \frac{r^{(1)}}{g} + \dots + \frac{r^{(N)}}{g \cdots g^{(N-1)}},
\end{equation}
which is established recursively, we find that
\begin{align*}
 \sum_{i=0}^{\infty} \res_{\Xi^{(i)}} (\alpha n_j)
        & = \sum_{i=0}^{\infty} \res_{\Xi^{(i)}} \left(
        \frac{\alpha^{(N+1)}n_j}{g^{(0)}\cdots g^{(N)}}
        \right) \\
 & \quad {} - \res_I \left(\frac{r^{(N)}n_j}{g^{(0)} \cdots g^{(N-1)}}
        \right) - \cdots - \res_I(rn_j).
\end{align*}
Scondly Sinha \cite{SS97d}, Lem.~4.6.4, shows that for any real number
$\rho > 0$, we have
\[
\lim_{N \to \infty} \rho^{q^N} \!\sup_{\abs{t(\gamma)} =
\abs{\theta}^{q^N - q^{N-2}}}
\bigl\{\abs{\alpha^{(N)}n_j/(g^{(0)} \cdots g^{(N-1)})}\bigr\} = 0.
\]
Combining these two facts with Lemma~\ref{L:sum0}, Sinha determines
for $N$ sufficiently large that
\begin{equation} \label{E:est}
 \sum_{\gamma \in \bU} \res_{\gamma} \left( \frac{\alpha^{(N +
                  1)}n_j}{g^{(0)} \cdots g^{(N)}} \right)
        = \sum_{i=0}^{\infty} \res_{\Xi^{(i)}} \left( 
                  \frac{\alpha^{(N+1)}n_j}{g^{(0)}\cdots
                  g^{(N)}} \right) = 0,
\end{equation}
since on $\bU$ the poles of $\alpha^{(N)}n_j/(g^{(0)} \cdots
g^{(N-1)})$ lie along $\Xi^{(0)} + \Xi^{(1)} + \cdots$.
\end{remark_num}

\subsection{Periods of Soliton $t$-Modules} \label{SS:per}

In the next important step, Sinha determines the kernel of the map
$\RExp : \cR \to E_{\ef}(\C)$ from \eqref{E:diag}, and \emph{a
fortiori} obtains the period lattice of $E_{\ef}$.

Consider the infinite product
\[
 c_g = \frac{1}{g^{(0)} g^{(1)} \cdots},
\]
which converges in the space of rigid analytic maps $\Gamma(\fU,
\cO_{\fX}(-W+\Xi))$.  Clearly,
\[
        c_g - \frac{c_g^{(1)}}{g} = 0,
\]
so $c_g \in \cR$.  From Remark~\ref{R:exp2} it then follows that $c_g$
is in the kernel of the exponential map $\RExp$, and in fact Sinha
shows that
\[
\Ker(\RExp) = \eB_{\ef} \cdot c_g.
\]
The period lattice $\Lambda_{\ef}$ of $E_{\ef}$ is then the
image of $\Ker(\RExp)$ under $\RLie : \cR \to \Lie(E_{\ef})$.
Thus $\Lambda_{\ef}$ is an $\eA$-module of rank $[\eK_{\ef}:\ek]$,
since it is a free $\eB_{\ef}$-module of rank~$1$.

Following Sinha, we now choose a $\Ct$-basis for $M_{\ef}$ which
is amenable to computing the period lattice explicitly.  Define
\[
\varepsilon(\ea) = \begin{cases}
1, & \text{if $\deg(\ea) = \deg(\ef) - 1$,} \\
0, & \text{if $\deg(\ea) < \deg(\ef) - 1$.}
\end{cases}
\]
Note for $\ea \in \cI$ that $\varepsilon(\ea)$ is the multiplicity of
$[\ea]\xi$ in $\Xi - W$.  Sinha then chooses a $\Ct$-basis $\{ n_{\ea}
\}_{\ea \in \cI_+}$ for $M_{\ef}$ so that $n_{\ea}$ is defined
over $\ok$;
\[
\begin{array}{c}
        n_{\ea} \in \Gamma(\bU, W - \Xi + [\ea]\xi); \\
        \res_{[\ea]\xi} ((n_{\ea})/(t - \theta)^{\varepsilon(\ea)}) =
        1;
\end{array}
\]
and
\[
        \text{$n_{\ea}$ has} \left\{ \begin{array}{ll}
                         \left. \begin{array}{p{2.5truein}}
                         \text{a pole of order $1$ at $[\ea]\xi$,} \\
                         \text{no pole at $[\eb]\xi$ if $\eb \neq \ea$,
                         $\varepsilon(\eb) = 0$,} \\ 
                         \text{a zero at $[\eb]\xi$ if
                                $\varepsilon(\eb) = 1$,}
                         \end{array} \right\}
                         & \text{if $\varepsilon(\ea) = 0$;} \\ [20pt]
                         \left. \begin{array}{p{2.5truein}}
                         \text{no pole or zero at $[\ea]\xi$,} \\
                         \text{no pole at $[\eb]\xi$ if
                                $\varepsilon(\eb) = 0$,} \\
                         \text{a zero at $[\eb]\xi$ if $\eb \neq \ea$,
                         $\varepsilon(\eb) = 1$,}
                         \end{array} \right\}
                         & \text{if $\varepsilon(\ea) = 1$.}
        \end{array}
        \right.
\]
Thus determining $\Lambda_{\ef}$ amounts to computing
\[
\pi_{\ea} \assign \res_{\Xi}(c_g n_{\ea}) = \res_{[\ea]\xi}
(c_g n_{\ea}),
\]
since the image of arbitrary $\eb \cdot c_g$ in $\Lambda_{\ef}
\subset \Lie(E_{\ef})$ is then
\begin{equation} \label{E:pers}
\RExp(\eb \cdot c_g)  = \begin{pmatrix}
  \vdots \\
  \sigma_{\ea}(\iota(\eb)) \pi_{\ea} \\
  \vdots
\end{pmatrix}_{\ea \in \cI_+}.
\end{equation}
The function $\eb$ is evaluated at the point $[\ea]\xi$, for which we
see $\eb([\ea]\xi) = \sigma_a(\iota(\eb))$.  Note that this shows that
the action of $d\Phi(\eB_+)$ on $\Lie(E_{\ef})$ is the conjugate
action $\bss(\eB_+)$ of $\eK_+/\ek$, where $\bss = {\oplus}
\sigma_{\ea}|_{\eK_+}$ for $\ea \in \cI_+$.

We now compute $\res_{[\ea]\xi} (c_g n_{\ea})$.  From \eqref{E:WXiI}
and \eqref{E:dvg} it follows that the function $c_g = (g^{(0)}g^{(1)}
\cdots )^{-1}$ has poles among the points $[\ea]\xi$ of $\Xi$ exactly
when $\varepsilon(\ea) = 1$.  Thus if $\varepsilon(\ea) = 0$, then by
Proposition~\ref{P:gint},
\begin{align}
\pi_{\ea} = \res_{[\ea]\xi} (c_g n_{\ea})  &= c_g([\ea]\xi)
  \notag \\
  &= g([\ea]\xi)^{-1} \prod_{n \in A_{+}} \left( 1 + \frac{a}{fn}
  \right)^{-1} \label{E:pia} \\
  &= g([\ea]\xi)^{-1} \frac{a}{f} \Gmm{\frac{a}{f}}. \notag
\end{align}
Likewise, if $\varepsilon(\ea) = 1$, then by a similar calculation
\[
\pi_{\ea} = \res_{[\ea]\xi} (c_g n_{\ea}) = \res_{[\ea]\xi}
(n_{\ea}/g) \frac{a}{f} \Gmm{\frac{a}{f}}.
\]
We see right away for all $\ea \in \cI_+$ that
\begin{equation} \label{E:piasim}
\pi_{\ea} \thicksim \Gmm{\frac{a}{f}},
\end{equation}
because $n_{\ea}$, $g$, and $\xi$ are defined over $\ok$.  Sinha
further computes these algebraic factors and obtains the following
explicit result after some careful analysis.

\begin{theorem}[Sinha {\cite{SS97d}, {\S}5.3.9}] \label{T:Spia}
  Let $a$, $f \in A$ be monic with $\deg(a) < \deg(f)$ and $(a,f) =
  1$.  Then
\[
\pi_{\ea} = \begin{cases}
   {\displaystyle \Gmm{\frac{a}{f}},}
   & \text{if $\deg(a) = \deg(f) - 1$,} \\ [10pt]
   {\displaystyle \frac{a}{f} \Gmm{\frac{a}{f}},}
   & \text{if $\deg(a) < \deg(f) - 1$.}
\end{cases}
\]
\end{theorem}

\begin{remark_num} \label{R:ECM}
{From} the preceding choice of coordinates on $\Lie(E_{\ef})$, in
particular \eqref{E:pers}, we see that $E_{\ef}$ is a H-B-D module
with real multiplications by $\eB_+$.  Furthermore, if we let
\[
\cS_{\ef} = \{ \sigma_{\ea} \in \Gal(\eK_{\ef}/\ek) : \ea \in \cI_+ \},
\]
then $\cS_{\ef}$ is an extension of $\Gal(\eK_+/\ek)$ to
$\Gal(\eK_{\ef}/\ek)$.  From \eqref{E:pers} we see that $E_{\ef}$ is a
$t$-module of CM-type $(\eK_{\ef}, \cS_{\ef})$ with complex
multiplications by $\eB_{\ef}$ and conjugate action $\bss_{\ef}
\assign \bss_{\cS_{\ef}}$.
\end{remark_num}

\subsection{Quasi-Periodic Soliton $t$-Modules and Quasi-Periods}

We now apply the ideas underlying Sinha's construction to the
quasi-periodic extensions of soliton $t$-modules defined
Section~\ref{S:qp}.  Since $d\Phi_{\ef}(t)$ for the soliton
$t$-module $E_{\ef}$ has no nilpotent part, i.e.\ $(t -
\theta)M_{\ef} \subset \tau M_{\ef}$, we see from
Proposition~\ref{P:qpdim} that
\[
        \dim_{\C} \Der_0(\Phi_{\ef}) = d \quad \text{and} \quad
        \dim_{\C} H_{sr}(\Phi_{\ef}) = r - d,
\]
where $d$ and $r$ are the dimension and rank of $E_{\ef}$ defined in
\eqref{E:d} and \eqref{E:r}.  Moreover, we know that
\begin{equation} \label{E:Dsr}
        H_{sr} (\Phi_{\ef}) \simeq \frac{\tau
        M_{\ef}}{(t-\theta)M_{\ef}},
\end{equation}
and our immediate task is to find a convenient basis for this space
(defined over $\ok$) and compute the associated quasi-periodic
functions and quasi-periods.

\subsubsection{Quasi-periodic Functions}

As in Section~\ref{SS:bider}, given an element of the $t$-motive $m
\in \tau M_{\ef}$, we have an associated biderivation $\bsd
\assign \bsd_m : \eA \to (\Ct \tau)^d$ defined by
\[
        \bsd(t) = (\delta_1(t), \dots, \delta_d(t)),
\]
where $m = \sum \delta_j(t) n_j$ for a fixed $\Ct$-basis $\{n_1, \dots,
n_d \}$ of $M_{\ef}$.

We now proceed to express the associated quasi-period function
$\Fd : \Lie(E_{\ef}) \to \C$ in terms of residues, as in \eqref{E:diag}.
Consider the map $\RF : \cR \to \C$ defined by
\[
        \RF : \alpha \mapsto \sum_{i=1}^{\infty} \res_{\Xi^{(i)}}
        \left( \frac{\alpha m}{t - \theta} \right).
\]
If $s \in \Gamma(\bU, \cO_{\bX}(-W))$, then $sm/(t-\theta) \in
\Gamma(\bU,\Omega_{\bX}(-[\cI]\xi - \Xi))$. Thus
$\RF(s) = 0$, and so $\RF$ is well-defined.  Moreover, $\RF$
essentially gives the quasi-periodic function for $\bsd$ according to
the following proposition.

\begin{proposition} \label{P:RF}
For $\bsd = \bsd_m$, the corresponding quasi-periodic function is
given by
\[
F_{\bsd} \assign \RF \circ \RLie^{-1}.
\]
Thus, by definition, the following diagram commutes:
\[
        \xymatrix{ & \cR \ar[dr]^{\RF} \ar[dl]_{\RLie\, =\,
                        \res_{\Xi}} & \\ \Lie(E_{\ef}) \ar[rr]^{\Fd} &
                        & \C. }
\]
\end{proposition}

\begin{proof}
  A. We first demonstrate that $\RF$ satisfies the functional equation
  of $\Fd$, proceeding as in Remark~\ref{R:exp1}.  Using $m = \sum
  \delta_j(t) n_j$, we see that
\begin{align*}
  \sum_{i=1}^{\infty} \res_{\Xi^{(i)}} \left( \frac{t\alpha m}{t -
      \theta} \right) & = \sum_{i=1}^{\infty} \res_{\Xi^{(i)}} \left(
    \frac{(t - \theta + \theta)\alpha m}{t - \theta} \right)
  \\
  & = \sum_{i=1}^{\infty} \res_{\Xi^{(i)}} (\alpha m) +
  \sum_{i=1}^{\infty} \res_{\Xi^{(i)}}
  \left( \frac{\theta \alpha m}{t - \theta} \right). \\
  \intertext{Note by \eqref{E:dvg} and \eqref{E:tau}, that $m \in \tau
    M_{\ef} = \tau\Gamma(\bU,\Omega_{\bX}(W)) = \Gamma(\bU,
    \Omega_{\bX}(W - \Xi))$, and so $\res_{\Xi^{(0)}} (\alpha m) = 0$.
    Therefore, the above calculation continues}
& = \sum_{i=0}^{\infty} \res_{\Xi^{(i)}}
  \left( \alpha \sum_j \delta_j(t) n_j \right) + \sum_{i=1}^{\infty}
  \res_{\Xi^{(i)}}
  \left( \frac{\theta \alpha m}{t - \theta} \right) \\
  & = \bsd(t) (\RExp(\alpha)) + \theta \RF(\alpha),
\end{align*}
where the last equality follows exactly as in \eqref{E:rtau}, using
the definitions of $\RExp$ and $\RF$.

B. Next we show that the composite $\RF \circ \RLie^{-1}:
\Lie(E_{\ef}) \to \C^d$ is an entire function.  For that we choose a
system $z_1, \dots, z_d$ of coordinates for $\Lie(E_{\ef})$.  Given
$\alpha \in \cR$, set $z_j \assign \res_{\Xi^{(0)}}(\alpha n_j)$, $j =
1, \dots, d$.  Thus from the definitions of the various twists
involved, we see that
\begin{equation} \label{E:zjqi}
        z_j^{q^i} = \res_{\Xi^{(i)}} (\alpha^{(i)} n_j^{(i)}).
\end{equation}
Our goal is to express $\res_{\Xi^{(i)}} (\alpha m/(t- \theta))$ as a
$\C$-linear combination of $z_j^{q^i}$ with coefficients which
decrease sufficiently rapidly with $i$.

We make some preliminary observations.  Let $r \assign \alpha -
\frac{\alpha^{(1)}}{g} \in \Gamma(\bU, \cO_{\bX}(-W + \Xi))$.  
Recall \eqref{E:itwist} and use the fact 
\[
\frac{r^{(j)}}{g^{(0)} \cdots g^{(j-1)}} \in
\Gamma(\bU,\cO_{\bX}(-W + {\textstyle \sum_{l=0}^j \Xi^{(l)}})),
\]
which is verified recursively, to see that
\begin{equation} \label{E:aai}
        \res_{\Xi^{(i)}} \left( \frac{\alpha m}{t - \theta} \right)
                = \res_{\Xi^{(i)}} \left( \frac{\alpha^{(i)}}{g^{(0)}
                \cdots g^{(i-1)}}\, \frac{m}{t - \theta} \right).
\end{equation}
Next, we write $m$ as
\begin{equation} \label{E:m}
        m = \sum_j \delta_j(t) n_j = \sum_{j=1}^d \sum_{s=1}^u
                \gamma_{sj} \tau^s n_j, \quad \gamma_{sj} \in \C.
\end{equation}
Combining \eqref{E:aai} and \eqref{E:m} we then have
\begin{equation} \label{E:start}
\res_{\Xi^{(i)}} \left( \frac{\alpha m}{t - \theta} \right)
  = \sum_{k=1}^d \sum_{s=1}^u \res_{\Xi^{(i)}} \left(
  \frac{\alpha^{(i)}}{g^{(0)} \cdots g^{(i-1)}} \,
  \frac{\gamma_{sk}}{t- \theta}\, \tau^s n_k \right).
\end{equation}

Finally, since for each $k$, $1 \leq k \leq d$, we have $(t - \theta)
n_k \in \tau M_{\ef}$, if we set $\bn = (n_1,\dots,n_d)^{tr}$, we can
write
\[
(t - \theta) \bn = B_1\tau \bn + \dots + B_v\tau^v \bn,
\]
where each $B_k \in \Mat_{d \times d}(\C)$.  Then we divide by $t -
\theta$ and substitute this expression in for $\bn$ in the term $B_1
\tau \bn$ to obtain
\[
\bn = B_2^{[1]} \tau^2 \bn + \dots + B_{v+1}^{[1]} \tau^{v+1} \bn,
\]
where each $B_i^{[1]}\tau^i$ has entries which are the sum of terms of
the form $b_i/(t-\theta)\tau^i$ or $b_1(t-\theta)^{-1} \tau \cdot
((t-\theta)^{-1}b_{i-1}) \tau^{i-1}$, where the $b_i$ are coefficients
of $B_i$.  Continuing in this way, we can progressively eliminate the
first $j-1$ degree terms to find an expression
\begin{equation} \label{E:bn}
\bn = B_{j+1}^{[j]} \tau^{j+1} \bn + \dots + B_{v+j}^{[j]} \tau^{v+j}
  \bn,
\end{equation}
where each $B_w^{[j]}\tau^w$ has coefficients which are sums of terms
of the form
\begin{equation} \label{E:coeff}
\left(\frac{b_{e_1}}{t-\theta}\right)\tau^{e_1}
\left(\frac{b_{e_2}}{t-\theta}\right)\tau^{e_2}
\left(\frac{b_{e_3}}{t-\theta}\right)\tau^{e_3} \dots
\left(\frac{b_{e_l}}{t-\theta}\right)\tau^{e_l}
\end{equation}
for which $e_1 + e_2+ \dots + e_{l-1} + e_l = w$, for $1 \le e_i \le
v$, each $b_e$ is an entry of the constant matrix $B_e$, and $j+1 \le
w \le j+v$.

Now we substitute the expressions from \eqref{E:bn} with $j = i - s - 1$
into the terms of \eqref{E:start} involving $\tau^s n_k$.  According
to \eqref{E:tau}, when we multiply $n_k$ by the preceding expression
\eqref{E:coeff} and apply $\tau^s$, we obtain

\begin{multline} \label{E:RFterms}
  \left(\frac{b_{e_1}}{t-\theta}\right)^{(s)}
  \left(\frac{b_{e_2}}{t-\theta}\right)^{(s+e_1)}
  \left(\frac{b_{e_3}}{t-\theta}\right)^{(s+ e_1 + e_2)} \cdots
  \left(\frac{b_{e_l}}{t-\theta}\right)^{(s+ e_1 + \dots +e_{l-1})}
\times \\
g^{(0)} \cdots g^{(w+s-1)}n^{(w+s)}_k.
\end{multline}

When $w > i-s$, the above differential multiplied by 
\[
\frac{\alpha^{(i)}}{g^{(0)}\dots g^{(i-1)}} \frac{1}{t - \theta}
\]
is regular along $\Xi^{(i)}$, since then the uncancelled $g^{(i)}$
occurs in the numerator.  Therefore the contribution towards the
residue at $\Xi^{(i)}$ of such terms is nil, and we can concentrate on
the terms with $w = i - s$.  We obtain that
\[
        \res_{\Xi^{(i)}} \left(\frac{\alpha m}{t - \theta} \right)
                = \sum_{j=1}^d \res_{\Xi^{(i)}} ( G_{j,i} \alpha^{(i)}
                n_j^{(i)} ),
\]
where $G_{j,i}$ is a rational function in $t-\theta$ and certain
twists $t-\theta^{(e)}$ on $\bX$ without poles along $\Xi^{(i)}$.
Because $G_{j,i}$ is constant on $\Xi^{(i)}$, we conclude from
\eqref{E:zjqi} that
\begin{equation}
\res_{\Xi^{(i)}} \left(\frac{\alpha m}{t - \theta} \right) =
        \sum_{j=1}^d G_{j,i}(\xi^{(i)}) z_j^{q^i}.
\end{equation}
It follows that
\begin{equation} \label{E:cpst}
        \RF \circ \RLie^{-1} (z_1, \dots, z_d) =
                \sum_{i=1}^{\infty} \sum_{j=1}^d
                G_{j,i}(\xi^{(i)}) z_j^{q^i},
\end{equation}
and so, using the fact that we have already shown that this function
satisfies the proper functional equation, by Proposition \ref{P:qpfn}
it will be shown equal to $\Fd$ once we verify that the right-hand
side is entire.

To this end we need to estimate $\abs{G_{j,i}(\xi^{(i)})}$, which we
do by estimating the terms appearing in \eqref{E:RFterms}.  Since we
have the bound
\[
\frac{\left| b_{e_k}^{(s+ e_1 + \dots +e_{k-1})} \right| }
{ | \theta^{q^i}-\theta^{(s+e_1+ \dots + e_{k-1})} | } \le
\frac{B^{q^{s+\dots+e_{k-1}}}}{|\theta|^{q^i}} 
\]
on the factors occurring in \eqref{E:RFterms}, where $B$ is an upper
bound for the absolute values of the entries of the matrices $B_j$, we
see that
\[
\left| G_{j,i}(\xi^{(i)})\right| \le 
  \frac{C}{\abs{\theta}^{q^i}} \, \left| \frac{B^{q + q^2 + \dots +
        q^{i-1}}}{(\theta^{q^i})^\frac{i-s}{v}} \right| \le
        \frac{C}{\abs{\theta}^{q^i}} \,
        \frac{B^{q^i}}{\abs{\theta}^{\frac{i-s}{v}q^{i}}} \le
\frac{C_0^{q^i}}{\abs{\theta}^{\frac{i}{v}q^{i}}},
\]
where $C$ is an upper bound on the coefficients $\gamma_{sj}$ in
\eqref{E:m}.  Here we have used the fact that $s + e_1 + \dots + e_l
= i$ in \eqref{E:coeff}, whereas $1 \le e_j \le v$, $ 1 \le j \le l$.
{From} this estimate it is clear that \eqref{E:cpst} is entire, since 
$\abs{\theta} > 1$.  One can also show that the power series in
\eqref{E:cpst} is entire by appealing to Proposition 2.1.4 of
\cite{And86}. 
\end{proof}

\subsubsection{Quasi-Periods}

As defined in Section \ref{SS:qpf}, the quasi-periods of $E_{\ef}$
coming from $\bsd$ are the values $\Fd(\Lambda_{\ef}) \subset \C,$
where $\Lambda_{\ef}$ is the period lattice of $E_{\ef}$.  As in
Proposition \ref{P:qpext}, these quasi-periods are then the
coordinates of periods of the quasi-periodic extension of $E_{\ef}$
corresponding to $\bsd$.

The quasi-periods of any inner biderivation $\bsd$ are linear
combinations of the coordinates of the periods of $E_{\ef}$ with
coefficients from any field of definition of $E_{\ef}$ and $\bsd$.
Thus for questions of linear independence over $\bar{k}$ involving
both periods and quasi-periods, we may as well assume that our
biderivations are strictly reduced and defined over $\bar{k}$.  Our
first task is to find a convenient $\C$-basis for $H_{sr}(\Phi_{\ef})$
or, equivalently by \eqref{E:Dsr}, a basis for $\tau M_{\ef}/
(t-\theta)M_{\ef}$ defined over $\bar{k}$.

\begin{lemma} \label{L:na}
There exist representatives $n_a, a \in \cI \smallsetminus \cI_+$
for a $\C$-basis of
$\tau M_{\ef}/(t-\theta)M_{\ef}$ for which
\begin{enumerate}
\item[(a)] $n_{\ea}$ is defined over $\ok$;
\item[(b)] $n_{\ea} \in \Gamma(\bU,\Omega_{\bX}(W - 2[\cI]\xi + 2[\ea]\xi))$;
\item[(c)] $\res_{[\ea]\xi} (n_{\ea}/(t-\theta)) = 1$;
\item[(d)] $\res_{[\ea]\xi} (n_{\eb}/(t-\theta)) = 0$, if $\ea \ne \eb$.
\end{enumerate}
\end{lemma}

\begin{proof}
  The $\C$-linear maps $\left\{ \rho_{\ea}: m \mapsto \res_{[\ea]\xi}
    \left( \frac{m}{t - \theta} \right) : \ea \in \cI \smallsetminus
    \cI_+ \right\}$ on
\[
\tau M_{\ef} = \tau \Gamma(\bU,\Omega_{\bX}(W)) = \Gamma(\bU,
  \Omega_{\bX}(W - \Xi)),
\]
are linearly independent over $\C$.  Indeed by Riemann-Roch, for $j$
sufficiently large $\Gamma(\bX, \cO_{\bX}(jI - 2[\cI]\xi + [a]\xi))
\subsetneq \Gamma(\bX, \cO_{\bX}(jI - 2[\cI]\xi + 2[a]\xi))$ and
$\Gamma(\bX, \Omega_{\bX}(jI + W - \Xi - [a]\xi)) \subsetneq
\Gamma(\bX, \Omega_{\bX}(jI + W - \Xi))$.  Further the $\rho_a$ are
trivial on $(t - \theta)M_{\ef}$.  Thus we can choose $\{n_{\ea} \}$
as representatives of a basis modulo $(t-\theta)M_{\ef}$ dual to the
maps $\rho_{\ea}$.  That the
$n_{\ea}$ are defined over $\ok$ follows from the fact that the
$t$-motive $M_{\ef}$ is obtained by extending scalars on a $t$-motive
which is initially defined over $\ok$ (as in Theorem~\ref{T:Mot}).
\end{proof}

We now fix a basis $\{ n_{\ea} : \ea \in \cI \smallsetminus \cI_+\}$
for $\tau M_{\ef}/(t-\theta)M_{\ef}$ as in the above lemma, thus
obtaining a basis $\{ \bsd_{\ea} \}$ for $H_{sr}(\Phi_{\ef})$.
For $b \in \eB_f$, let
\[
  \bsl_b \assign \begin{pmatrix} \vdots \\
    \sigma_a(\iota(b))\pi_{\ea} \\
    \vdots \end{pmatrix}_{\ea \in \cI_+}
\]
be a period in $\Lie(E_{\ef})$ as in \eqref{E:pers}.  We define for
each $\ea \in \cI \smallsetminus \cI_+$,
\[
\eta_{a,b} \assign \eta_{a}(\bsl_b) \assign F_{\bsd_a}(\bsl_b)
\]
to be the corresponding quasi-periods.

\begin{theorem} \label{T:qp}
  Fix a basis for $\tau M_f/(t-\theta)M_f$ defined over $\ok$ as in
  Lemma~\ref{L:na}.  For each $a \in \cI \smallsetminus \cI_+$ and $b
  \in \eB_f$,
\[
\eta_{a,b} = F_{\bsd_a}(\bsl_b) = \sigma_a(b) c_{\ea}
\Gmm{\frac{a}{f}}, \quad c_a \in \ok^{\times}.
\]
\end{theorem}

\begin{remark}
  As pointed out by the referee, the explicit formulas for solitons in
  Anderson \cite{And92}, Sinha \cite{SS97d}, and Thakur \cite{Th92}
  provide a basis for determining the values of $c_a$ precisely, in
  the spirit of Theorem \ref{T:Spia}.  See Section \ref{SS:simp} for a
  special case.
\end{remark}

\begin{remark}
For an arbitrary $\Phi_f$-biderivation defined over $\ok$, the
corresponding quasi-periods are $\ok$-linear combinations of Gamma
values in
\[
  \Gamma_f \assign \{ \Gmm{a/f} : \deg(a) < \deg(f), (a,f) =1 \}.
\]
Indeed, the strictly quasi-periodic extensions $Q_f$ of $E_f$ with $j
= d - r(E_f)$ as described in Proposition \ref{P:minqp} provide
examples.  Furthermore, Theorem~\ref{T:qp} shows that by appropriately
choosing a basis for $H_{sr}(\Phi_f)$, we guarantee that the periods
of this quasi-periodic extension have coordinates which are simply
non-zero algebraic multiples of all the values in $\Gamma_f$.
\end{remark}

\begin{remark}
  It is possible to define more general soliton functions on $\bX$.
  Namely, given $\es \in \eA$ with $\deg(\es) < \deg(\ef)$, Sinha
  defines a function $\phi_{\es}$ on $\bX$ which is a certain
  pull-back of $\phi$ (see \cite{SS97d}, {\S}2.2.8).  We can proceed
  with our various concerns using this function $\phi_{\es}$ instead
  of $\phi$ and construct a corresponding $t$-module $E_{\ef,\es}$.
  Nevertheless, this $t$-module can be shown to be isogenous to $E_g$
  for some monic $g$ dividing $f$, so for considerations of
  transcendence and linear independence over $\ok$ for Gamma values,
  we gain nothing new.
\end{remark}

\begin{proof}[Proof of Theorem~\ref{T:qp}]
  Fix $\ea \in \cI \smallsetminus \cI_+$ and $b \in \eB_f$.  In the
  following, for precision we distinguish between the function $b \in
  \eB_f$ on $\bX$ and the scalar $\iota(b) \in B_f$.  As $\bsl_b =
  \RLie(b\cdot c_g)$ from \eqref{E:diag}, it follows from
  Proposition~\ref{P:RF} that
\begin{equation} \label{E:eta}
\eta_{a,b} = \cR F_{\bsd_a}(b\cdot c_g)
= \sum_{i=1}^{\infty} \res_{\Xi^{(i)}} \left(
  \frac{b c_g n_{\ea}}{t - \theta} \right).
\end{equation}
{From} Lemma~\ref{L:na} it follows that $b n_{\ea}/(t-\theta) \in
\Gamma(\bU,\Omega_{\bX}(W - [\cI]\xi + 2[\ea]\xi))$.  Thus as $c_g \in
\Gamma(\fU,\cO_{\fX}(-W+\Xi))$, the remark immediately following
equation \eqref{E:alph} shows that the poles of $b c_g
n_{\ea}/(t-\theta)$ contained in $\bU$ lie along the support of the
divisor
\[
[\ea]\xi + \Xi^{(0)} + \Xi^{(1)} + \cdots.
\]
Moreover, we establish
\begin{equation} \label{E:res}
\res_{[\ea]\xi} \left( \frac{b c_g n_{\ea}}{t - \theta} \right) 
+ \res_{\Xi^{(0)} + \Xi^{(1)} + \cdots} \left(
\frac{b c_g n_{\ea}}{t - \theta} \right) = 0.
\end{equation}
This equality follows from Lemma~\ref{L:sum0}.  The required estimates
are obtained exactly as in Sinha~\cite{SS97d}, Lem.~4.6.4, (taking
$\alpha = c_g$), and the sum on the left of \eqref{E:res} is taken
over all of the poles of $b c_g n_{\ea}/(t-\theta)$ contained in
$\bU$.  Combining \eqref{E:eta} and \eqref{E:res} we find that
\begin{equation} \label{E:efin}
\eta_{a,b} = \cR F_{\bsd_a}(b\cdot c_g) = -\res_{[\ea]\xi +
\Xi^{(0)}} \left( \frac{b c_g n_{\ea}}{t - \theta} \right).
\end{equation}
As $b n_{\ea}/(t-\theta) \in \Gamma(\bU, \Omega_{\bX}(W - [\cI]\xi +
2[\ea]\xi))$, it follows that $b c_g n_a/(t-\theta)$ is regular at
$[a']\xi$ for $a' \in \cI$, $a' \neq a$.  Lemma~\ref{L:na}bc combined
with the calculation of \eqref{E:pia} shows that
\[
\res_{[\ea]\xi} \left( \frac{b c_g n_{\ea}}{t - \theta} \right) =
\sigma_a(\iota(b)) g([\ea]\xi)^{-1} \frac{a}{f} \Gmm{\frac{a}{f}}.
\]
Because $g$ and $[\ea]\xi$ are defined over $\ok$, the constant
$c_{\ea} \assign - g([\ea]\xi)^{-1} \frac{a}{f} \in \ok$ and $c_{\ea}
\neq 0$.
\end{proof}

\subsection{Sub-$t$-modules and Connections with Bracket Relations}
\label{SS:solbrk}

In this section we investigate the correspondence between the bracket
relations on special values of the Gamma function of
Section~\ref{SS:brk} and the presence of sub-$t$-modules in soliton
$t$-modules $E_{\ef}$.  Since the soliton $t$-modules are of CM-type,
the results of Section~\ref{SS:CMsub} apply.

Let $f$ be monic and let $a \in \cI$.  Recalling in \eqref{E:Gal} that
we identify $\Gal(\eK_f/\ek)$ with $\cI$, we define the following
subset of $\Gal(\eK_f/\ek)$:
\begin{equation} \label{E:Fa}
  \eF(\ea) \assign \{ \sigma_{\es} \in \Gal(\eK_{\ef}/\eK) :
  \Gamma(s/f) \approx \Gamma(a/f) \}.
\end{equation}
Note first that $\eF(1)$ is in fact a subgroup of
$\Gal(\eK_{\ef}/\ek)$: by the bracket relations, with $m_1 = 1$, $m_s
= -1$ and all other entries of $\bm$ equal to $0$, we see that
$\sigma_{\es} \in \eF(1)$ if and only if for all representatives $\eu$
of elements of $(\eA/\ef)^{\times}$ we have
\[
\text{$\eu\es \bmod \ef$ is monic $\Longleftrightarrow$ $\eu \bmod \ef$
is monic,}
\]
where $a \bmod f$ denotes the remainder of $a$ after division by $f$.
This condition is certainly closed under multiplication.  Similarly,
for each $a \in \cI$, we find that $\eF(\ea)$ is a coset of $\eF(1)$,
i.e.\ $\eF(\ea) = \eF(1)\sigma_{\ea}$.

The set $\cS_{\ef} \subset \Gal(\eK_{\ef}/\ek)$ from
Remark~\ref{R:ECM} is defined by
\begin{equation} \label{E:Sf}
\cS_{\ef} = \{ \sigma_{\ea} \in \Gal(\eK_{\ef}/\ek) :
\text{$\ea \bmod \ef$ is monic} \},
\end{equation}
and so we find that $\cS_{\ef}$ is the union of cosets of $\eF(1)$.
Let $\eL_{\ef} \subset \eK_{\ef}$ be the fixed field of $\eF(1)$.
Theorem~\ref{T:Esub} then shows that there is a sub-$t$-module
$H_{\ef}$ of CM-type $(\eL_{\ef},\cS_{\ef}|_{\eL_{\ef}})$ such that
$E_{\ef}$ is isogenous to $H_{\ef}^m$, where $m =
[\eK_{\ef}:\eL_{\ef}]$.

\begin{proposition} \label{P:Gmult}
The relation on $\cI$ induced by $\Gamma(a/f) \approx \Gamma(b/f)$
gives a decomposition of $\cI$ into disjoint subsets of cardinality
$m$.
\end{proposition}

\begin{proof}
$\abs{\eF(a)} = \abs{\eF(1)} = [\eK_f:\eL_f] = m$.
\end{proof}

\begin{lemma}
$H_{\ef}$ is a simple $t$-module.
\end{lemma}

\begin{proof}
Since $H_{\ef}$ is itself a $t$-module of CM-type, it is isogenous to
a power of a simple $t$-module of CM-type.  Thus there is a smallest
field $\eL \subset \eL_{\ef} \subset \eK_{\ef}$ which satisfies the
criteria of Lemma~\ref{L:crit}.  Let $\eF$ be the subgroup of
$\Gal(\eK_{\ef}/\ek)$ corresponding to $\Gal(\eK_{\ef}/\eL)$.
Certainly $\eF(1) \subset \eF$.  However, if $\sigma_{\es} \in \eF$,
then since $\cS_{\ef}$ is the union of cosets of $\eF$, it must be the
case for all $\eu \in (\eA/\ef)^{\times}$ that $\eu\es \bmod \ef$ is
monic if and only if $\eu \bmod \ef$ is monic.  Thus $\eF =
\eF(1)$, $\eL = \eL_f$, and $H_f$ is simple.
\end{proof}

The following proposition follows from Proposition~\ref{P:Gmult} and
the simplicity of $H_f$.

\begin{proposition} \label{P:sbrk}
The soliton $t$-module $E_{\ef}$ has a proper sub-$t$-module if and
only if there exist distinct $\ea$, $\eb \in \cI$ such that
$\Gamma(a/f) \approx \Gamma(b/f)$.
\end{proposition}

We now consider special values of the Gamma function at fractions
having different denominators.  When necessary we write for any $n \in
\eA_+$, $\cI_n = \{ a \in \eA : \deg(a) < \deg(n), (a,n) = 1 \}$, and
similarly for $\eF_n(a)$ corresponding to the group in \eqref{E:Fa}.

\begin{theorem} \label{T:BrkIs}
Let $\ef$ and $\eg$ be monic and distinct.  The following are
equivalent.
\begin{enumerate}
\item[(a)] There exist non-trivial $t$-module homomorphisms
$E_{\ef} \to E_{\eg}$.
\item[(b)] $H_{\ef}$ and $H_{\eg}$ are isogenous.
\item[(c)] There exist $\ea \in \cI_{\ef}$ and $\eb \in \cI_{\eg}$ such
that $\Gamma(a/f) \approx \Gamma(b/g)$.
\end{enumerate}
\end{theorem}

\begin{proof}
Certainly (a) and (b) are equivalent according to the discussion in
Section~\ref{SS:CMsub}.  Assuming that $H_{\ef}$ and $H_{\eg}$ are
isogenous, then by Theorem~\ref{T:subt} the CM-field $\eL$ of
$H_{\ef}$ and $H_{\eg}$ is a subfield of $\eK_{\ef} \cap
\eK_{\eg}$ and simultaneously satisfies the criteria of
Lemma~\ref{L:crit} for both $K_{\ef}$ and $K_{\eg}$.  Furthermore, if
we let $\cS_{\eL}$ be the preferred embeddings for the CM-type of
$H_{\ef}$, then for some $\sigma_{\eb} \in \Gal(\eK_{\eg}/\ek)$ we
have
\begin{equation} \label{E:cS}
  \cS_{\eL} = \cS_{\ef}|_{\eL} =
  (\cS_{\eg}|_{\eL})\sigma_{\eb}^{-1}|_{\eL}.
\end{equation}
Now let $\eum$ be the least common multiple of $\ef$ and $\eg$, and
let $\eu \in \cI_{\eum}$.  We claim that
\begin{equation} \label{E:umon}
  \text{$\eu \bmod{\ef}$ is monic $\Longleftrightarrow$
  $\eu\eb \bmod{\eg}$ is monic.}
\end{equation}
Indeed, suppose $\eu \bmod{\ef}$ and $\zeta\eu\eb \bmod{\eg}$ are
monic with $\zeta \in \Fq^{\times}$.  It follows that $\sigma_{\eu}
\in \cS_{\ef}$ and $\sigma_{\zeta\eu\eb} \in \cS_{\eg}$.
However, by \eqref{E:cS} we can choose $\sigma_{\ey} \in
\cS_{\eg}$ so that $\sigma_{\ey}\sigma_{\eb}^{-1}|_{\eL} =
\sigma_{\eu}|_{\eL}$.  Then certainly $\sigma_{\ey}\sigma_{\eb}^{-1}|_{\eL_+} 
= \sigma_{\zeta\eu\eb}\sigma_{\eb}^{-1}|_{\eL_+}$, and so by the
hypotheses on $\eL$ from Lemma~\ref{L:crit}, we must have
\[
\sigma_{\ey}|_{\eL} = \sigma_{\zeta\eu\eb}|_{\eL}.
\]
Therefore $\sigma_{\zeta\eu}|_{\eL} =
\sigma_{\ey}\sigma_{\eb}^{-1}|_{\eL} = \sigma_{\eu}|_{\eL}$, implying
that $\zeta = 1$.  Using \eqref{E:umon} we find that $\Gamma(1/f) \approx
\Gamma(b/g)$, completing (b) implies (c).

Now suppose that $\Gamma(a/f) \approx \Gamma(b/g)$ for some $a \in
\cI_f$ and $b \in \cI_g$.  For any $\eu \in
\cI_{\eum}$ it follows that $\Gamma(au/f) \approx \Gamma(bu/g)$.  Thus
we can specify without loss of generality that $\Gamma(1/f) \approx
\Gamma(b/g)$ for some $\eb \in \cI_{\eg}$.  Let $\eG \subset
\Gal(\eK_{\eum}/\ek)$ be the subset
\[
\eG \assign \{ \sigma_{\es} \in \Gal(\eK_{\eum}/\ek) : \Gamma(s/f)
\approx \Gamma(1/f) \}.
\]
Certainly $\eG$ is a subgroup of $\Gal(\eK_{\eum}/\ek)$ by the bracket
relations.  We claim that
\[
\eG = \{ \sigma_{\es} \in \Gal(\eK_{\eum}/\ek) : \Gamma(bs/g)
\approx \Gamma(b/g) \}. 
\]
Indeed, if $\eu \in (\eA/\eum)^{\times}$ and $\sigma_{\es} \in \eG$,
then
\[
  \xymatrix{ \textnormal{$ub \bmod g$ monic}
  \ar@{<=>}[r]^{\eqref{E:umon}}
  \ar@{<=>}[d]_{\sigma_s|_{K_g} \in \eF_g(1)} &
  \textnormal{$u \bmod f$ monic}
  \ar@{<=>}[d]^{\sigma_s|_{K_f} \in \eF_f(1)} \\
  \textnormal{$ubs \bmod g$ monic}
  \ar@{<=>}[r]_{\eqref{E:umon}} &
  \textnormal{$us \bmod f$ monic.}
  }
\]
Furthermore, as in \eqref{E:Fa},
\[
\eG|_{\eK_{\ef}} = \eF_{\ef}(1) \quad
\text{and} \quad
\eG|_{\eK_{\eg}} = \eF_{\eg}(1).
\]
We let $\eL$ be the fixed field of $\eG$, and thus $\eL$ is the
CM-field of both $H_{\ef}$ and $H_{\eg}$.  Because
$(\eA/\eum)^{\times} \to (\eA/\ef)^{\times}$ and $(\eA/\eum)^{\times}
\to (\eA/\eg)^{\times}$ are surjective, if follows from \eqref{E:Sf}
that
\begin{align*}
  \cS_{\ef} &= \{ \sigma_{\eu}|_{\eK_{\ef}} : \text{$\eu \in
    (\eA/\eum)^{\times}$ and $\eu \bmod \ef$ is monic} \}, \\
  \cS_{\eg} &= \{ \sigma_{\eu}|_{\eK_{\eg}} : \text{$\eu \in
    (\eA/\eum)^{\times}$ and $\eu \bmod \eg$ is monic} \}.
\end{align*}
Because $\Gamma(1/f) \approx \Gamma(b/g)$, it follows that
\[
\sigma_{\eu}|_{\eK_{\ef}} \in \cS_{\ef} \Longleftrightarrow
\sigma_{\eu}|_{\eK_{\eg}} \sigma_{\eb} \in \cS_{\eg}.
\]
Thus $\cS_{\ef}|_{\eL} = \cS_{\eg}|_{\eL}
\sigma_{\eb}^{-1}|_{\eL}$, and by Theorem~\ref{T:subt} $H_f$ is
isogenous to $H_g$.
\end{proof}

\begin{remark}
  The proofs of the above results make no use of the period
  computations for $E_{\ef}$ and $Q_{\ef}$ performed in the previous
  sections.  With these period computations in hand, we now point out
  that Yu's Theorem of the Sub-$t$-module provides another explanation
  for the direction that the existence of bracket relations of the
  form $\Gamma(a/f) \approx \Gamma(b/g)$ guarantees that $H_{\ef}$ and
  $H_{\eg}$ are isogenous.  These methods will be the main focus of
  the next section.
\end{remark}

% indep.tex

\section{Linear Independence Results} \label{S:indep}

\subsection{Yu's Theorem of the Sub-$t$-Module}

The following result is indispensable for all our transcendence
considerations. 

\begin{theorem}{\rm (Yu, \cite{Yu97}, Thm.\ 3.3)}  \label{T:Yusubt} Let the
  $t$-module $E = (\Phi,\Ga^d)$ be defined over $\ok$. Let $\bu \in
  \Lie(E)$ with $\Exp(\bu) \in \ok^d$.  Let $V$ be the smallest
  $\C$-vector space defined over $\ok$ which contains $\bu$ and for
  which
\begin{equation} \label{E:closure}
  d\Phi(t) V \subseteq V.
\end{equation}
Then $V = \Lie(H)$ for some sub-$t$-module $H$ of $E$.
\end{theorem}

This result is a precise analogue of W\"ustholz's Subgroup Theorem
\cite{wue89}, with one very important deviation, namely the condition
\eqref{E:closure}.  Because of \eqref{E:closure}, Yu's theorem
\ref{T:Yusubt} does not directly imply the $\ok$-linear independence
of the coordinates of the point $\bu \in \Lie(E)$ except when the
nilpotent part $N$ of $d\Phi(t) = \theta I_d + N$ is actually zero.
Therefore we show that $N = 0$ for Sinha's soliton $t$-modules.  This
fact and our analysis of the structure of such $t$-modules allows the
complete analysis of the $\ok$-linear independence of $1$, $\tpi$ and
values $\Gamma(a/f)$ where $a, f \in \eA_+$ with $(a,f) = 1$.  Since
the other special Gamma values $\Gamma(r), r \in \ek$ occur as
coordinates of periods of quasi-periodic extensions, and $N = 0$ in
that case as well, we use our analysis of the structure of products of
minimal extensions of $t$-modules to obtain the desired general
results of $\ok$-linear independence by, in a sense, reducing to the
soliton base case.  

\subsection{Proofs of Results}

In this section, we apply the results of preceding sections to obtain
linear independence statements for minimal quasi-periodic extensions
of simple $t$-modules, for $t$-modules of CM-type in general, and for
soliton $t$-modules in particular.  As the setting becomes more and
more specialized, we shall see that the assertions become more
specific. In this section, for a $t$-module $E$ defined over $\ok$, we
denote by $Q_E$ a strictly quasi-periodic extension of maximal
dimension and defined over $\ok$.

\begin{theorem} Let $H$ be a simple $t$-module defined
  over $\ok$ in which $d\Phi_H(t) = \theta I_d$.  Let $Q \assign Q_H =
  (\Psi, \Ga^{d + j})$ be a strictly quasi-periodic extension of $H$
  with corresponding quasi-periodic functions $F_1, \dots, F_j$.  Let
  $\bu = (u_1,\dots,u_d) \in \C^d$ be non-zero with $\Exp_H(\bu) \in
  \ok^d$.  Then the quantities
\begin{equation} \label{E:qplinind}
u_1, \dots, u_d, F_1(\bu), \dots, F_j(\bu)
\end{equation}
are $\ok$-linearly independent.
\end{theorem}

\begin{proof} 
  If the values of \eqref{E:qplinind} are linearly dependent over
  $\ok$, then, by Yu's Theorem of the Sub-$t$-Module \ref{T:Yusubt},
    the point corresponding to \eqref{E:qplinind} lies in the tangent
    space at the origin of a sub-$t$-module $R$ of $Q$.
  
Proposition \ref{P:minqp} shows that $Q$ is a minimal extension of
$H$.  Therefore, by minimality, $\bu$ is contained in the tangent
space of a proper sub-$t$-module $S$ of $H$.  As $H$ is simple, $S$ is
zero, contrary to our choice of a non-zero $\bu$.
\end{proof}

The proof of this result generalizes to admit several $H$ and various
$\bu$.  We say that the points $\bu_1, \dots, \bu_n \in \Lie(E)$ are
{\em linearly independent over $\End(E)$} when the only endomorphisms
$e_1,\dots,e_n \in \End(E)$ with $de_1 \bu_1 + \dots + de_n \bu_n = 0$
are $e_1 = \dots = e_n = 0$.

\begin{corollary}
  Let $H_1,\dots,H_n$ denote non-isogenous simple $t$-modules defined
  over $\ok$ and with each $d\Phi_i(t) = \theta I_{d_i}$.  For each $i
  = 1, \dots, n$, let the set ${\mathcal U}_i = \{ \bu_{i1}, \dots,
  \bu_{i\ell_i} \}$ be linearly independent over $\End(H_i)$,
  $\bu_{i\ell} = (u_{i\ell 1},\dots,u_{i\ell d_i})$ with each
  $\Exp_{H_i}(\bu_{i\ell}) \in \Ga(\ok)^{d_i},$ and let $F_{i1}, \dots,
  F_{is_i}$ defined over $\ok$ be the quasi-periodic functions
  corresponding to a maximal strictly quasi-periodic extension
  $Q_{H_i}$ over $\ok$.  Then the quantities
\[
u_{i\ell m},\  F_{jh}(\bu_{j\ell^\prime})
\]
are $\ok$-linearly independent.
\end{corollary}

\begin{proof}
  By Yu's Theorem of the Sub-$t$-Module \ref{T:Yusubt}, any linear
  dependence relation would give a proper sub-$t$-module $R$ of $\prod
  Q_{H_i}^{\ell_i}$ with $\Lie(R)$ containing the point whose
  coordinates are given by tuples of the above values.  Since strictly
  quasi-periodic extensions are minimal, we know by Proposition 1 of
  \cite{wdb01} that $R$ projects onto a proper sub-$t$-module $S$ of
  $\prod H_i^{\ell_i}$.
  
  In that case, we can apply Yu's Kolchin-type result \cite{Yu97},
  Thm.~1.3, to conclude that, for some fixed $i$, there are non-zero
  endomorphisms $\Theta_{\ell} \in \End(H_i), \ell = 1, \dots, \ell_i$
  such that the projection of $S$ is contained in the sub-$t$-module
  of $H_i^{\ell_i}$ given by $\sum_{\ell} \Theta_{\ell} \bx_{\ell} =
  0$.

In particular, $\sum d\Theta_{\ell} \bu_{i\ell} = 0$.  This identity is
contradicted by our hypothesis that $\bu_{i1},\dots,\bu_{i\ell_i}$ are
linearly independent over $\End(H_i)$.
\end{proof}

Recall that, according to Theorem~\ref{T:Esub}, every $t$-module
$E$ of CM-type is isogenous to a power of a simple sub-$t$-module $H$
(also of CM-type).  So the above result tells us in particular that the
coordinates of periods of $Q_E$ are given by those of $Q_H$.  However
in the CM setting, we can be even more specific.

\begin{theorem} \label{T:qpspan} Let $E$ be a $t$-module of CM-type
  $(\eK,\cS)$ defined over $\ok$.  Say that $E$ is isogenous to $H^m$
  with $H$ simple, defined over $\ok$, and of CM-type
  $(\eL,\cS|_{\eL})$.  Let $Q_E$ and $Q_H$ denote maximal strictly
  quasi-periodic extensions of $E$ and $H$, respectively, defined over
  $\ok$.  Then the $\ok$-vector space $V_E$ spanned by the coordinates
  of all periods of $Q_E$ has as basis the coordinates of {\em any}
  non-zero period of $Q_H$.
\end{theorem}

\begin{proof}  
  Note from Corollaries \ref{C:qpprod}, \ref{C:isogspan} that $V_E$
  is spanned by the coordinates of the periods from $Q_H$.  
  We use the CM structure to show that $V_{Q_H}$ is spanned by the
  coordinates of any non-zero period of $Q_H$:
  
Define $\bss_H \assign \bss_{\cS}|_L$, in the notation of Section
\ref{SS:CM}. By Remark \ref{R:CMisog}, we may assume that the period
lattice of $H$ is $\bss_H(\eB_H)\bsl$ for some $\bsl$.  We need to
verify that, for each $\bsd$ among the $\Phi_H$-biderivations
$\bsd_1,\dots,\bsd_j$ underlying the quasi-periodic extension $Q_H$,
$F_{\bsd}(\bss_H(\eB_H)\bsl)$  lies in the $\ok$-span of the coordinates
of $\bsl$ and the $F_{\bsd_i}(\bsl)$.

Recall that $\Phi_H$ extends from $\eA$ to $\eB_H$ in such a way that
\[
\Phi_H(t)\Phi_H(\eb) = \Phi_H(\eb)\Phi_H(t)
\]
and thus in particular $\bss_H(\eb)d\Phi_H(t) = d\Phi_H(t)\bss_H(\eb)$.
Now let $\bsd$ be a $\Phi_H$-biderivation.  Then the fact that
   \begin{align*}
F_{\bsd}(\bss_H(\eb)d\Phi_H(t) \bz)  
                 &= F_{\bsd}(d\Phi_H(t) \bss_H(\eb) \bz) \\
                 & = \theta F_{\bsd}(\bss_H(\eb) \bz ) + \bsd(t)
                 \Exp_H(\bss_H(\eb)\bz) \\
                 & = \theta F_{\bsd}(\bss_H(\eb) \bz ) + \bsd(t)
                 \Phi_H(\eb)\Exp_H(\bz) 
  \end{align*}
  shows thus that $F_{\bsd}(\bss_H(\eb) \bz)$ is itself the
  quasi-periodic function associated to the $\Phi_H$-biderivation
  $\Phi(b)_\ast \bsd$.  Therefore the values of this function at $\bsl
  \in \Lambda$ will lie in the $\ok$-span of the coordinates of any
  non-zero period.  Hence $V_{Q_H}$ is spanned by the coordinates of
  any non-zero period of $Q_H$.
  
  Since the $t$-action on $\Lie(Q_H)$ is scalar,  Yu's Theorem of
  the Sub-$t$-Module \ref{T:Yusubt}, implies that any $\ok$-linear
  relation on the coordinates of a fixed period $\bsl$ actually will
  hold on $\Lie(R)$ for some proper sub-$t$-module $R$ of $Q_H$.
  However $Q_H$ is a minimal extension of the $t$-module $H$.  So the
  period onto which $\bsl$ projects in $H$ would lie in a proper
  sub-$t$-module of $H$.  Since $H$ is simple, $\bsl$ projects onto
  $0$; however such a $\bsl = 0$, contrary to our hypothesis.
  Therefore the coordinates of $\bsl$ and the quasi-periods
  $F_{\bsd_i}(\bsl)$ form an $L$-basis for $V_H$, as claimed.
\end{proof}

The preceding theorem also extends to several $t$-modules at once.

\begin{theorem} \label{T:CMbasis}  
  Let $E_1, \dots, E_n$ be $t$-modules of CM-type defined over $\ok$.
  For each $i$, we use the following notation:
\begin{enumerate}
\item[(a)] $E_i \thicksim H_i^{m_i}$, with $H_i$ simple, defined over
  $\ok$.
\item[(b)] $P_i$ denotes the set of coordinates of a non-zero period
  of $Q_{H_i}$.
\item[(c)] $V_i$ denotes the $\ok$-vector space spanned by all the
  non-zero coordinates of periods of $Q_{E_i}$.
 \end{enumerate}
 If the $H_i$ are pair-wise non-isogenous, then $\cup_{i = 1}^n P_i$
 is a $\ok$-basis of $V_1 + \dots + V_n$.
\end{theorem}

\begin{proof}
  We know from the preceding result that $P_i$ is a $\ok$-basis for
  the $\ok$-vector space $V_i, i = 1, \dots, n.$ Since each $Q_{H_i}$
  is a minimal extension of $H_i$, then by Lemma 1 of \cite{wdb01},
  $\prod Q_{H_i}$ is a minimal extension of $\prod H_i$.  We know by
  Yu's Theorem of the Sub-$t$-Module \ref{T:Yusubt} that any
  $\ok$-linear relation on $\cup P_i$ gives rise to a proper
  sub-$t$-module of $\prod Q_{H_i}$, which, by minimality, projects
  onto a proper sub-$t$-module $R$ of $\prod H_i$.  By the simplicity
  and non-isogeneity of the $H_i$, we know that the only proper
  sub-$t$-modules for $\prod H_i$ have trivial projections onto some
  factor, say $H_j$.  But then the underlying period of $H_j$ must be
  zero.  That forces the whole period of $Q_{H_j}$ to vanish, contrary
  to our choice of $P_j$.
\end{proof}

The above considerations apply to soliton $t$-modules.  But we can be
completely precise in this case.  For the following theorem we recall
that the notation $\Gamma(r_1) \approx \Gamma(r_2)$ means that a
bracket relation implies that $\Gamma(r_1)/\Gamma(r_2) \in \ok$.

\begin{theorem}
Let $q > 2$.  The numbers
\[
\{ 1, \tpi \} \cup \text{\rm rep}\{ \Gamma(r) : r \in k \smallsetminus
A\}/\approx 
\]
are $\ok$-linearly independent, where the notation means that we take
any set of representatives of the equivalence classes for the relation
$\approx$ on the set $\{ \Gamma(r) : r \in k \smallsetminus A \}$.
\end{theorem}

\begin{proof}
It is clear that the claim holds for every choice of representatives
if it holds for any choice. 
  If $r = b + a/f,\ a,b \in A, f \in A_+$,  then $\Gamma(r) \approx
  \Gamma(a/f)$, since $a \bmod f = (bf + a) \bmod f$.  Therefore we
  may always choose representatives of $\approx$ of the form
  $\Gamma(a/f),\ a \in \cI_f, f \in A_+$.

  Proposition \ref{P:Gmult} says that all $\approx$-equivalence
  classes among these values have cardinality $m$ when $E_f \thicksim
  H_f^m$.  Thus there are exactly $\dim Q_{H_f}$ classes among them.
  According to equation \eqref{E:piasim}, Theorem \ref{T:qp}, and
  Theorem \ref{T:qpspan}, these representatives span $V_{E_f}$, a
  $\ok$-space of dimension $\dim Q_{H_f}$.  Therefore any choice
  rep$(Q_{H_f})$ of $\approx$-representatives among the $\Gamma(a/f)$
  gives a $\ok$-basis for $V_{E_f} = V_{H_f}$.
  
  Thus one choice of $\text{\rm rep}\{ \Gamma(r) : r \in k
  \smallsetminus A\}/\approx $ will be a disjoint union of
  rep$(Q_{H_f})$, taken over non-isogenous $H_f$.  Consequently, the
  claim of the theorem is equivalent to the statement that $1, \tpi$
  and coordinates of non-zero periods from non-isogenous $Q_{H_f}$ are
  $\ok$-linearly independent.  We proceed to prove this assertion.
  
  Yu's Theorem of the Sub-$t$-Module \ref{T:Yusubt}, shows
  that, if there were a $\ok$-linear relation on the coordinates of
  such periods, then there would have to be a proper sub-$t$-module
  $R$ of a finite product of the form
\[
Q \assign \Ga \times C \times Q_{H_{f_1}} \times \dots \times Q_{H_{f_n}}
\]
for which $\Lie(R)$ would contain the point ${\mathbf q} \assign (1,
\tpi, \bsl_{f_1}, \dots, \bsl_{f_n})$, where each $\bsl_{f}$ is some
non-zero period of $Q_{H_f}$, $ f = f_1, \dots, f_n$ and $C$ denotes
the Carlitz module.  Here the $Q_{H_f}$ occur exactly when some
$\Gamma(a/f), a \in \cI_f$ is involved in the supposed $\ok$-linear
relation.  Moreover, since $q > 2$, $\Ga$ and/or $C$ appear exactly
when $1$ and/or $\tpi$ are involved in the relation.  Cf.\ 
\cite{Th91}, \S 6.  For ease of exposition, we simply assume that to
be the case here.

By Lemma 1 of \cite{wdb01}, the above product is minimal.  Therefore
$R$ projects onto a proper sub-$t$-module $S$ of the corresponding
product
\[
E \assign \Ga \times C \times {H_{f_1}} \times \dots \times {H_{f_n}}
\]
for which $\Lie(S)$ would contain the projection ${\mathbf p}$ of the
point of ${\mathbf q}$ of $\Lie(R)$.  Again the point ${\mathbf p}$ in
$\Lie(S)$ projects non-trivially onto the factors because, as we saw
in Proposition \ref{P:qpext}, the non-zero periods of $Q_H$ are
produced from -- and project to -- non-zero periods of $H$.

For the conclusion of the proof we keep in mind the following three
remarks:
\begin{enumerate}
\item[(a)] The underlying simple $H_f$, $H_g$ have been taken to be
non-isogenous.

\item[(b)] The Carlitz module $C$ has period an algebraic multiple of
  $\tpi$, and $C$ does not have CM.  Therefore it is a simple
  $t$-module which is not isogenous to any of the soliton $t$-modules
  $E_f$.
  
\item[(c)] The trivial $t$-module $\Ga$ is also simple and not isogenous to
  $C$ nor any $H_f$.
\end{enumerate}
Now 
\[
\Lie(E) = \Lie(\Ga) \times \Lie(C) \times \Lie(H_{f_1}) \times \dots
\times \Lie(H_{f_n}). 
\]
Since the factors of $E$ are non-isogenous and simple, the only proper
sub-$t$-modules have tangent spaces which project trivially to the
tangent space of at least one factor.  Thus $\Lie(S)$ does so, as $S$
is a proper sub-$t$-module.  However this is contradicted by the facts
that ${\mathbf p} \in \Lie(S)$ and our choice of product $Q$ ensures
that ${\mathbf p} \in \Lie(S)$ has a non-zero entry in every factor of
$\Lie(E)$.  We conclude that there is no non-trivial $\ok$-linear
relation on the set in question.
\end{proof}

\begin{remark}  We note that, in the results of this section,  we can
  adjoin coordinates of linearly independent logarithms of algebraic
  points, e.g.\ periods, of other simple $t$-modules defined over
  $\ok$, as long as the $t$-modules are not isogenous to each other
  nor to the $H_f$ nor the Carlitz module.
\end{remark}

\subsection{Proof of Corollary \ref{C:moredenoms}}
As Corollary \ref{C:nonintegral} is immediate from the Main Theorem,
we need only consider Corollary \ref{C:moredenoms}.  We break the
proof up into several cases: we may without loss of generality
consider $\ef$ and the $\ef_i$ to be monic.

\subsubsection{Case $f = f_1$}
Let $\eF \subset (\eA/f)^\times$ correspond to the Galois group
$\Gal(\eK_f/\eL)$.  Assume that condition (b) of Lemma \ref{L:crit}
holds for $\eF$, where $\cS_f \simeq \cI_+$.  Then for any $b \in
\eF$,
\[
b \cI_+ = \cI_+,
\]
as $\cI_+$ is a union of cosets of $\eF$.  Now 
\[
S \assign \sum_{a \in \cI_+} \ea = 1,
\]
as the sum of all monic polynomials over $\Fq$ of fixed positive
degree with all but the constant term fixed is easily seen to vanish.
Thus
\[
b = b S = S = 1,
\]
and $\eF$ is trivial.  We conclude from part (b) of
Theorem~\ref{T:Esub} that $E_f$ is simple.  See
Shimura~\cite{ShimCM}, p.~64, for the analogue of this case for abelian
varieties with CM by $\mathbb{Q}(\zeta_p)$.  As $E_f$ is simple and
$f$ is irreducible, Proposition \ref{P:sbrk} and Theorem \ref{T:MT}
show that all the quantities $\Gamma(a/f)$, $\deg a < \deg f$ are
$\ok$-linearly independent.

\subsubsection{Case $f = f_1^{e_1}$, $e_1 \ge 2$}  
The proof begins as before, except that now we must exclude the monic
multiples of $f_1$ of degree less than $\deg f$ from the sum $S$.  The
argument given in the first case shows that the multipliers involved
in these multiples add up to $1$ and therefore the multiples
themselves sum to $f_1$.  Hence
\[
S  \assign \sum_{\ea \in \cI_+} \ea = 1 - \ef_1.
\]
As before, any $\eb \in \eF$ satisfies
\[
\eb S \equiv S \mod \ef,
\]
and, since $(S,\ef) = 1$, we conclude $\eb = 1$.  Thus as in the
previous case, $E_{\ef}$ is simple, and the $\Gamma(a/f)$, $\deg a <
\deg f$, are $\ok$-linearly independent.

However this does not cover the values $\Gamma(a/f)$ where $f_1
\mid a$.  For them consider the $t$-modules $E_{f_1}, E_{f_1^2},
\dots, E_{f_1^{e_1-1}}$, which are simple and according to their
dimensions are non-isogenous.  Thus from Theorems
\ref{T:BrkIs} and \ref{T:MT} we obtain the linear independence of all
$\Gamma(a/f)$ for $0 \le \deg a < \deg f$, whether $(a,f) = 1 $ or
not.

\subsubsection{General case $\ef = \ef_1^{e_1}\dots \ef_m^{e_m}$.}

Inclusion-exclusion gives that
\[
S \assign \sum_{a \in I_+} a \equiv (1 - \ef_1) \cdots (1 - \ef_m) \mod \ef
\]
(with equality unless $\ef = \ef_1 \cdots \ef_m$).  Thus by hypothesis
on the $\ef_i$, $(S,\ef) = 1$, whereas $\eb S = S$.  So $\eb = 1$ and $E_{\ef}$
is simple.  This accounts for the linear independence of all
$\Gamma(a/f)$, $a \in \cI_{\ef}$.

The values $\Gamma(a/f)$ with $(a,f) \ne 1$ are included by induction
through the remark that, for non-constant proper monic divisors $\eg \mid
\ef$, the various $E_{\eg}$ are simple and, since they have distinct
CM-fields, are non-isogenous to each other and to $E_{\ef}$.  The
corollary follows.

% examples.tex

\section{Examples} \label{S:exmp}

Many of the calculations below on Gamma values are due to
Sinha~\cite{SS95}, {\S}VI.3.2, and Thakur~\cite{Th91}, {\S}9.  The
reader interested in the explicit computation of the Anderson-Coleman
soliton functions should consult the specific examples of
Coleman~\cite{rc88} and Sinha~\cite{SS97d}, {\S}3.3, and the general
methods of Thakur~\cite{Th99}.  The examples below demonstrate in
particular the correspondence between $\ok$-linear relations on Gamma
values and the structure of the underlying soliton $t$-modules
developed in this paper.

\subsection{The Simplest Case} \label{SS:simp}
Let $\ef = t$.  In this case $\cI_+ = \{ 1 \}$, and so the $t$-module
$E_{\ef}$ is simply a Drinfeld module as it has dimension $1$.

Here $\eK_{\ef} = \ek(z)$ and $\eB_{\ef} = \ek[z]$, where $z \assign
\zeta_t = \sqrt[q-1]{-t}$, and thus the curve $\bX$ is isomorphic to
$\PP^1/\C$.  Note that the point $\xi \in \bX(\C)$ corresponds to the
zero of $z - \zeta_\theta$.  Now the soliton function is
\[
g_{\ef} = 1 - \frac{z}{\zeta_{\theta}},
\]
and by \eqref{E:WXiI} we have $W = 0$, $\Xi = \xi$ and $I =
\infty$.  The $t$-motive $M_{\ef}$ is then
\[
M_{\ef} = \Gamma(\bU,\Omega_{\bX}) = \C[z]dz.
\]
As defined in Section~\ref{SS:per} we let $n_1 =
\zeta_\theta^{q-2}\,dz$ and then note by \eqref{E:tau} that $zn_1 =
\left( \zeta_\theta \tau^0 +
\frac{\zeta_\theta}{\theta^{q - 2}}\tau \right) n_1$.  Hence
\begin{equation} \label{E:Phiz}
  \Phi(z) = \zeta_\theta \tau^0 + \frac{\zeta_\theta}{\theta^{q-2}} \tau.
\end{equation}
Thus $E_{\ef}$ is isomorphic over $\ok$ to the Carlitz module for the
polynomial ring $\Fq[z]$.  By Theorem~\ref{T:Spia} we see that $\pi_1
= \Gamma(1/\theta)$.  On the other hand, if $\tpi_\theta$ is the
period of the Carlitz module for $\Fq[z]$, then we find from
\eqref{E:Phiz} that
\begin{equation} \label{E:pi1pi}
\pi_1 = \Gmm{\frac{1}{\theta}} = \frac{\theta}{\sqrt[q - 1]{\theta
\zeta_\theta}} \, \tpi_\theta.
\end{equation}
When $q=2$, we note that $\tpi_\theta = \tpi$ and in fact $\pi_1
\thicksim \tpi$.

We now compute the quasi-periods for $E_{\ef}$.  For each $\ell \in
\Fq^{\times}$, $\ell \neq 1$, we choose $n_{\ell}$ as in
Lemma~\ref{L:na}.  Namely we let
\[
n_{\ell} \assign \frac{\ell}{\zeta_{\theta}^{q-2}} \left(
  \frac{t - \theta}{z - \ell \zeta_{\theta}} \right)^2 dz.
\]
By Theorem~\ref{T:qp} and in particular \eqref{E:efin}, we see that
$\eta_{\ell} = - \res_{[\ell]\xi}(c_g n_{\ell}/(t-\theta))$.
Proceeding as in \eqref{E:pia}, we obtain the quasi-periods
\[
\eta_{\ell} = \frac{1}{\ell - 1} \cdot \frac{\ell}{\theta} \,
\Gmm{\frac{\ell}{\theta}}, \quad \ell \in \Fq^{\times}, \ell \neq 1.
\]
Note that Corollary~\ref{C:powerdenom} shows that the numbers
$\Gamma(\ell/\theta)$, $\ell \in \Fq^{\times}$, are $\ok$-linearly
independent.

\subsection{A Non-Uniform Example: $\ef = t(t-1)$}

We saw in the proof of Corollary \ref{C:moredenoms} that, in order
that $E_{\ef}$ be non-simple when $\ef = \ef_1^{e_1}\dots
\ef_m^{e_m}$, we must have some $\ef_i | (1 - \ef_j)$.  In this
example, we will see that this necessary condition is not uniformly
sufficient even in the simplest case, namely $\ef = t(t-1)$.

We look for possible subgroups $\eF$ of $\cI_+$.  Since
$\ef$ is quadratic, we consider monic linear polynomials $\eb = t + \ell
\in \eF$, $\ell \in \Fq$.  Then 
\[
\eb^2 \in \eF \Leftrightarrow 
\begin{cases}
1 + 2\ell = 0 & \text{if $\ell^2 = 1$,} \\
2\ell = 0 & \text{otherwise.}
\end{cases}
\]

In the first case, $4 = 1$, so $p = 3$.  Since $(\eb,\ef) = 1$, $a \ne
0$, so in the second case $p = 2$.  Thus, for $p > 3$, the $t$-module
$E_{\ef}$ is simple.

Let us examine the first case more closely when $p = 3$ under the
assumption that $E_{\ef}$ is non-simple.  Then $\ell = 1$, so $\eF =
\{1, t+1\}.$ Now as
\begin{equation} \label{E:F_3}
  \cI_+ = \{ 1 \} \cup \{t + a : a \in \Fq, a \ne 0, 1\} 
\end{equation}
is a union of cosets of $\eF$, we know that $(t+1) \cI_+ = \cI_+$.  In
particular for any element $t + a \neq t + 1$ of $\cI_+$, $(t+a)(t+1)
\in \cI_+ \setminus \{ 1 \}$.  This means that $2 + a = 1$, and so $t
+ a = t - 1$.  But $(\ef, t-1) \neq 1$.  So in fact, if $E_{\ef}$ is
non-simple and $p = 3$, then $\cI_+ = \eF = \{1, t+1\}$, and obviously
{from} \eqref{E:F_3}, $q = 3$.  Thus by Theorem~\ref{T:Esub},
$E_{\ef}$ is isogenous to a power of a Drinfeld module.  The lattice
for this Drinfeld module can be taken to be the ring of integers in
the fixed field $\ek(\zeta_t)$ of $\eF$, which in this case is
${\mathbb F}_3[\zeta_t] = {\mathbb F}_3[\sqrt{-t}]$.

When $p = 2$, the $t$-module $E_{\ef}$ is $1$-dimensional if $q=2$.
For $q > 2$, if $t+a$, $t+b$ are elements of $\eF$, then the closure
of $\ef$ under multiplication requires that the product $(t+a)(t+b)$
be monic modulo $\ef$ of degree either one or zero.  In the first
case, $a = b$; in the second
\[
1 + a + b = 0, \quad ab = 1,
\]
i.e. $a^2 + a + 1 = 0$, i.e. $a \in {\mathbb F}_4 \setminus {\mathbb
F}_2$.  Thus $q=4$, and $E_{\ef}$ is isogenous to a power of the
Drinfeld module whose lattice is ${\mathbb F}_4[t,\sqrt[3]{t(t+1)}]$,
the ring of integers in the fixed field $\ek(\zeta_t\zeta_{t+1})$ of
$\eF$.  In this way we obtain the following result:

\begin{corollary}
The $t$-module $E_{t(t-1)}$ is simple except in the following two
cases:

\begin{enumerate}
\item[(a)] $E_{t(t-1)}$ is isogenous to $H_3^2$ when $q = 3$, where $H_3$
  is the Drinfeld $\mathbb{F}_3[t]$-module with lattice ${\mathbb
  F}_3[\sqrt{-t}]$;
\item[(b)] $E_{t(t+1)}$ is isogenous to $H_4^3$ when $q = 4$, where $H_4$
  is the Drinfeld $\mathbb{F}_4[t]$-module with lattice ${\mathbb
  F}_4[t,\sqrt[3]{t(t+1)}]$.
\end{enumerate}
\end{corollary}

We consider the exceptional cases $q=3$, $4$ in a bit more detail:

\emph{Case $q=3$.}  From Proposition~\ref{P:sbrk} we see that
\[
\Gmm{\frac{1}{\theta(\theta - 1)}} \thicksim \Gmm{\frac{\theta +
1}{\theta(\theta - 1)}}.
\]
Now we see that the sub-$t$-modules of $E_{\ef}$ are isogenous to
$E_t$.  Moreover, taking into account \eqref{E:pi1pi} we see
\[
\Gmm{\frac{1}{\theta}} \thicksim \Gmm{\frac{1}{\theta(\theta - 1)}} \thicksim
\Gmm{\frac{\theta + 1}{\theta(\theta - 1)}} \thicksim \tpi_{\theta}.
\]
All of these equivalences are confirmed by the bracket relations.  In
fact, if we let $\Exp_t$ be the exponential function of $E_t$, then it
can be shown that the exponential function $\Exp_{\ef}$ of $E_{\ef}$ is
\[
\Exp_{\ef} \begin{pmatrix} z_1 \\ z_2 \end{pmatrix}
= \begin{pmatrix}
  \alpha & \zeta_{\theta-1} \alpha \\
  \beta & -\zeta_{\theta-1} \beta
        \end{pmatrix}
\begin{pmatrix}
\Exp_t(-z_1/\alpha - z_2/\beta) \\
\Exp_t(-z_1/\zeta_{\theta-1}\alpha + z_2/\zeta_{\theta-1}\beta)
\end{pmatrix},
\]
where $\Gamma(1/\theta(\theta-1)) = \alpha
\theta(\theta-1)\Gamma(1/\theta)$ and
$\Gamma((\theta+1)/\theta(\theta-1)) = \beta \Gamma(1/\theta)$
(cf.~proof of Theorem~\ref{T:Esub}).

\emph{Case $q=4$.}
Here the monics in $\cI_+$ are $1$, $t + \ell$, $t + \ell^2$.  The
bracket relations as well as Proposition~\ref{P:sbrk} confirm that
\[
\Gmm{\frac{1}{\theta(\theta+1)}} \thicksim
\Gmm{\frac{\theta + \ell}{\theta(\theta+1)}} \thicksim
\Gmm{\frac{\theta + \ell^2}{\theta(\theta+1)}} \thicksim
\tpi_{H_4},
\]
where $\tpi_{H_4}$ is the period of $H_4$.

\subsection{Another Example} $\ef = t(t-1)(t+1)$ and
$q=3$.  In this example, the set of monic elements $\cI_+$ of
$(\eA/f)^{\times}$ consists of four elements:
\[
\cI_+ = \{ 1, t^2+1, t^2 + t - 1, t^2 - t - 1 \}.
\]
We check that $\cI_+$ is a subgroup of $(\eA/f)^{\times}$ and that the
fixed field of $\cI_+$ is the field
\[
\eL \assign \ek(\zeta_t \zeta_{t-1}\zeta_{t+1})
= \ek(\sqrt{-(t^3 - t)}).
\]
The ring of integers in $\eL$ forms a rank $2$ lattice in $\C$ which
corresponds to a Drinfeld $\mathbb{F}_3[t]$-module $\psi$ with CM by
the ring of integers of $\eL$.  The $t$-module $E_{\ef}$ is then
isogenous to $\psi^4$.  If we let $\tpi_{\psi}$ be a fundamental
period of the Drinfeld module $\psi$, then
\[
\Gmm{\frac{1}{\theta^3 - \theta}} \thicksim
\Gmm{\frac{\theta^2 + 1}{\theta^3 - \theta}} \thicksim
\Gmm{\frac{\theta^2 + \theta - 1}{\theta^3 - \theta}} \thicksim
\Gmm{\frac{\theta^2 - \theta - 1}{\theta^3 - \theta}} \thicksim
\tpi_{\psi}.
\]

% gammabib.tex

\end{document}